\documentclass[aos,preprint]{imsart}
\usepackage{amsthm,amsmath}
\usepackage{amsfonts}
\usepackage{amssymb}
\usepackage{graphicx}
\usepackage{epstopdf}
\usepackage{hyperref}
\usepackage{epsfig}
\RequirePackage[colorlinks,citecolor=blue,urlcolor=blue]{hyperref}

\startlocaldefs
\endlocaldefs

\begin{document}

\begin{frontmatter}

\title{Profile Estimation for Partial Functional Partially Linear Single-Index Model\thanksref{T1}}
\runtitle{Functional Single Index Model}
\thankstext{T1}{Part of the data used in preparation of this article were obtained from
the Alzheimer's Disease Neuroimaging Initiative (ADNI) database
(adni.loni.ucla.edu). As such, the investigators within the ADNI contributed
to the design and implementation of ADNI and/or provided data but did not
participate in analysis or writing of this report. A complete listing of ADNI
investigators can be found at: \url{http://adni.loni.ucla.edu/wp-content/
uploads/how\_to\_} \url{apply/ADNI\_Acknowledgement\_List.pdf}}

\begin{aug}
\author{\fnms{Qingguo} \snm{Tang}\thanksref{t1,u1}\ead[label=e1]{tangqguo@163.com}},
\author{\fnms{Linglong} \snm{Kong}\thanksref{t2,u2}\ead[label=e2]{lkong@ualberta.ca}},
\author{\fnms{David} \snm{Ruppert}\thanksref{t3,u3}\ead[label=e3]{dr24@cornell.edu}},
\and 
\author{\fnms{Rohana J.} \snm{Karunamuni}\corref{}\thanksref{t4,u2}\ead[label=e4]{R.J.Karunamuni@ualberta.ca}}

\thankstext{t1}{Supported by the National Social Science Foundation of China (16BTJ019) and Natural
Science Foundation of Jiangsu Province of China (Grant No. BK20151481).}
\thankstext{t2}{Supported by the Natural Sciences and
Engineering Research Council of Canada (NSERC) and the Canadian Statistical
Sciences Institute Collaborative Research Team (CANSSI-CRT).}
\thankstext{t3}{Supported by NSF grant
AST-1312903 and NIH grants P30 AG010129, K01 AG030514.}
\thankstext{t4}{Supported by the Natural Sciences and
Engineering Research Council of Canada (NSERC).}
\runauthor{Q. Tang et al.}

\affiliation{Nanjing University of Science and Technology\thanksmark{u1}, University of Alberta\thanksmark{u2} and Cornell University\thanksmark{u3}}

\address{School of Economics and Management\\
Nanjing University of Science and Technology\\
Nanjing, Jiangsu 210094\\
 China\\
\printead{e1}
}

\address{Department of Mathematical \\and Statistical Sciences\\
University of Alberta\\
Edmonton, AB T6G 2G1\\
Canada\\
\printead{e2}\\
\phantom{E-mail:\ }\printead*{e4}}

\address{School of Operations Research and Information Engineering\\
Cornell University\\
Ithaca, NY 14853\\
 USA \\
\printead{e3}}

\end{aug}

\begin{abstract}
This paper studies a \textit{partial functional partially
linear single-index model} that consists of a functional linear component as
well as a linear single-index component. This model generalizes many
well-known existing models and is suitable for more complicated data
structures. However, its estimation inherits the difficulties and complexities
from both components and makes it a challenging problem, which calls for new
methodology. We propose a novel profile B-spline method to estimate the
parameters by approximating the unknown nonparametric link function in the
single-index component part with B-spline, while the linear slope function in
the functional component part is estimated by the functional principal
component basis. The consistency and asymptotic normality of the parametric
estimators are derived, and the global convergence of the proposed estimator
of the linear slope function is also established. More excitingly, the latter
convergence is optimal in the minimax sense. A two-stage procedure is
implemented to estimate the nonparametric link function, and the resulting
estimator possesses the optimal global rate of convergence. Furthermore, the
convergence rate of the mean squared prediction error for a predictor is also
obtained. Empirical properties of the proposed procedures are studied through
Monte Carlo simulations. A real data example is also analyzed to illustrate
the power and flexibility of the proposed methodology.
\end{abstract}



\end{frontmatter}

\section{Introduction}

Functional data analysis has generated increasing interest in recent years in
many areas, including biology, chemometrics, econometrics, geophysics, medical
sciences, meteorology, etc. Functional data are made up of repeated
measurements taken as curves, surfaces or other objects varying over a
continuum, such as the time and space. In many experiments, such as clinical
diagnosis of neurological diseases from the brain imaging data, functional
data appear as the basic unit of observations. As a natural extension of the
multivariate data analysis, functional data analysis provides valuable
insights into these experiments, taking into account the underlying smoothness
of high-dimensional covariates and provides new approaches for solving
inference problems. One may refer to the monographs of Ramsay and Silverman
\cite{r26,r27}, Ferraty and Vieu \cite{r10} 
and Horv\'{a}th and Kokoszka \cite{r12} 
for
a general overview on functional data analysis.

Motivated by more complicated data structures, which appeal to more
comprehensive, flexible and adaptable models, in this paper we investigate the
following \textit{partial functional partially linear single-index model}:
\[
Y=\int_{\mathcal{T}}a(t)X(t)dt+W^{T}\pmb{\alpha}_{0}+g(Z^{T}\pmb{\beta}_{0}%
)+\varepsilon,\tag{1.1}
\]
where $X(t)$ is a random function defined on some bounded interval
$\mathcal{T}$, $a(t)$ is an unknown square integrable slope function on
$\mathcal{T}$, $W$ is a $q\times1$ vector of covariates, $\pmb{\alpha}_{0}$ is
a $q\times1$ unknown coefficient vectors, $Z\in R^{d}$ is a $d\times1$ vector
of covariates, $\pmb{\beta}_{0}$ is a $d\times1$ coefficient vector to be
estimated, $g$ is an unknown link function and $\varepsilon$ is a random error
with mean zero that is independent of the covariates $(X(t),W,Z)$.

Model (1.1) is more flexible and can deal with more complicated data
structures. To the best of our knowledge, this model has not been fully
studied in the literature yet. It consists of a functional linear component as
well as a linear single-index component. This model generalizes many
well-known existing models and is suitable for more complicated data
structures. However, its estimation inherits some difficulties and
complexities from both components and makes it a challenging problem, which
calls for new methodology. We propose a novel profile B-spline method to
estimate the parameters by approximating the unknown nonparametric link
function in the single-index component part with B-spline, while the linear
slope function in the functional component part is estimated by the functional
principal component basis.

More specially, model (1.1) can be interpreted from two perspectives. First,
it generalizes the partial functional linear models
\[
Y=\int_{\mathcal{T}}a(t)X(t)dt+W^{T}\pmb{\alpha}_{0}+\varepsilon,\tag{1.2}
\]
by adding a nonparametric component, $g(Z^{T}\pmb{\beta}_{0}),$ with an
unknown univariate link function $g.$ This single-index term reduces the
dimensionality from the multivariate predictors to a univariate index
$Z^{T}\pmb{\beta}_{0}$ and avoids the curse of dimensionality, while still
capturing important features in high-dimensional data. Furthermore, since a
nonlinear link function $g$ is applied to the index $Z^{T}\pmb{\beta}_{0},$
interactions between the covariates $Z$ can be modeled. The standard
functional linear model \cite{r5,r3,r11} 
with scalar response $Y$ has the same form as model (1.2)
without the linear part. In general, $X(t)$ can be a multivariate functional
variable, but here we shall only focus on the univariate case. The main
interest is estimation of functional coefficient $a(t)$ based on a sample
$(X_{1},Y_{1}),...,(X_{n},Y_{n})$ generated from model (1.2). There are number
of articles in the literature discussing the slope estimation in model (1.2)
using methods such as the penalized spline method \cite{r5},
the
functional principal component analysis \cite{r41,r3,r11,r44} 
and the functional partial least
squares method \cite{r9},
 among others.

Second, model (1.1) can be considered as a generalization of the partially
linear single-index model \cite{r4,r43},
\[
Y=g(Z^{T}\pmb{\beta}_{0})+W^{T}\pmb{\alpha}_{0}+\varepsilon,\tag{1.3}
\]
with an addition of functional covariates $X(t).$ The partially linear
single-index model (1.3) was first explored by Carroll et al. \cite{r4}.
In fact,
the autors considered a more generalized version, where a known
link function is employed in the regression function, while model (1.3)
assumes an identity link function. Model (1.3) has also been studied by many
other authors, including Xia et al. \cite{r40}
Xia and H\"{a}rdle \cite{r39},
Liang
et al. \cite{r18}
and Wang et al. \cite{r37}
to name a few.

To tackle the challenging estimation problem, our innovation is to propose a
\textit{profile B-spline} (PBS) method to estimate the unknown parameters
$(\pmb{\alpha}_{0}^{T},\pmb{\beta}_{0}^{T})^{T}$ by employing a B-spline
function to approximate the unknown link function $g$ and using the functional
principal component analysis (FPCA) to estimate the slope function $a(t)$.
Under some regularity conditions, we prove the consistency and asymptotic
normality of the proposed estimators. We also establish a global rate of
convergence of the estimator of $a(t)$, and it is shown to be optimal in the
minimax sense of Hall and Horowitz \cite{r11}.
Based on the estimators of
parameters, another B-spline function is employed to approximate the function
$g$ and then the optimal global convergence rate of the approximation is
established. We also obtain convergence rates of the mean squared prediction
error for a predictor. For model (1.3), Yu and Ruppert \cite{r43} 
studied
asymptotic properties of their estimators of $(\pmb{\alpha}_{0}^{T}%
,\pmb{\beta}_{0}^{T})^{T}$ under the condition that the link function $g$
falls in a finite-dimensional spline space. We note here that the asymptotic
properties of all our estimators are derived under the assumption that $g$ can
be well approximated by spline functions with increasing the number of knots.

To gain more flexibility and partly motivated by applications, a number of
other models based on the standard functional linear model have been studied
in the literature, including the partial functional linear regression model
(1.2) \cite{r29,r30,r34},
a generalized
functional linear model \cite{r22,r7}, 
single and multiple index functional regression models \cite{r6,r21}
and a functional partial linear single-index model \cite{r36}
among others.

The paper is organized as follows. Section 2 describes the proposed profile
estimation method. Section 3 presents asymptotic results of our estimator. In
Section 4, we conduct simulation studies to examine the finite sample
performance of the proposed procedures. In Section 5, the proposed method is
illustrated by analyzing a diffusion tensor imaging (DTI) data set from the
Alzheimer's Disease Neuroimaging Initiative (ADNI) database
(adni.loni.ucla.edu). Finally, Section 6 contains some concluding remarks. All
proofs are relegated to the Appendix.

\section{Profile B-spline estimation}

Let $Y$ be a real-valued response variable and $\{X(t):t\in\mathcal{T}\}$ be a
mean zero second-order (i.e., $EX(t)^{2}<\infty$ for all $t\in\mathcal{T})$
stochastic process with sample paths in $L_{2}(\mathcal{T})$, the set of all
square integrable functions on $\mathcal{T}$, where $\mathcal{T}$ is a bounded
closed interval. Let $\langle\cdot,\cdot\rangle$ and $\Vert\cdot\Vert$ denote
the $L_{2}(\mathcal{T})$ inner product and norm, respectively. Denote the
covariance function of the process $X(t)$ by $K(s,t)=\mathrm{cov}(X(s),X(t))$.
We suppose that $K(s,t)$ is positive definite. Then $K(s,t)$ admits a spectral
decomposition in terms of strictly positive eigenvalues $\lambda_{j}$:
\[
K(s,t)=\sum_{j=1}^{\infty}\lambda_{j}\phi_{j}(s)\phi_{j}(t),\ \ \ s,t\in
\mathcal{T},\tag{2.1}
\]
where $\lambda_{j}$ and $\phi_{j}$ are eigenvalue and eigenfunction pairs of
the linear operator with kernel $K$, the eigenvalues are ordered so that
$\lambda_{1}\geq\lambda_{2}\geq\cdots>0$ and eigenfunctions $\phi_{1},\phi
_{2},\cdots$ form an orthonormal basis for $L_{2}(\mathcal{T})$. This leads to
the Karhunen-Lo\'{e}ve representation $X(t)=\sum_{j=1}^{\infty}\xi_{j}\phi
_{j}(t),$ where $\xi_{j}=\int_{\mathcal{T}}X(t)\phi_{j}(t)dt$ are uncorrelated
random variables with mean zero and variance $E\xi_{j}^{2}=\lambda_{j}$. Let
$a(t)=\sum_{j=1}^{\infty}a_{j}\phi_{j}(t)$. Then model (1.1) can be written
as
\[
Y=\sum_{j=1}^{\infty}a_{j}\xi_{j}+W^{T}\pmb{\alpha}_{0}+g(Z^{T}\pmb{\beta}_{0}%
)+\varepsilon.\tag{2.2}
\]
By (2.2), we have
\[
a_{j}=E\{[Y-(W^{T}\pmb{\alpha}_{0}+g(Z^{T}\pmb{\beta}_{0}))]\xi_{j}%
\}/\lambda_{j}.\tag{2.3}
\]
Let $(X_{i}(t),W_{i},Z_{i},Y_{i}),i=1,\cdots,n$, be independent realizations
generated from model (1.1). Then the empirical versions of
$K$ and of its spectral decomposition are
\[
\hat{K}(s,t)=\frac{1}{n}\sum_{i=1}^{n}X_{i}(s)X_{i}(t)=\sum_{j=1}^{\infty}%
\hat{\lambda}_{j}\hat{\phi}_{j}(s)\hat{\phi}_{j}(t),\ \ \ s,t\in
\mathcal{T}.\tag{2.4}
\]
Analogously to the case of $K$, $(\hat{\lambda}_{j},\hat{\phi}_{j})$ are
(eigenvalue, eigenfunction) pairs for the linear operator with kernel $\hat
{K}$, ordered such that $\hat{\lambda}_{1}\geq\hat{\lambda}_{2}\geq\ldots
\geq0$. We take $(\hat{\lambda}_{j},\hat{\phi}_{j})$ and $\hat{\xi}%
_{ij}=\langle X_{i},\hat{\phi}_{j} \rangle$ to be the estimators of
$(\lambda_{j},\phi_{j})$ and $\xi_{ij},$ respectively, and take
\[
\tilde{a}_{j}=\frac{1}{n\hat{\lambda}_{j}}\sum_{i=1}^{n} \left[ Y_{i}%
-(W_{i}^{T}\pmb{\alpha}_{0}+g(Z_{i}^{T}\pmb{\beta}_{0}))\right] \hat{\xi}%
_{ij}\tag{2.5}
\]
to be the estimator of $a_{j}$.

In order to estimate $g$, we adapt spline approximations. We assume that
$\Vert\pmb{\beta}_{0}\Vert=1$ and that the last element $\beta_{0d}$ of
$\pmb{\beta}_{0}$ is positive, to ensure identifiability. Let
$\pmb{\beta}_{-d}=(\beta_{1},\ldots,\beta_{d-1})^{T}$ and $\pmb{\beta}_{0,-d}%
=(\beta_{01},\ldots,\beta_{0(d-1)})^{T}$. Since $\beta_{0d}=\sqrt
{1-(\beta_{01}^{2}+\cdots+\beta_{0(d-1)}^{2})}>0$, there exists a constant
$\rho_{0}\in(0,1)$ such that $\pmb{\beta}_{0}\in\Theta_{\rho_{0}%
}=\{\pmb{\beta}=(\beta_{1},\ldots,\beta_{d})^{T}:\beta_{d}=\sqrt{1-(\beta
_{1}^{2}+\cdots+\beta_{d-1}^{2})}\geq\rho_{0}\}$. Suppose that the
distribution of $Z$ has a compact support set $\mathcal{D}$. Denote $U_{\ast
}=\inf_{z\in\mathcal{D},\pmb{\beta}\in\Theta_{\rho_{0}}}z^{T}\pmb{\beta}$ and
$U^{\ast}=\sup_{z\in\mathcal{D},\pmb{\beta}\in\Theta_{\rho_{0}}}%
z^{T}\pmb{\beta}$. We first split the interval $[U_{\ast},U^{\ast}]$ into
$k_{n}$ subintervals with knots $\{U_{\ast}=u_{n0}<u_{n1}<\cdots<u_{nk_{n}%
}=U^{\ast}\}$. For fixed $\pmb{\beta}$, suppose $u_{n(l-1)}<\inf
_{z\in\mathcal{D}}z^{T}\pmb{\beta}\leq u_{nl}<u_{n(l+k_{\pmb{\beta}})}\leq
\sup_{z\in\mathcal{D}}z^{T}\pmb{\beta}<u_{n(l+k_{\pmb{\beta} }+1)}$. Let
$U_{\pmb{\beta}}=u_{nl}$ and $U^{\pmb{\beta}}=u_{n(l+k_{\pmb{\beta}})}$. For
any fixed integer $s\geq1$, let $S_{k_{\pmb{\beta}}}^{s}(u)$ be the set of
spline functions of degree $s$ with knots $\{U_{\pmb{\beta}}=u_{nl}%
<u_{n(l+1)}<\cdots<u_{n(l+k_{\pmb{\beta}})}=U^{\pmb{\beta}}\}$; that is, a
function $f(u)$ belongs to $S_{k_{\pmb{\beta}}}^{s}(u)$ if and only if $f(u)$
belongs to $C^{s-1}[u_{nl},u_{n(l+k_{\pmb{\beta}})}]$ and its restriction to
each $[u_{nk},u_{n(k+1)})$ is a polynomial of degree at most $s$. Put
\[
B_{k\beta}(u)=(\tilde{u}_{nk}-\tilde{u}_{n(k-s-1)})[\tilde{u}_{n(k-s-1)}%
,\ldots,\tilde{u}_{nk}](w-u)_{+}^{s},\ \ \ \ k=1,\ldots,K_{\beta},\tag{2.6}
\]
where $K_{\pmb{\beta}}=k_{\pmb{\beta}}+s$, $[\tilde{u}_{n(k-s-1)}%
,\ldots,\tilde{u}_{nk}]f$ denotes the $(s+1)$th-order divided difference of
the function $f$, $\tilde{u}_{nk}=u_{nl}$ for $k=-s,\ldots,-1$, $\tilde
{u}_{nk}=u_{n(l+k)}$ for $k=0,1,\ldots,k_{\pmb{\beta}}$, and $\tilde{u}%
_{nk}=u_{nk_{\pmb{\beta}}}$ for $k=k_{\pmb{\beta} }+1,\ldots,K_{n}$. Then
$\{B_{k\pmb{\beta}}(u)\}_{k=1}^{K_{\pmb{\beta}}}$ form a basis for
$S_{k_{\pmb{\beta}}}^{s}(u)$.

For fixed $\pmb{\alpha}$ and $\pmb{\beta}$, we use $\sum_{j=1}^{m}\tilde
{a}_{j}\hat{\xi}_{j}$ to approximate $\sum_{j=1}^{\infty}a_{j}\xi_{j}$ in
(2.2) and use $\sum_{k=1}^{K_{\pmb{\beta}}}b_{k}B_{k\pmb{\beta}}(u)$ to
approximate $g(u)$ for $u\in\lbrack U_{\pmb{\beta}},U^{\pmb{\beta}}]$. We then
estimate $g(\cdot)$ by minimizing
\begin{align*}
\sum_{i=1}^{n}\Bigg\{Y_{i}-&\sum_{j=1}^{m}\frac{\hat{\xi}_{ij}}{n\hat
{\lambda}_{j}}\sum_{l=1}^{n}\Bigg[  Y_{l}-W_{l}^{T}\pmb{\alpha}-\sum
_{k=1}^{K_{\pmb{\beta}}}b_{k}B_{k\pmb{\beta}}(Z_{l}^{T}\pmb{\beta})\Bigg]
\hat{\xi}_{lj}-\\
&W_{i}^{T}\pmb{\alpha}-\sum_{k=1}^{K_{\pmb{\beta}}}%
b_{k}B_{k\pmb{\beta}}(Z_{i}^{T}\pmb{\beta})\Bigg\} ^{2}\tag{2.7}
\end{align*}
with respect to $b_{1},\ldots,b_{K_{\pmb{\beta}}}$, where $m$ is a smoothing
parameter which denotes a frequency cut-off. Define $\tilde{\xi}_{il}%
=\sum_{j=1}^{m}\hat{\xi}_{ij}\hat{\xi}_{lj}/\hat{\lambda}_{j}$, $\tilde{Y}%
_{i}=Y_{i}-\frac{1}{n}\sum_{l=1}^{n}Y_{l}\tilde{\xi}_{il}$, $\tilde{W}%
_{i}=W_{i}-\frac{1}{n}\sum_{l=1}^{n}W_{l}\tilde{\xi}_{il}$ and $\tilde
{B}_{k\pmb{\beta}}(Z_{i}^{T}\pmb{\beta})=B_{k\pmb{\beta}}(Z_{i}^{T}%
\pmb{\beta})-\frac{1}{n}\sum_{l=1}^{n}B_{k\pmb{\beta}}(Z_{l}^{T}%
\pmb{\beta})\tilde{\xi}_{il}$. Then (2.7) can be written as
\[
\sum_{i=1}^{n}\left\{  \tilde{Y}_{i}-\tilde{W}_{i}^{T}\pmb{\alpha}-\sum
_{k=1}^{K_{\pmb{\beta}}}b_{k}\tilde{B}_{k\pmb{\beta}}(Z_{i}^{T}%
\pmb{\beta})\right\}  ^{2}.\tag{2.8}
\]
Denote $\tilde{\pmb{B}}_{\pmb{\beta}}(Z_{i}^{T}\pmb{\beta})=(\tilde
{B}_{1\pmb{\beta}}(Z_{i}^{T}\pmb{\beta}),\ldots,\tilde{B}_{K_{\pmb{\beta}}%
\pmb{\beta}}(Z_{i}^{T}\pmb{\beta}))^{T}$, $\tilde{\pmb{B}}%
(\pmb{\beta})=(\tilde{\pmb{B}}_{\pmb{\beta}}(Z_{1}^{T}\pmb{\beta}),\ldots
,$ $\tilde{\pmb{B}}_{\pmb{\beta}}(Z_{n}^{T}\pmb{\beta}))^{T}$, $\tilde
{\pmb{Y}}=(\tilde{Y}_{1},\ldots,\tilde{Y}_{n})^{T}$, $\tilde{\pmb{W}}%
=(\tilde{W}_{1},\ldots,\tilde{W}_{n})^{T}$ and $\pmb{b}=(b_{1},\ldots
,b_{K_{\pmb{\beta}}})^{T}$. If $\tilde{\pmb{B}}^{T}(\pmb{\beta})\tilde
{\pmb{B}}(\pmb{\beta})$ is invertible, then the estimator $$\tilde
{\pmb{b}}(\pmb{\alpha},\pmb{\beta})=(\tilde{b}_{1}%
(\pmb{\alpha},\pmb{\beta}),\ldots,\tilde{b}_{K_{\pmb{\beta}}}%
(\pmb{\alpha},\pmb{\beta}))^{T}$$ of $\pmb{b}$ is given by
\[
\tilde{\pmb{b}}(\pmb{\alpha},\pmb{\beta})=\left\{  \tilde{\pmb{B}}%
^{T}(\pmb{\beta})\tilde{\pmb{B}}(\pmb{\beta})\right\}  ^{-1}\tilde
{\pmb{B}}^{T}(\pmb{\beta})(\tilde{\pmb{Y}}-\tilde{\pmb{W}}%
\pmb{\alpha}).\tag{2.9}
\]
We solve the following minimization problem
\[
\min_{\pmb{\alpha},\pmb{\beta}}\left\{  \tilde{\pmb{Y}}-\tilde{\pmb{W}}%
\pmb{\alpha}-\tilde{\pmb{B}}(\pmb{\beta})\tilde{\pmb{b}}%
(\pmb{\alpha},\pmb{\beta})\right\}  ^{T}\left\{  \tilde{\pmb{Y}}%
-\tilde{\pmb{W}}\pmb{\alpha}-\tilde{\pmb{B}}(\pmb{\beta})\tilde{\pmb{b}}%
(\pmb{\alpha},\pmb{\beta})\right\}  \tag{2.10}
\]
to obtain the estimators $\hat{\pmb{\alpha}}$ and $\hat{\pmb{\beta}}$. A
Newton-Raphson algorithm can be applied for the minimization. An estimator of
$\pmb{b}$ is obtained by solving the following minimization problem
\[
\hat{\pmb{b}}=\min_{\pmb{b}}\sum_{i=1}^{n}\left\{  \tilde{Y}_{i}-\tilde{W}%
_{i}^{T}\hat{\pmb{\alpha}}-\pmb{b}^{T}\tilde{\pmb{B}}_{\hat{\pmb{\beta}}%
}(Z_{i}^{T}\hat{\pmb{\beta}})\right\}  ^{2},\tag{2.11}
\]
and then $\hat{\pmb{b}}$ is given by
\[
\hat{\pmb{b}}=\tilde{\pmb{b}}(\hat{\pmb{\alpha}},\hat{\pmb{\beta}})=\left\{
\tilde{\pmb{B}}^{T}(\hat{\pmb{\beta}})\tilde{\pmb{B}}(\hat{\pmb{\beta}}%
)\right\}  ^{-1}\tilde{\pmb{B}}^{T}(\hat{\pmb{\beta}})(\tilde{\pmb{Y}}%
-\tilde{\pmb{W}}^{T}\hat{\pmb{\alpha}}).\tag{2.12}
\]
Let $\tilde{g}(u)=\sum_{k=1}^{K_{\hat{\pmb{\beta}}}}\hat{b}_{k}B_{k\hat
{\pmb{\beta}}}(u)$ for $u\in\lbrack U_{\hat{\pmb{\beta}}},U^{\hat
{\pmb{\beta}}}]$. We then choose a new tuning parameter $\tilde{m}$ and an
estimator of $a(t)$ given by $\hat{a}(t)=\sum_{j=1}^{\tilde{m}}\hat{a}_{j}%
\hat{\phi}_{j}(t)$ with
\[
\hat{a}_{j}=\frac{1}{n\hat{\lambda}_{j}}\sum_{i=1}^{n}\left\{  Y_{i}-W_{i}%
^{T}\hat{\pmb{\alpha}}-\tilde{g}(Z_{i}^{T}\hat{\pmb{\beta}})\right\}  \hat
{\xi}_{ij}.\tag{2.13}
\]

In order to construct an estimator of $g$ that achieves the optimal rate of
convergence, we select new knots and new B-spline basis based on the
estimators $\hat{\pmb{\alpha}}$ and $\hat{\pmb{\beta}}.$ Let $\{U_{\hat
{\pmb{\beta}}}=\bar{u}_{n0}<\bar{u}_{n1}<\cdots<\bar{u}_{nk^{\ast}_{n}%
)}=U^{\hat{\pmb{\beta}}}\}$ be new knots and $\{B^{\ast}_{k}(u)\}_{k=1}%
^{K^{\ast}_{n}}$ be a new basis, where $K^{\ast}_{n}=k^{\ast}_{n}+s$. Then
$B^{\ast}_{k\pmb{\beta}}(Z_{i}^{T}\pmb{\beta})$, $\pmb{B}^{\ast}_{\pmb{\beta}}(Z_{i}%
^{T}\pmb{\beta})$ and $\pmb{B}^{\ast}(\pmb{\beta})$ are defined similarly as
$\tilde{B}_{k\pmb{\beta}}(Z_{i}^{T}\pmb{\beta})$, $\tilde{\pmb{B}}%
_{\pmb{\beta}}(Z_{i}^{T}\pmb{\beta})$ and $\tilde{\pmb{B}}(\pmb{\beta}),$
respectively. We then solve the following minimization problem
\[
\min_{\pmb{b}^{\ast}}\sum_{i=1}^{n}\left\{ \tilde{Y}_{i}-\tilde{W}_{i}^{T}%
\hat{\pmb{\alpha} }-{\pmb{b}^{\ast}}^{T}\pmb{B}^{\ast}_{\hat{\pmb{\beta}}}%
(Z_{i}^{T}\hat{\pmb{\beta}})\right\} ^{2},\tag{2.14}
\]
to obtain an estimator of $\pmb{b}^{\ast}$, where $\pmb{b}^{\ast}=(b_{1}%
,\ldots,b_{K^{\ast}_{n}})^{T}$. If ${\pmb{B}^{\ast}}^{T}(\hat{\pmb{\beta}}%
)\pmb{B}^{\ast}(\hat{\pmb{\beta}})$ is invertible, then an estimator of
$\pmb{b}^{\ast}$ is given by
\[
\hat{\pmb{b}}^{\ast}=\pmb{b}^{\ast}(\hat{\pmb{\alpha}},\hat{\pmb{\beta}})=\left\{
{\pmb{B}^{\ast}}^{T}(\hat{\pmb{\beta}})\pmb{B}^{\ast}(\hat{\pmb{\beta}})\right\}
^{-1}{\pmb{B}^{\ast}}^{T}(\hat{\pmb{\beta}})(\tilde{\pmb{Y}}-\tilde{\pmb{W}}%
^{T}\hat{\pmb{\alpha} }).\tag{2.15}
\]
The second stage estimator of $g(u)$ is then equal to $\hat{g}(u)=\sum
_{k=1}^{K^{\ast}_{n}}\hat{{b}}_{k}^{\ast}B^{\ast}_{k}(u)$ for $u\in\lbrack
U_{\hat{\pmb{\beta}}},U^{\hat{\pmb{\beta}}}]$.

To implement our estimation method, some appropriate values for $m$,
$\tilde{m}$, $k_{n}$ and $K^{\ast}_{n}$ are necessary. From our simulation in
Section 4 below, we observe that the parametric estimators $\hat
{\pmb{\alpha}}$ and $\hat{\pmb{\beta}}$ are not sensitive to the choices of
$m$ and $k_{n}$, they can be chosen subjectively. In the simulation in Section
4, we also choose $h_{0}=n^{-1/(2s-1)}$ with $s=3$, where $s$ is defined in
Assumption 4 in Section 3 below. The value for tuning parameter $\tilde{m}$
can be selected by information criteria BIC, which is given by
\[
BIC(\tilde{m})=\log\left\{  \frac{1}{n}\sum_{i=1}^{n}\left(  Y_{i}-W_{i}%
\hat{\pmb{\alpha}}-\sum_{j=1}^{\tilde{m}}\hat{a}_{j}\hat{\xi}_{ij}-\tilde
{g}(Z_{i}^{T}\hat{\pmb{\beta}})\right)  ^{2}\right\}  +\frac{\log(n)\tilde{m}%
}{n}.
\]
Large values of $BIC$ indicate either poor fidelity to the data or overfitting
because $\tilde{m}$ is too large. A value for $K^{\ast}_{n}$ can also be
selected by the following BIC information criteria:
\[
BIC(K^{\ast}_{n})=\log\left\{  \frac{1}{n}\sum_{i=1}^{n}\left(  \tilde{Y}%
_{i}-\tilde{W}_{i}\hat{\pmb{\alpha}}-\hat{\pmb{b}}^{\ast T}\pmb{B}^{\ast}%
_{\hat{\pmb{\beta}}}(Z_{i}^{T}\hat{\pmb{\beta}})\right)  ^{2}\right\}
+\frac{\log(n)K^{\ast}_{n}}{n}.
\]
In practice, the proposed estimation method is implemented using the following steps:

\textbf{Step 1.} Choose an $m$ and fit a partial functional linear model; that
is, solve the minimization problem (2.8) with the link function $g$ replaced
by a linear function to obtain initial values $\hat{\pmb{\alpha}}^{(0)}$ and
$\hat{\pmb{\beta}}_{1}^{(0)}$. Then set $\hat{\pmb{\beta}}^{(0)}%
=\hat{\pmb{\beta}}_{1}^{(0)}/\Vert\hat{\pmb{\beta}}_{1}^{(0)}\Vert$, and
multiply it by $-1$ if necessary.

\textbf{Step 2.} Construct the B-spline basis $\{B_{k\hat
{\pmb{\beta}}^{(0)}}(u)\}_{k=1}^{K_{\hat{\pmb{\beta}}^{(0)}}}$ based on 
the computed  $U_{\hat{\pmb{\beta}}^{(0)}}$ and $U^{\hat
{\pmb{\beta}}^{(0)}}$. Then obtain
$\tilde{b}(\hat{\pmb{\alpha}}^{(0)},\hat{\pmb{\beta}}^{(0)})$ from (2.9) and
solve the minimizing problem (2.10) to obtain the estimators $\hat
{\pmb{\alpha}}$ and $\hat{\pmb{\beta}}$.

\textbf{Step 3.} Compute $\hat{b}$ and $\hat{a}_{j}$ from (2.12) and (2.13),
respectively, and obtain the estimator $\hat{a}(t)$.

\textbf{Step 4.} Compute $U_{\hat{\pmb{\beta}}}$ and $U^{\hat{\pmb{\beta}}},$
and construct the basis $\{B^{\ast}_{k}(u)\}_{k=1}^{K^{\ast}_{n}}$. Then obtain
the estimator $\hat{\mathbf{b}}^{\ast}$ from (2.15) and obtain the estimator $\hat
{g}(u)$.

\medskip

\textbf{Remark\ 2.1.} In practical applications, $X(t)$ is only discretely
observed. Without loss of generality, suppose for each $i=1,\ldots,n$,
$X_{i}(t)$ is observed at $n_{i}$ discrete points $0=t_{i1}<\ldots<t_{in_{i}%
}=1$. Then linear interpolation functions or spline interpolation functions
can be used for the estimators of $X_{i}(t)$.

\textbf{Remark\ 2.2} \ Though the basis function $B_{k\pmb{\beta}}(u)$ depends
on $\pmb{\beta}$, we see from (2.6) that the total number of all the different
$B_{k\pmb{\beta}}(u)$ is not more than $(s+1)k_{n}$. In certain practical
applications where the sample size $n$ is not large enough and $h_{0}$ is not
small enough, one can choose $U_{\pmb{\beta}}=\inf_{z\in\mathcal{D}}%
z^{T}\pmb{\beta}$ and $U^{\pmb{\beta}}=\sup_{z\in\mathcal{D}}z^{T}\pmb{\beta}$
and construct the basis $\{B_{k\pmb{\beta}}(u)\}_{k=1}^{K_{\pmb{\beta}}}$ with
knots $\{U_{\pmb{\beta}}<u_{n(l+1)}<\cdots<u_{n(l+k_{\pmb{\beta}}%
-1)}<U^{\pmb{\beta}}\}$ to make full use of the data. That is, the intervals
$[u_{nl},u_{n(l+1)}]$ and $[u_{n(l+k_{\pmb{\beta}}-1)},u_{n(l+k_{\pmb{\beta}}%
)}]$ are replaced by $[U_{\pmb{\beta}},u_{n(l+1)}]$ and
$[u_{n(l+k_{\pmb{\beta}}-1)},U^{\pmb{\beta}}],$ respectively.

\section{Asymptotic properties}

In this section we establish the asymptotic normality and convergence rates of
the estimators proposed in the previous section. Before stating main results,
we first state a few assumptions that are necessary to prove the theoretical results.

\textbf{Assumption 1.}\ $E(Y^{4})<+\infty$ and $\int_{\mathcal{T}}%
E(X^{4}(t))dt<\infty$. $E(\xi_{j}|Z^{T}\pmb{\beta})=0$ and $E(\xi_{i}\xi
_{j}|Z^{T}\pmb{\beta})=0$ for $i\neq j,\ i,j=1,2,\ldots;$ and $\pmb{\beta}\in
\Theta_{\rho_{0}}$. For each $j\geq1$, $E(\xi_{j}^{2r}|Z^{T}\pmb{\beta})\leq
C_{1}\lambda_{j}^{r}$ for $r=1,2$, where $C_{1}>0$ is a constant. For any
sequence $j_{1},\ldots,j_{4}$, $E(\xi_{j_{1}}\ldots\xi_{j_{4}}|Z^{T}%
\pmb{\beta})$ $=0$ unless each index $j_{k}$ is repeated.

\textbf{Assumption 2.} There exists a convex function $\varphi$ defined on the
interval $[0,1]$ such that $\varphi(0) = 0$ and $\lambda_{j}=\varphi(1/j)$ for
$j\geq1$.

\textbf{Assumption 3.}\ For Fourier coefficients $a_{j}$, there exist
constants $C_{2}>0$ and $\gamma>3/2$ such that $|a_{j}|\leq C_{2}j^{-\gamma}$
for all $j\geq1$.

\textbf{Assumption 4.} The function $g(u)$ is a $s$-times continuously
differentiable function such that $|g^{(s)}(u^{\prime}) - g^{(s)}(u) |\leq
C_{3}|u^{\prime} -u|^{\varsigma}$, for $U_{\ast}\leq u^{\prime},u\leq U^{\ast
}$ and $p=s+\varsigma>3$, with constants $0<\varsigma\leq1$ and $C_{3}>0.$ The
knots $\{U_{\ast}=u_{n0}<u_{n1}<\cdots<u_{nk_{n}}=U^{\ast}\}$ satisfy that
$h_{0}/\min_{1\leq k\leq k_{n}}h_{nk}\leq C_{4}$, where $h_{nk}=u_{nk}%
-u_{n(k-1)},h_{0}=\max_{1\leq k\leq k_{n}}h_{nk}$ and $C_{4}>0$ is a constant.

\textbf{Assumption 5. }$nh_{0}^{2p}\rightarrow0$, $n^{-1/2}m\lambda_{m}%
^{-1}\rightarrow0$, $n^{-1}m^{4}\lambda_{m}^{-1}h_{0}^{-6}\log m\rightarrow0$
and $m^{-2\gamma}h_{0}^{-2}\rightarrow0$.

\textbf{Assumption 5'. }$m\rightarrow\infty$, $h_{0}\rightarrow0$,
$n^{-1/2}m\lambda_{m}^{-1}\rightarrow0$, $n^{-1}m^{4}\lambda_{m}^{-1}%
h_{0}^{-2}\log m$ $\rightarrow0$ and $(nh_{0}^{3})^{-1}(\log n)^{2}\rightarrow0$.

\textbf{Assumption 6. }The distribution of $Z$ has a compact support set
$\mathcal{D}$. The marginal density function $f_{\pmb{\beta}}(u)$ of
$Z^{T}\pmb{\beta}$ is bounded away from zero and infinity for $u\in\lbrack
U_{\pmb{\beta}},U^{\pmb{\beta}}]$ and satisfies that $0<c_{1}\leq
f_{\pmb{\beta}}(u)\leq C_{5}<+\infty$ for $\pmb{\beta}$ in a small
neighborhood of $\pmb{\beta}_{0}$ and $u\in\lbrack U_{\pmb{\beta}_{0}%
},U^{\pmb{\beta}_{0}}]$, where $c_{1}$ and $C_{5}$ are two positive constants.

\textbf{Assumption 7. }$W=(W_{1},\ldots,W_{q})^{T}$, $W_{r}=\check{W}%
_{r}+V_{r}$, $\check{W}_{r}=\sum_{j=1}^{\infty}w_{rj}\xi_{j}$ and
$|w_{rj}|\leq C_{6}j^{-\gamma}$ for all $j\geq1$ and $r=1,\ldots,q$, where
$C_{6}>0$ is a constant. $V=(V_{1},V_{2},\ldots,V_{q})^{T}$ is independent of
$\{\xi_{j}, \, j=1,\ldots\}$ and $E(\Vert V\Vert^{4})<+\infty$.

Under Assumption 4, according to Corollary 6.21 of Schumaker (1981, p.227),
there exists a spline function $g_{0}(u)=\sum_{k=1}^{K_{\pmb{\beta}_{0}}%
}b_{0k}B_{k\pmb{\beta}_{0}}(u)$ and a constant $C_{7}>0$ such that, for $k=0,1,\ldots,s$,
\[
\sup_{u\in\lbrack U_{\pmb{\beta}_{0}},U^{\pmb{\beta}_{0}}]}|R^{(k)}(u)|\leq
C_{7}h_{0}^{p-k}\tag{3.1}
\]
 where $R(u)=g(u)-g_{0}(u)$. Let $\pmb{B}_{\pmb{\beta}}%
(u)=(B_{1\pmb{\beta}}(u),\ldots,B_{K_{\pmb{\beta}}\pmb{\beta}}(u))^{T}$ and
$\pmb{b}_{0}=(b_{01},\ldots,b_{0K_{\pmb{\beta}_{0}}})^{T}$. Define
\[%
\begin{array}
[c]{ll}%
G(\pmb{\alpha},\pmb{\beta}) & =\left\{(\pmb{\alpha}-\pmb{\alpha}_{0})^{T}%
E(VV^{T})-2\pmb{b}_{0}^{T}%
E[\pmb{B}_{\pmb{\beta}_{0}}(Z^{T}\pmb{\beta}_{0})V^{T}%
]\right\}(\pmb{\alpha}-\pmb{\alpha}_{0})\\
& +\pmb{b}_{0}^{T}\Gamma(\pmb{\beta}_{0},\pmb{\beta}_{0})\pmb{b}_{0}-\Pi
^{T}(\pmb{\alpha},\pmb{\beta})\Gamma^{-1}(\pmb{\beta},\pmb{\beta})\Pi
(\pmb{\alpha},\pmb{\beta})+\sigma^{2},
\end{array}
\tag{3.2}
\]
where $\Gamma(\pmb{\beta}_{1},\pmb{\beta}_{2})=(\gamma_{kk^{\prime}%
}(\pmb{\beta}_{1},\pmb{\beta} _{2}))_{K_{\pmb{\beta}_{1}}\times
K_{\pmb{\beta}_{2}}}$ with $\gamma_{kk^{\prime}}(\pmb{\beta}_{1}%
,\pmb{\beta}_{2})=E[B_{k\pmb{\beta}_{1}}(Z^{T}\pmb{\beta}_{1})$ $B_{k^{\prime
}\pmb{\beta}_{2}}(Z^{T}\pmb{\beta}_{2})]$ and $\Pi
(\pmb{\alpha},\pmb{\beta})=\Gamma(\pmb{\beta},\pmb{\beta}_{0})\pmb{b}_{0}%
-E[\pmb{B}_{\pmb{\beta}}(Z^{T}\pmb{\beta})V^{T}](\pmb{\alpha}-\pmb{\alpha}_{0}%
)$. Put $\pmb{\theta} =(\pmb{\alpha}^{T},\pmb{\beta}^{T})^{T}$,
$\pmb{\theta}_{-d}=(\pmb{\alpha}^{T},\pmb{\beta}_{-d}^{T})^{T}$,$\ \hat
{\pmb{\theta}}_{-d}=(\hat{\pmb{\alpha}}^{T},\hat{\pmb{\beta}}_{-d}^{T})^{T}$
and $\pmb{\theta}_{0,-d}=(\pmb{\alpha}_{0}^{T},\pmb{\beta}_{0,-d}^{T})^{T}$.
Define
\[
G^{\ast}(\pmb{\theta}_{-d})=G^{\ast}(\pmb{\alpha},\pmb{\beta}_{-d}%
)=G(\pmb{\alpha},\beta_{1},\ldots,\beta_{d-1},\sqrt{1-\Vert\pmb{\beta}_{-d}%
\Vert^{2}})
\]
and its Hessian matrix $H^{\ast}(\pmb{\theta}_{-d})=\frac{\partial^{2}%
}{\partial\pmb{\theta}_{-d}\partial\pmb{\theta}_{-d}^{T}}G^{\ast
}(\pmb{\theta}_{-d}).$

\textbf{Assumption 8.} $G^{\ast}(\pmb{\theta}_{-d})$ is locally convex at
$\pmb{\theta}_{0,-d}$ such that for any $\varepsilon>0$, there exists some
$\epsilon>0$ such that $\Vert\pmb{\theta}_{-d}-\pmb{\theta}_{0,-d}%
\Vert<\varepsilon$ holds whenever $|G^{\ast}(\pmb{\theta}_{-d})-G^{\ast
}(\pmb{\theta}_{0,-d})|<\epsilon$. Furthermore, the Hessian matrix $H^{\ast
}(\pmb{\theta}_{-d})$ is continuous in some neighborhood of
$\pmb{\theta}_{0,-d}$ and $H^{\ast}(\pmb{\theta}_{0,-d})>0$.

\textbf{Assumption 9.} The knots $\{U_{\hat{\pmb{\beta}}}=\bar{u}_{n0}<\bar
{u}_{n1}<\cdots<\bar{u}_{n\vec{k}_{n})}=U^{\hat{\pmb{\beta}}}\}$ satisfy that
$h/\min_{1\leq k\leq\vec{k}_{n}}\bar{h}_{nk}\leq C_{8}$, where $\bar{h}%
_{nk}=\bar{u}_{nk}-\bar{u}_{n(k-1)},h=\max_{1\leq k\leq\vec{k}_{n}}\bar
{h}_{nk}$ and $C_{8}>0$ is a constant. Further, $h\rightarrow0$ and
$n^{-1}m^{4}\lambda_{m}^{-1}h^{-4}\log m\rightarrow0$.

Assumptions 1 and 3 are standard conditions for functional linear models; see,
e.g., Cai and Hall \cite{r3} 
and Hall and Horowitz \cite{r11}.
Assumption 2 is
slightly less restrictive than (3.2) of Hall and Horowitz \cite{r11}. 
The quantity
$p$ in Assumption 4 is the order of smoothness of the function $g(u)$.
Assumptions 5 and 5' can be easily verified and will be further discussed
below. Assumption 6 ensures the existence and uniqueness of the spline
estimator of the function $g(u)$. If the marginal density $f_{\pmb{\beta}}(u)$
of $Z^{T}\pmb{\beta}$ is uniformly continuous for $\pmb{\beta}$ in some
neighborhood of $\pmb{\beta}_{0},$ then the second part of Assumption 6 is
easily satisfied by modifying the knots. Assumption 8 ensures the existence
and uniqueness of the estimator of $\pmb{\theta}_{0,-d}$ in a neighborhood of
$\pmb{\theta}_{0,-d}$.

\textbf{Remark\ 3.1.} \ If $\lambda_{j}\sim j^{-\delta}$, $m\sim n^{\iota}$
and $h_{0}\sim n^{-\tau}$, then Assumption 5 holds when $\iota<\min
(1/(2(1+\delta)),1/(\delta+4))$ and $1/(2p)<\tau<(1-\iota(\delta+4))/6$, where
$\delta>1$, $\iota>0$ and $\tau>0$ are constants and the notation $a_{n}\sim
b_{n}$ means that the ratio $a_{n}/b_{n}$ is bounded away from zero and infinity.

The next theorem gives the consistency and convergence rate of the estimators
of $\pmb{\alpha}_{0}$ and $\pmb{\beta}_{0,-d}$.

\textbf{Theorem\ 3.1.}\ (i) Suppose that Assumptions 1 to 4, 5', 6 and 7 hold,
and that $G^{\ast}(\pmb{\theta}_{-d})$ is locally convex at
$\pmb{\theta}_{0,-d}$. Then, as $n\rightarrow\infty$,
\[
\hat{\pmb{\alpha}}\overset{P}{\rightarrow}\pmb{\alpha}_{0},\ \ \ \hat
{\pmb{\beta}}_{-d}\overset{P}{\rightarrow}\pmb{\beta}_{0,-d},\tag{3.3}
\]
where $\overset{P}{\rightarrow}$ means convergence in probability.

(ii) \ Suppose that Assumptions 1 to 8 hold. Then
\[
\hat{\pmb{\alpha}}-\pmb{\alpha}_{0}=o_{p}(h_{0}),\ \ \hat{\pmb{\beta}}%
_{-d}-\pmb{\beta}_{0,-d}=o_{p}(h_{0}).\tag{3.4}
\]

In order to establish the asymptotic distributions of the estimators
$\hat{\pmb{\alpha}}$ and $\hat{\pmb{\beta}}_{-d}$, we first introduce some
notation. Define
\[
G_{n}(\pmb{\theta})=G_{n}(\pmb{\alpha},\pmb{\beta})=\frac{1}{n}\sum_{i=1}%
^{n}\left\{ \tilde{Y}_{i}-\tilde{W}_{i}^{T}\pmb{\alpha}-\sum_{k=1}%
^{K_{\pmb{\beta}}}\tilde{b}_{k}(\pmb{\alpha},\pmb{\beta})\tilde{B}%
_{k\pmb{\beta}}(Z_{i}^{T}\pmb{\beta})\right\} ^{2}.\tag{3.5}
\]
If $u_{n(l-1)}<\inf_{z\in\mathcal{D}}z^{T}\pmb{\beta}_{0}<u_{nl},$ then by
(3.4) we have $U_{\hat{\pmb{\beta}}}=U_{\pmb{\beta}_{0}}=u_{nl}$ for
sufficiently large $n$. If $\inf_{z\in\mathcal{D}}z^{T}\pmb{\beta}_{0}=u_{nl}%
$, then we modify $u_{nl}$ such that $\inf_{z\in\mathcal{D}}z^{T}%
\pmb{\beta}_{0}<u_{nl}$, and also we then have $U_{\hat{\pmb{\beta}}%
}=U_{\pmb{\beta}_{0}}=u_{nl}$. Similarly, if $\sup_{z\in\mathcal{D}}%
z^{T}\pmb{\beta}_{0}=u_{n(l+k_{\pmb{\beta}})}$, then we modify
$u_{n(l+k_{\pmb{\beta}})}$ such that $u_{n(l+k_{\pmb{\beta}})}<\sup
_{z\in\mathcal{D}}z^{T}\pmb{\beta}_{0}$, and then we have $U^{\hat
{\pmb{\beta}}}=U^{\pmb{\beta}_{0}}=u_{n(l+k_{\pmb{\beta}})}$. Therefore, if
necessary, we first modify the knots $\{u_{nk}\}_{k=0}^{k_{_{n}}}$, so that
there exists a neighborhood $\delta^{\ast}(\pmb{\beta}_{0,-d};r^{\ast})$ of
$\pmb{\beta}_{0,-d}$ such that $U_{\pmb{\beta}}=U_{\pmb{\beta}_{0}}$,
$U^{\pmb{\beta}}=U^{\pmb{\beta}_{0}}$ for $\pmb{\beta}\in\delta^{\ast
}(\pmb{\beta}_{0,-d};r^{\ast})$ and $\hat{\pmb{\beta}}\in\delta^{\ast
}(\pmb{\beta}_{0,-d};r^{\ast})$ for sufficiently large $n$. Let $K_{n}%
=K_{\pmb{\beta}_{0}}$, $B_{k}(u)=B_{k\pmb{\beta}_{0}}(u)$ and $\tilde{B}%
_{k}(u)=\tilde{B}_{k\pmb{\beta}_{0}}(u)$. For $\pmb{\beta}\in\delta^{\ast
}(\pmb{\beta}_{0,-d};r^{\ast})$, we have $K_{\pmb{\beta}}=K_{n}$,
$B_{k}(u)=B_{k\pmb{\beta}}(u)$ and $\tilde{B}_{k}(u)=\tilde{B}_{k\pmb{\beta}}%
(u)$. Further, we have
\[%
\begin{array}
[c]{ll}%
G_{n}(\pmb{\alpha},\pmb{\beta}) & =\frac{1}{n}\sum_{i=1}^{n}\left\{ \tilde
{Y}_{i}-\tilde{W}_{i}^{T}\pmb{\alpha}-\sum_{k=1}^{K_{n}}\tilde{b}%
_{k}(\pmb{\alpha},\pmb{\beta})\tilde{B}_{k}(Z_{i}^{T}\pmb{\beta})\right\}
^{2}\\
& =\frac{1}{n}\left\{ \tilde{\pmb{Y}}-\tilde{\pmb{W}}\pmb{\alpha}-\tilde
{\pmb{B}}(\pmb{\beta})\tilde{\pmb{b}}(\pmb{\alpha},\pmb{\beta})\right\}
^{T}\left\{ \tilde{\pmb{Y}}-\tilde{\pmb{W}}\pmb{\alpha}-\tilde{\pmb{B}}%
(\pmb{\beta})\tilde{\pmb{b}}(\pmb{\alpha},\pmb{\beta})\right\} ,
\end{array}
\]
$G_{n}(\pmb{\theta}_{-d},\pmb{b})=G_{n}(\pmb{\alpha},\pmb{\beta}_{-d}%
,\pmb{b})=\frac{1}{n}\Big\{ \tilde{\pmb{Y}}-\tilde{\pmb{W}}%
\pmb{\alpha}-\tilde{\pmb{B}}(\pmb{\beta}_{-d})\pmb{b}\Big\} ^{T}\Big\{
\tilde{\pmb{Y}}-\tilde{\pmb{W}}\pmb{\alpha}-\tilde{\pmb{B}}(\pmb{\beta}_{-d}%
)\pmb{b}\Big\} $, where
\[
\tilde{\pmb{B}}(\pmb{\beta}_{-d})=\tilde{\pmb{B}}(\beta_{1},\ldots,\beta
_{d-1},\sqrt{1-(\beta_{1}^{2}+\ldots+\beta_{d-1}^{2})}).
\]
Since $(\hat{\pmb{\alpha}},\hat{\pmb{\beta}})$ is the minimizer of
$G_{n}(\pmb{\alpha},\pmb{\beta})$, then $(\hat{\pmb{\alpha}},\hat
{\pmb{\beta}}_{-d},\hat{\pmb{b}})$ is the minimizer of $G_{n}%
(\pmb{\alpha},\pmb{\beta}_{-d},\pmb{b})$, where $\hat{\pmb{b}}=\tilde
{\pmb{b}}(\hat{\pmb{\theta}}_{-d})=\tilde{\pmb{b}}(\hat{\pmb{\alpha}}%
,\hat{\pmb{\beta}}_{-d})=\left\{ \tilde{\pmb{B}}^{T}(\hat{\pmb{\beta}}%
_{-d})\tilde{\pmb{B}}(\hat{\pmb{\beta}}_{-d})\right\} ^{-1}(\tilde
{\pmb{B}}^{T}\hat{\pmb{\beta}}_{-d})$ $(\tilde{\pmb{Y}}-\tilde{\pmb{W}}%
\hat{\pmb{\alpha}})$. Hence,
\begin{align*}
\frac{\partial G_{n}(\pmb{\alpha},\pmb{\beta}_{-d},\pmb{b})}{\partial
\pmb{\alpha}}\bigg|_{(\pmb{\alpha},\pmb{\beta}_{-d},\pmb{b})=(\hat
{\pmb{\alpha}},\hat{\pmb{\beta}}_{-d},\hat{\pmb{b}})}&=&\\&-&\frac{2}{n}%
\tilde{\pmb{W}}^{T}\left\{ \tilde{\pmb{Y}}-\tilde{\pmb{W}}\hat{\pmb{\alpha}}%
-\tilde{\pmb{B}}(\hat{\pmb{\beta}}_{-d})\hat{\pmb{b}}\right\} =0\tag{3.6}
\end{align*}
\begin{align*}
\frac{\partial G_{n}(\pmb{\alpha},\pmb{\beta}_{-d},\pmb{b})}{\partial\beta
_{r}}\bigg|_{(\pmb{\alpha},\pmb{\beta}_{-d},\pmb{b})=(\hat{\pmb{\alpha}}%
,\hat{\pmb{\beta}}_{-d},\hat{\pmb{b}})}&=&\\&-&\frac{2}{n}\left\{ \tilde
{\pmb{Y}}-\tilde{\pmb{W}}\hat{\pmb{\alpha}}-\tilde{\pmb{B}}(\hat
{\pmb{\beta}}_{-d})\hat{\pmb{b}}\right\} ^{T}\dot{\tilde{B}}_{r}%
(\hat{\pmb{\beta}}_{-d})\hat{\pmb{b}}=0\tag{3.7}
\end{align*}
for $r=1,\ldots,d-1$, where $\dot{\tilde{B}}_{r}(\pmb{\beta}_{-d}%
)=\frac{\partial\tilde{\pmb{B}}(\pmb{\beta}_{-d})}{\partial\pmb{\beta}_{r}}$.
Set $\dot{G}_{n}(\pmb{\theta}_{-d},\pmb{b})=\frac{\partial G_{n}%
(\pmb{\theta}_{-d},\pmb{b})}{\partial\pmb{\theta}_{-d}}=(\frac{\partial
}{\partial\pmb{\alpha}}G_{n}(\pmb{\alpha},\pmb{\beta}_{-d},\pmb{b})^{T}%
$,\newline$\frac{\partial}{\partial\pmb{\beta}_{-d}}G_{n}%
(\pmb{\alpha},\pmb{\beta}_{-d},\pmb{b})^{T})^{T}$. Then from (3.6) and (3.7)
and using a Taylor expansion, we obtain
\[
\dot{G}_{n}(\pmb{\theta}_{0,-d},\tilde{\pmb{b}}(\pmb{\theta}_{0,-d}))+\ddot
{G}_{n}(\pmb{\theta}_{-d}^{\ast},\tilde{\pmb{b}}(\pmb{\theta}_{-d}^{\ast
}))(\hat{\pmb{\theta}}_{-d}-\pmb{\theta}_{0,-d})=0,\tag{3.8}
\]
where $\ddot{G}_{n}(\pmb{\theta}_{-d},\tilde{\pmb{b}}(\pmb{\theta}_{-d}%
))=\frac{\partial}{\partial\pmb{\theta}_{-d}}\dot{G}_{n}(\pmb{\theta}_{-d}%
,\tilde{\pmb{b}}(\pmb{\theta}_{-d}))$ is a $(q+d-1)\times(q+d-1)$ matrix and
$\pmb{\theta}_{-d}^{\ast}$ is between $\hat{\pmb{\theta}}_{-d}$ and
$\pmb{\theta}_{0,-d}$. Let $\Omega_{0}=(\varpi_{kr})_{(q+d-1)\times(q+d-1)}$
with
\[
\varpi_{kr}=E(V_{k}V_{r})-E[\pmb{B}(Z^{T}\pmb{\beta}_{0})V_{k}]^{T}\Gamma
^{-1}(\pmb{\beta}_{0},\pmb{\beta}_{0})E[\pmb{B}(Z^{T}\pmb{\beta}_{0}%
)V_{r}],\ \ k,r=1,\ldots,q,
\]%
\[
\varpi_{k(q+r)}=E[\dot{\pmb{B}}_{r}(Z^{T}\pmb{\beta}_{0})V_{k}]^{T}%
\pmb{b}_{0}-E[\pmb{B}(Z^{T}\pmb{\beta}_{0})V_{k}]^{T}\Gamma^{-1}%
(\pmb{\beta}_{0},\pmb{\beta}_{0})H_{r}(\pmb{\beta}_{0},\pmb{\beta}_{0}%
)\pmb{b}_{0},
\]
$\varpi_{(q+r)k}=\varpi_{k(q+r)}$ for $k,=1,\ldots,q;r=1,\ldots,d-1$, and
\[
\varpi_{(q+k)(q+r)}=\pmb{b}_{0}^{T}\left\{ R_{rk}(\pmb{\beta}_{0}%
,\pmb{\beta}_{0})-H_{r}^{T}(\pmb{\beta}_{0},\pmb{\beta}_{0})\Gamma
^{-1}(\pmb{\beta}_{0},\pmb{\beta}_{0})H_{k}(\pmb{\beta}_{0},\pmb{\beta}_{0}%
)\right\} \pmb{b}_{0}%
\]
for $k,r=1,\ldots,d-1$, where $\pmb{B}(Z^{T}\pmb{\beta})=(B_{1}(Z^{T}%
\pmb{\beta}),\ldots,B_{K_{n}}(Z^{T}\pmb{\beta}))^{T}$, and $\dot{\pmb{B}}%
_{r}(Z^{T}\pmb{\beta})=\frac{\partial\pmb{B}(Z^{T}\pmb{\beta})}{\partial
\beta_{r}}$, $H_{r}(\pmb{\beta},\pmb{\beta}^{\prime})$ and $R_{rk}%
(\pmb{\beta},\pmb{\beta}^{\prime})$ are $K_{n}\times K_{n}$ matrices whose
$(l,l^{\prime})$th elements are $E[B_{l}(Z^{T}\pmb{\beta})\dot{B}_{l^{\prime
}r}(Z^{T}\pmb{\beta}^{\prime})]$ and $E[\dot{B}_{lr}(Z^{T}\pmb{\beta})\dot
{B}_{l^{\prime}k}(Z^{T}\pmb{\beta}^{\prime})],$ respectively, and $\dot
{B}_{lr}(Z^{T}\pmb{\beta})=\frac{\partial B_{l}(Z^{T}\pmb{\beta})}%
{\partial\beta_{r}}$.

\textbf{Theorem\ 3.2.}\ Suppose that Assumptions 1 to 8 hold and that
$\Omega_{0}$ is invertible. Then we have
\[
\sqrt{n}\Omega_{0}^{1/2}(\hat{\pmb{\theta}}_{-d}-\pmb{\theta}_{0,-d}%
)\rightarrow_{d}N(0,\sigma^{2}I_{q+d-1}),\tag{3.9}
\]
where $I_{q+d-1}$ is the $(q+d-1)\times(q+d-1)$ identity matrix.

Next we establish the convergence rates of the estimators $\hat{a}(t)$ and
$\hat{g}(u)$.

\textbf{Theorem\ 3.3.}\ \ Assume that Assumptions 1 to 8 hold and that
$\tilde{m}\rightarrow\infty$, $n^{-1/2}\tilde{m}^{2}\lambda_{\tilde{m}}%
^{-1}\log\tilde{m}\rightarrow0$. Then
\[
\int_{\mathcal{T}}\left\{ \hat{a}(t)-a(t)\right\} ^{2}dt=O_{p}\Big(\frac
{\tilde{m}}{n\lambda_{\tilde{m}}}+\frac{\tilde{m}}{n^{2}\lambda_{\tilde{m}%
}^{2}}\sum_{j=1}^{\tilde{m}}\frac{j^{3}a_{j}^{2}}{\lambda_{j}^{2}}+\frac
{1}{n\lambda_{\tilde{m}}}\sum_{j=1}^{\tilde{m}}\frac{a_{j}^{2}}{\lambda_{j}%
}+\tilde{m}^{-2\gamma+1}\Big).\tag{3.10}
\]

If $\lambda_{j}\sim j^{-\delta}$, $\tilde{m}\sim n^{1/(\delta+2\gamma)}$,
$\gamma>2$ and $\gamma>1+\delta/2$, then $\sum_{j=1}^{\tilde{m}}j^{3}a_{j}%
^{2}\lambda_{j}^{-2}\leq C_{9}(\log\tilde{m}+\tilde{m}^{2\delta+4-2\gamma)})$
and $\sum_{j=1}^{\tilde{m}}a_{j}^{2}\lambda_{j}^{-1}<+\infty$, where $C_{9}$
is a positive constant. Then we have the following corollary.

\textbf{Corollary\ 3.1.}\ Under Assumptions 1 to 8, if $\lambda_{j}\sim
j^{-\delta}$, $\tilde{m}\sim n^{1/(\delta+2\gamma)}$ and $\gamma
>\min(2,1+\delta/2)$, then it follows that
\[
\int_{\mathcal{T}}\left\{ \hat{a}(t)-a(t)\right\} ^{2}dt=O_{p}\left(
n^{-(2\gamma-1)/(\delta+2\gamma)}\right) .\tag{3.11}
\]

The global convergence result (3.11) indicates that the estimator $\hat{a}(t)$
attains the same convergence rate as those of the estimators of Hall and
Horowitz \cite{r11},
which are optimal in the minimax sense.

From Theorem 3.2, we have $\Vert\hat{\pmb{\beta}}-\pmb{\beta}_{0}\Vert
=O_{p}(n^{-1/2})$. Then for sufficiently large $n$, $U_{\hat{\pmb{\beta}}}=
U_{\pmb{\beta}_{0}}$ and $U^{\hat{\pmb{\beta}}}= U^{\pmb{\beta}_{0}}$.

\textbf{Theorem\ 3.4}.\ Suppose that Assumptions 1 to 9 hold. Then,
\[
\int_{U_{\pmb{\beta}_{0}}}^{U^{\pmb{\beta}_{0}}}\left\{ \hat{g}%
(u)-g(u)\right\} ^{2}du=O_{p}\left( (nh)^{-1}+h^{2p}\right) .\tag{3.12}
\]
Further, if $h=O(n^{-1/(2p+1)})$ in Assumption 9, then
\[
\int_{U_{\pmb{\beta}_{0}}}^{U^{\pmb{\beta}_{0}}}\left\{ \hat{g}%
(u)-g(u)\right\} ^{2}du=O_{p}\left( n^{-2p/(2p+1)}\right) .\tag{3.13}
\]

\textbf{Remark\ 3.2.} Under Assumptions 1-8 and from a proof of similar to
that of Theorem 3.4, one can obtain
\[
\int_{U_{\pmb{\beta}_{0}}}^{U^{\pmb{\beta}_{0}}}\left\{ \tilde{g}%
(u)-g(u)\right\} ^{2}du=O_{p}\left( (nh_{0})^{-1}+h_{0}^{2p}\right)
=O_{p}\left( (nh_{0})^{-1}\right) .
\]
Due to the fact that $nh_{0}^{2p}\rightarrow0$, $\tilde{g}(u)$ does not attain
the global convergence rate of $O_{p}(n^{-2p/(2p+1)}),$ which is the optimal
rate for nonparametric models. In fact, the assumption that $nh_{0}%
^{2p}\rightarrow0$ is made in order to make the bias of the estimator
$\hat{\pmb{\beta}}_{-d}$ in Theorem 3.2 negligible. This results in slower
global convergence rate for the estimator $\tilde{g}(u)$.

Let $\mathcal{S}=\{(Y_{i},X_{i},W_{i},Z_{i}):i=1,\ldots,n\}$. If 
$(Y_{n+1},X_{n+1},W_{n+1},Z_{n+1})$ is a new vector of outcome and predictor
variables taken from the same population as that of the data $\mathcal{S}$ and
are independent of $\mathcal{S}$, then the \textit{mean squared}
\textit{prediction error} (MSPE) of $\hat{Y}_{n+1}$ is given by
\[%
\begin{array}
[c]{ll}%
\mbox{MSPE} & =E\Big[\Big\{\int_{\mathcal{T}}\hat{a}(t)X_{n+1}(t)dt+W_{n+1}%
^{T}\hat{\pmb{\alpha}}+\hat{g}(Z_{n+1}\hat{\pmb{\beta}})\\
& -\left( \int_{\mathcal{T}}a(t)X_{n+1}(t)dt+W_{n+1}^{T}\pmb{\alpha}_{0}%
+g(Z_{n+1}\pmb{\beta}_{0})\right) \Big\}^{2} \Big|\mathcal{S}\Big].
\end{array}
\]

\textbf{Theorem\ 3.5}. Under Assumptions 1 to 4 and 6 to 9, if $\lambda_{j}\sim
j^{-\delta}$, $\tilde{m}\sim n^{1/(\delta+2\gamma)}$, where $\gamma
>\min(2,1+\delta/2)$, $h_{0}\sim n^{-\tau}$ with $1/(2p)<\tau<(\gamma
-2)/(3(\delta+2\gamma))$ and $h=O(n^{-1/(2p+1)})$, then it follows that
\[
\mbox{MSPE}=O_{p}\left( n^{-(\delta+2\gamma-1)/(\delta+2\gamma)}%
)+O_{p}(n^{-2p/(2p+1)}\right) .\tag{3.14}
\]
Furthermore, if $\delta+2\gamma=2p+1$ then
\[
\mbox{MSPE}=O_{p}\left( n^{-(\delta+2\gamma-1)/(\delta+2\gamma)}\right)
.\tag{3.15}
\]

\textbf{Remark\ 3.3.} In Theorem 3.5, it is assumed that $h_{0}\sim n^{-\tau}$
and $1/(2p)<\tau<(\gamma-2)/(3(\delta+2\gamma))$. If $\delta+2\gamma=2p+1$,
then the conditions that $p>\gamma$ and $\gamma>5+3/(2p)$ are required. The
preceding conditions hold when $p>\gamma\geq5.3$.

\section{Simulation results}

In this section we present two Monte Carlo simulation studies to evaluate the
finite-sample performance of the proposed estimator. The data are generated
from the following models
\[
Y_{i}=\int_{\mathcal{T}}a(t)X_{i}(t)dt+\alpha_{0}W_{i}+\sin\left(  \pi
(Z_{i}^{T}\pmb{\beta}_{0}-E\right)  /(F-E))+\varepsilon_{i},\tag{4.1}
\]%
\[
Y_{i}=\int_{\mathcal{T}}a(t)X_{i}(t)dt+\alpha_{1}W_{i1}+\alpha_{2}%
W_{i2}-2Z_{i}^{T}\pmb{\beta}_{0}+5+\varepsilon_{i},\tag{4.2}
\]
with $\mathcal{T}=[0,1]$ and the trivariate random vectors $Z_{i}$'s have
independent components following the uniform distribution on $[0,1]$. In model
(4.1), $\alpha_{0}=0.3$, $\pmb{\beta}_{0}=(1,1,1)^{T}/\sqrt{3}$, $E=\sqrt
{3}/2-1.645/\sqrt{12}$ and $F=\sqrt{3}/2+1.645/\sqrt{12}$. We let $W_{i}=0$
for odd $i$ and $W_{i}=1$ for even $i$, and the $\varepsilon_{i}$'s are
independent errors following $N(0,0.5^{2})$. We take $a(t)=\sum_{j=1}%
^{50}a_{j}\phi_{j}(t)$ and $X_{i}(t)=\sum_{j=1}^{50}\xi_{ij}\phi_{j}(t)$,
where $a_{1}=0.3$ and $a_{j}=4(-1)^{j+1}j^{-2},\ j\geq2$; $\phi_{1}(t)\equiv1$
and $\phi_{j}(t)=2^{1/2}\cos((j-1)\pi t),\ j\geq2$; the $\xi_{ij}$'s are
independently and normally distributed with $N(0,j^{-\delta})$. In model
(4.2), $\alpha_{1}=-2$, $\alpha_{2}=1.5$, $\pmb{\beta}_{0}=(1,2,2)^{T}/3$ and
$X_{i}(t)=\sum_{j=1}^{50}\xi_{ij}\phi_{j}(t)$, the $\xi_{ij}$'s are
independently and normally distributed with $N(0,\lambda_{j})$, where
$\lambda_{1}=1$, $\lambda_{j}=0.22^{2}(1-0.0001j)^{2}$ if $2\leq j\leq4$,
$\lambda_{5j+k}=0.22^{2}((5j)^{-\delta/2}-0.0001k)^{2}$ for $j\geq1$ and
$0\leq k\leq4$. Further, $W_{ik}=\check{W}_{ik}+V_{ik}$ and $\check{W}%
_{ik}=\sum_{j=1}^{50}kj^{-2}\xi_{ij}$ for $k=1,2$. The $V_{ik}$'s are
independently and normally distributed with $N(-1,2^{2})$ and $N(2,3^{2}),$
respectively, and independent of $\xi_{ij}$. Finally, the error terms
$\varepsilon_{i}$'s in both (4.1) and (4.2) are independent $N(0,1)$ random variables.

For the functional linear part of model (4.1), the eigenvalues of the operator
$K$ are well-spaced, while the latter part of model (4.1) was investigated by
Carroll et al. \cite{r4}
and Yu and Ruppert \cite{r43}
In model (4.2), the
eigenvalues of the operator $K$ are closely spaced, while the link function
$g(u)=-2u+5$ is a linear function. All our results are reported based on the
average over 500 replications for each setting. In each sample, we first use a
linear function to replace $g(u)$ and use the least squares estimates for the
partial functional linear model as an initial estimator. The function $g(u)$
is approximated using a cubic spline with equally spaced knots. We note from
our simulation results (see Table 3) that parametric estimators are not
sensitive to the choices of parameters $h_{0}$ and $m$. Here we take
$h_{0}=c_{0}n^{-1/5}$ and $m=5,$ with $c_{0}=1.$ When we compute the
estimators of $g(u)$ and $a(t)$, the parameters $K_{n}$ and  $m$ are selected respectively by the BIC given in Section 2.

\begin{table}[hb]
\begin{center}
\caption{Results of Monte Carlo experiments for model
(4.2). The biases and sds  of parametric estimators and MISE
 of $\hat{g}(u)$ and MISE of $\hat{a}(t)$.}
\begin{tabular}
[c]{lccccccc}%
\hline
&  & \multicolumn{3}{c}{n=100} & \multicolumn{3}{c}{n=200}\\\hline
&  & LSPFL & ORACLE & PBS & LSPFL & ORACLE & PBS\\\hline
$\hat{\alpha}_{0}$ & bias & -0.0019 & 0.0034 & -0.0008 & -0.0025 & 0.0002 &
0.0002\\
& sd & 0.0836 & 0.0330 & 0.0307 & 0.0565 & 0.0159 & 0.0122\\
$\hat{\beta}_{01}$ & bias & -0.3678 & -0.0066 & -0.0056 & -0.3365 & -0.0037 &
0.0006\\
& sd & 0.5445 & 0.0441 & 0.0464 & 0.5141 & 0.0202 & 0.0206\\
$\hat{\beta}_{02}$ & bias & -0.3780 & -0.0075 & -0.0031 & -0.3283 & -0.0041 &
-0.0018\\
& sd & 0.5449 & 0.0457 & 0.0479 & 0.5201 & 0.0263 & 0.0178\\
$\hat{\beta}_{03}$ & bias & -0.0771 & 0.0082 & 0.0016 & -0.0553 & 0.0058 &
-0.0001\\
& sd & 0.2695 & 0.0506 & 0.0599 & 0.2694 & 0.0307 & 0.0239\\
$\hat{g}(u)$ & MISE &  &  & 0.0090 &  &  & 0.0007\\
$\hat{a}(t)$ & MISE & 0.1205 & 0.0189 & 0.0218 & 0.0756 & 0.0082 &
0.0084\\\hline
\end{tabular}
\end{center}
\end{table}

Table 1 reports the biases and standard deviations (sd) of the profile
B-spline (PBS) estimators $\hat{\alpha}_{0}$, $\hat{\pmb{\beta}}_{0}%
=(\hat{\beta}_{01},\hat{\beta}_{02},\hat{\beta}_{03})^{T}$ and the mean
integrated squared errors (MISE) of the estimators $\hat{g}(u)$ and $\hat
{a}(t)$ for model (4.1) based on $\delta=1.5$ and sample sizes $n=100$, $200$.
Figure 1 displays the true curves and the mean estimated curves over 500
simulations with sample size $n=100$ of $g(u)$, $a(t)$ and their $95\%$
pointwise confidence bands. Table 2 reports the biases and standard deviations
(sd) of the estimators $\hat{\alpha}_{k}$ for $k=1,2$ and $\hat{\pmb{\beta}}%
_{1}=(\hat{\beta}_{11},\hat{\beta}_{12},\hat{\beta}_{13})^{T}$, and the mean
integrated squared errors (MISE) of the estimators $\hat{g}(u)$ and $\hat
{a}(t)$ for model (4.2) with $\delta=1.5$ and $n=100,200$. For comparison
purposes, Tables 1 and 2 also list the simulation results based on the least
squares partial functional linear (LSPFL) estimators, which are obtained by
using a linear function to approximate the link function $g$. Further, Table 1
also lists the simulation results based on the nonlinear least squares
(ORACLE) estimation method when the exact form of sinusoidal model is known.


\begin{table}[h]
\begin{center}
\caption{Results of Monte Carlo experiments for model
(4.2). The biases ($\times10^{-4}$) and sds ($\times10^{-4}$) of parametric estimators and MISE
($\times10^{-4}$) of $\hat{g}(u)$ and MISE of $\hat{a}(t)$.}
\begin{tabular}
[c]{lccccc}%
\hline
&  & \multicolumn{2}{c}{n=100} & \multicolumn{2}{c}{n=200}  \\\hline
&  & LSPFL & PBS & LSPFL & PBS  \\\hline
$\hat{\alpha}_{1}$ & bias (sd) & 0.078(6.815) & 0.100(6.870) & 0.186(4.415) &
0.173(4.435)  \\
$\hat{\alpha}_{2}$ & bias (sd) & -0.071(4.612) & -0.085(4.666) &
0.359(3.038) & 0.373(3.056)  \\
$\hat{\beta}_{11}$ & bias (sd) & -0.707(22.725) & -0.753 (23.162) &
0.942(14.762) & 0.816(14.896)  \\
$\hat{\beta}_{12}$ & bias (sd) & -1.670(18.370) & -1.720(18.347) &
0.711(11.939) & 0.735(11.906)  \\
$\hat{\beta}_{13}$ & bias (sd) & 1.936(17.630) & 2.007(17.655) &
-1.220(11.944) & -1.181(11.919)  \\
$\hat{g}(u)$ & MISE &  & 3.852 &  & 2.503  \\
$\hat{a}(t)$ & MISE & 0.0087 & 0.0096 & 0.0047 & 0.0044  \\\hline
\end{tabular}
\end{center}
\end{table}
%

We observe from Table 1 that the least squares partial functional linear
(LSPFL) method gives poor estimates, while our profile B-spline estimates are
far more accurate than the LSPFL estimates, and they can be as accurate as
those obtained from the ORACLE when the exact form of sinusoidal model is
known. Figure 1 shows that the difference between the true curves and the mean
estimated curves are barely visible, and it shows that the bias is very small
in the estimates. Furthermore, the $95\%$ pointwise confidence bands are
reasonably close to the true curve, showing a very little variation in the
estimates. Table 2 shows that, even if the unknown link function $g(u)$ is a
linear function, our profile B-spline estimates behave as good as the least
squares partial functional linear estimates. Both tables indicate that the
proposed profile B-spline method yields accurate estimates and outperforms the
least squares partial functional linear estimates when the link function is
nonlinear, and it is comparable to the least squares partial functional linear
estimates when the link function is a linear function.

\begin{figure}[ptb]
\begin{center}
\resizebox{125mm}{55mm}{\includegraphics{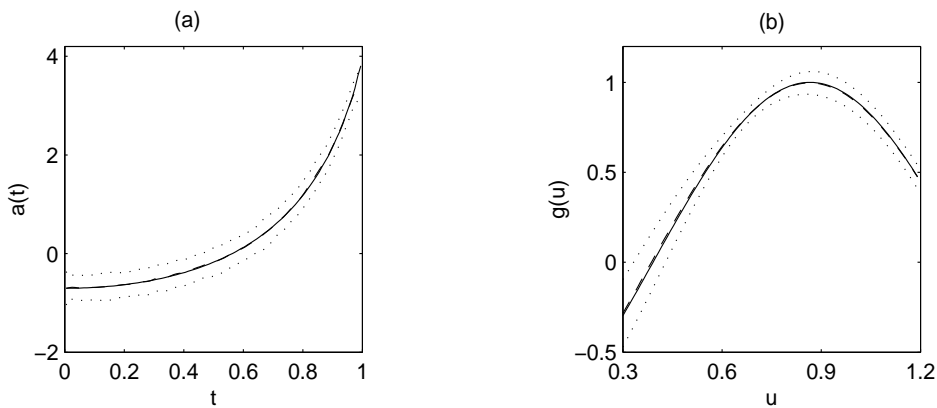}}
\caption{The actual and the mean estimated curves for $g(u)$
and $a(t)$ in model (4.1) with $n=100$ and the 95\% pointwise confidence
bands. (a) is the figure for $a(t)$ and (b) is the figure for $g(u)$. ---,
true curves; - - -, mean estimated curves; ..., 95\% pointwise confidence
bands.}
\end{center}
\end{figure}

To study the prediction performance of the proposed profile B-spline method,
we generated samples of $n=100,200$ from models (4.1) and (4.2) with
$\delta\in\{1.1,1.5,2\}$ for estimation, where $\delta$ is related to the
eigenvalue of the operator with kernel $K$. We also generated test samples of
size $300$ to compute the prediction mean absolute error (MAE) defined by
$MAE=\frac{1}{N}\sum_{i=1}^{N}|\tilde{Y}_{n+i}-\hat{Y}_{n+i}|$, where
$\tilde{Y}_{n+i}=\int_{\mathcal{T}}a(t)X_{n+i}(t)dt+W_{n+i}^{T}%
\pmb{\alpha}_{0}+g_{0}(Z_{n+i}^{T}\pmb{\beta}_{0})$ and $\hat{Y}_{n+i}%
=\int_{\mathcal{T}}\hat{a}(t)X_{n+i}(t)dt+W_{n+i}^{T}\hat{\pmb{\alpha}}%
+\hat{g}(Z_{n+i}^{T}\hat{\pmb{\beta}})$. Figures 2 and 3 display the boxplots
of $MAE$ based on 500 replications and $n=300$. We observe that the proposed
profile B-spline method shows good prediction performances for both models and
the MAEs are quite small even if $n=100$. 
Figure 2 also shows that the
$MAE$ decreases as $n$ increases or as $\delta$ increases. 

\begin{figure}[ht]
\begin{center}
\resizebox{125mm}{55mm}{\includegraphics{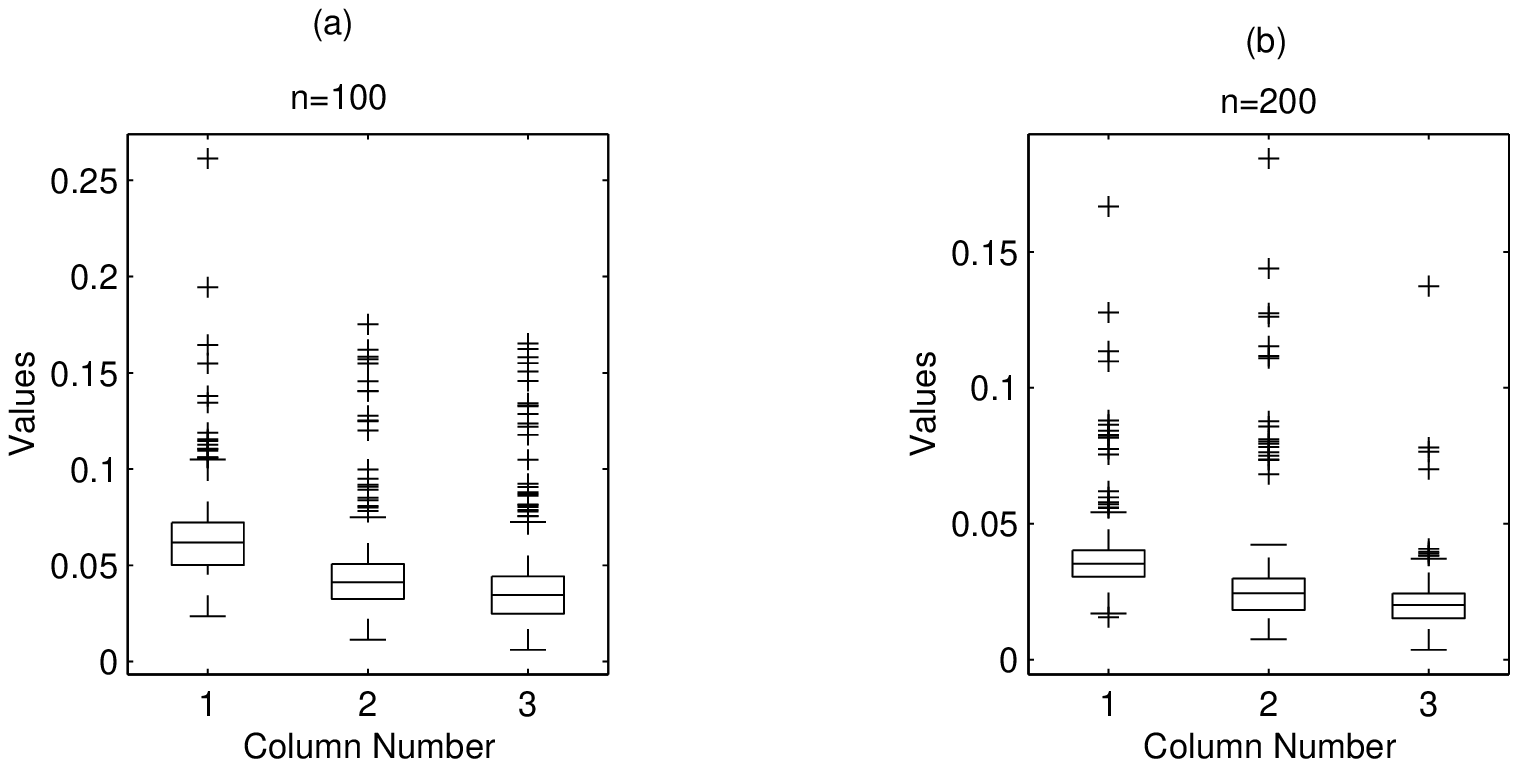}}
\resizebox{125mm}{55mm}{\includegraphics{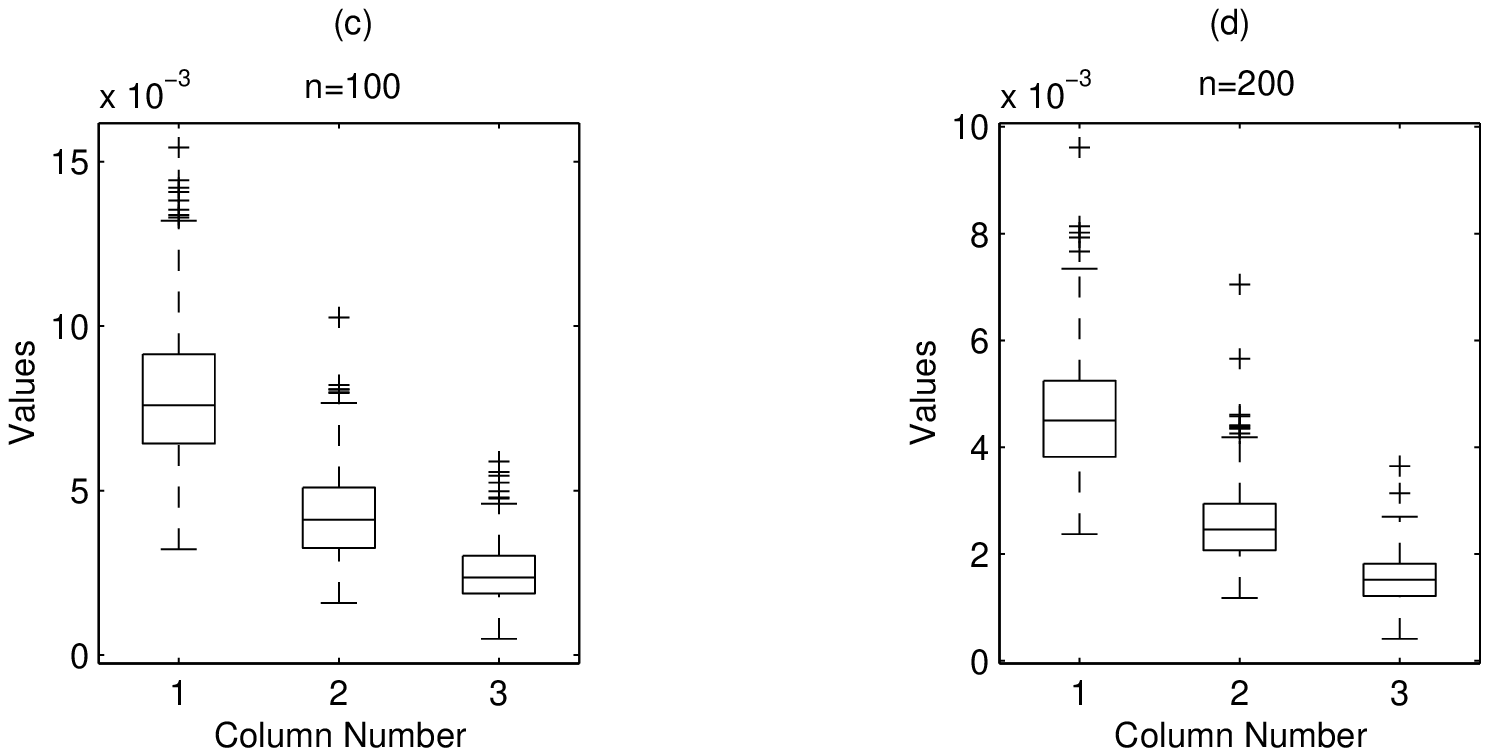}}
\caption{Boxplots of $MAE$ for models (4.1) ((a) and (b)) and (4.2) ((c) and (d)). Label 1 is
boxplot for $\delta=1.1$, 2 is boxplot for $\delta=1.5$ and 3 is boxplot for
$\delta=2$.}
\end{center}
\end{figure}


For different $m$ and $h_{0}$, Table 3 exhibits the MSEs of the estimators
$\hat{\alpha}_{0}$ and $\hat{\beta}_{01}$ for model (4.1) with $\delta=1.5$
and sample size $n=200$. We observe from Table 3 that MSEs of $\hat{\alpha
}_{0}$ and $\hat{\beta}_{01}$ are not very sensitive to the change of $m$ and
$h_{0},$ and the estimators of $\alpha_{0}$ and $\beta_{01}$ are efficient
under a broad range of values for $m$ and $h_{0}$. The MSEs of $\hat{\beta
}_{02}$ and $\hat{\beta}_{03}$ also show similar behaviors and are omitted here.

\begin{table}
\caption{MSE ($\times10^{-3}$) of $\hat{\alpha}_{0}$ and
$\hat{\beta}_{01}$ in model (4.1). The sample size is $n=200$.}
\begin{center}
{ 
\begin{tabular}
[c]{lcccccccccc}%
\hline
& $h_{0}$ & \multicolumn{9}{c}{$m$}\\\cline{3-11}
& & $2$ & $3$ & $4$ & $5$ & $6$ & $7$ & $8$ & $9$ &
$10$\\\hline
$\hat{\alpha}_{0}$ & 0.2 & 1.4 & 0.9 & 0.5 & 0.3 & 0.4 & 0.5 & 0.6 & 0.4 &
0.6\\
& 0.3 & 1.4 & 0.8 & 0.3 & 0.3 & 0.4 & 0.4 & 0.6 & 0.2 & 0.3\\
& 0.4 & 1.4 & 0.8 & 0.3 & 0.2 & 0.4 & 0.3 & 0.6 & 0.4 & 0.6\\
& 0.5 & 1.4 & 0.6 & 0.3 & 0.3 & 0.3 & 0.4 & 0.6 & 0.6 & 0.4\\
$\hat{\beta}_{01}$ & 0.2 & 1.8 & 2.0 & 1.0 & 0.5 & 0.9 & 1.1 & 0.8 & 0.9 &
1.6\\
& 0.3 & 1.5 & 1.5 & 0.7 & 0.6 & 1.4 & 0.3 & 1.1 & 0.7 & 1.4\\
& 0.4 & 1.5 & 1.5 & 0.7 & 0.6 & 1.1 & 0.3 & 1.1 & 1.1 & 1.4\\
& 0.5 & 2.5 & 1.6 & 1.5 & 0.8 & 0.6 & 1.0 & 1.6 & 1.2 & 1.0\\\hline
\end{tabular}
}
\end{center}
\end{table}

\section{Real data application}

In this section we analyze a real data set using the proposed method. For this
purpose we use the diffusion tensor imaging (DTI) data with 217 subjects from
the NIH Alzheimer's Disease Neuroimaging Initiative (ADNI) study. For more
information on how this data were collected etc., see
http://www.adni-info.org. The DTI data were processed by two key steps
including a weighted least squares estimation method \cite{r1,r46} 
to construct the diffusion tensors and a TBSS pipeline in FSL \cite{r31}
to register DTIs from multiple subjects to create a mean
image and a mean skeleton. This data have been recently analyzed by many
authors using different models; see, e.g., Yu et al. \cite{r42},
 Li et al. \cite{r17}
and the references therein.

Our interest is to predict mini-mental state examination (MMSE) scores, one of
the most widely used screening tests to provide brief and objective measures
of cognitive functioning for a long time. The MMSE scores have been seen as a
reliable and valid clinical measure quantitatively assessing the severity of
cognitive impairment. It was believed that the MMSE scores to be affected by
demographic features such as age, education and cultural background \cite{r35}
gender \cite{r25,r23},
and possibly some genetic factors, for example, AOPE polymorphic alleles \cite{r20}.

After cleaning the raw data that failed in quality control or had missing
data, we include totally 196 individuals in our analysis. The response of
interest $Y$ is the MMSE scores. The functional covariate is fractional
anisotropy (FA) values along the corpus callosum (CC) fiber tract with 83
equally spaced grid points, which can be treated as a function of arc-length.
FA measures the inhomogeneous extent of local barriers to water diffusion and
the averaged magnitude of local water diffusion \cite{r2}. 
The
scalar covariates of primary interests include gender ($W_{1}$), handedness
($W_{2}$), education level ($W_{3}$), genotypes for apoe4 ($W_{4},W_{5}$,
categorical data with 3 levels), age ($W_{6}$), ADAS13 ($Z_{1}$) and ADAS11
($Z_{2}$). The genotypes apoe4 is one of three major alleles of apolipoprotein
E (Apo-E), a major cholesterol carrier that supports lipid transport and
injury repair in the brain. APOE polymorphic alleles are the main genetic
determinants of Alzheimer disease risk \cite{r20}.
ADAS11 and ADAS13
are respectively the 11-item and 13-item versions of the Alzheimer\~{O}s
Disease Assessment Scale-Cognitive subscale (ADAS-Cog), which were originally
developed to measure cognition in patients within various stages of
Alzheimer's Disease \cite{r19,r48,r24}.

We study the following two models
\[%
\begin{array}
[c]{ll}%
Y & =\int_{0}^{1}a(t)X(t)dt+\alpha_{0}+\alpha_{1}W_{1}+\alpha_{2}W_{2}%
+\alpha_{3}W_{3}+\alpha_{4}W_{4}+\alpha_{5}W_{5}\\
& +\alpha_{6}W_{6}+\beta_{1}Z_{1}+\beta_{2}Z_{2}+\varepsilon,
\end{array}
\tag{5.1}
\]%
\[%
\begin{array}
[c]{ll}%
Y & =\int_{0}^{1}a(t)X(t)dt+\alpha_{1}W_{1}+\alpha_{2}W_{2}+\alpha_{3}%
W_{3}+\alpha_{4}W_{4}+\alpha_{5}W_{5}\\
& +\alpha_{6}W_{6}+g(\beta_{1}Z_{1}+\beta_{2}Z_{2})+\varepsilon,
\end{array}
\tag{5.2}
\]
where $W_{1}=1$ stands for male and $W_{1}=0$ stands for female, $W_{2}=1$
denotes right-handed and $W_{2}=0$ denotes left-handed, $W_{4}=1$ and
$W_{5}=0$ indicates type 0 for apoe4, $W_{4}=0$ and $W_{5}=1$ indicates type 1
for apoe4 and both $W_{4}=0$ and $W_{5}=0$ indicates type 2 for apoe4. The
functional component $X(t)$ is chosen as the centered fractional anisotropy
(FA) values so that $E[X(t)]=0$. Model (5.1) is a partial functional linear
model, while model (5.2) is partial functional linear single index model in
which ADAS13 ($Z_{1}$) and ADAS11 ($Z_{2}$) are index variables.

\begin{table}[h]
\caption{The parametric estimators for models (5.1) and
(5.2).}
\begin{center}
{ 
\begin{tabular}
[c]{lcccccccc}%
\hline
model & $\alpha_{1}$ & $\alpha_{2}$ & $\alpha_{3}$ & $\alpha_{4}$ &
$\alpha_{5}$ & $\alpha_{6}$ & $\beta_{1}$ & $\beta_{2}$\\\hline
(5.1) & 0.0758 & 0.4317 & 0.1105 & 0.6875 & 0.5581 & -0.0239 & -0.0429 &
-0.1865\\
(5.2) & -0.0754 & 0.1814 & 0.1138 & 0.5961 & 0.5245 & -0.0305 & 0.1957 &
0.9807\\\hline
\end{tabular}
}
\end{center}
\end{table}

The parametric and nonparametric components in the models are computed by the
procedure given in Section 2, with the nonparametric function $g(u)$ being
approximated by a cubic spline with equally spaced knots. Since the values of
$Z_{1}$ and $Z_{2}$ are large, we choose $h_{0}=5.0$ for model (5.2) and $m=3$
for parametric estimation. Table 4 exhibits the parametric estimators, and
Figure 4 shows the estimated curves of $a(t)$ and $g(u)$. For model (5.1),
$\hat{a}_{0}=28.9388$. The MSE of $Y$ for models (5.1) and (5.2) are 2.8684
and 2.7782, respectively, and can be further reduced for model (5.2) as the
number of knots increases.

From Table 4 and Figure 3, we observe that in both models MMSE is decreasing
in terms of ADAS13 and ADAS11. However, in Figure 3 this decline is found to
be nonlinear evidenced by the nonlinear trends of $g(u)$ in model (5.2). In
single index models (5.2), we found that MMSE is higher for female than male,
which is consistent with the results in the literature \cite{r25,r23},
while model (5.1) incorrectly finds the opposite.
Although we may not able to perform a formal test on model fitting, these
observations show the superiority of the single index model (5.2).

\begin{figure}[ptb]
\begin{center}

\resizebox{125mm}{50mm}{\includegraphics{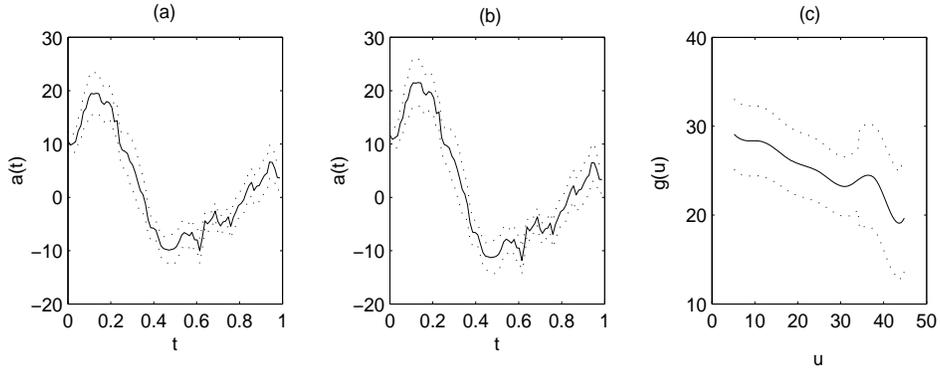}}
\caption{The solid lines are the estimated curves of $a(t)$
in model (5.1) (a), $a(t)$ in model (5.2) (b), and $g(u)$ in model (5.2) (c);
The doted lines are their corresponding $95\%$ point-wise confidence
intervals.}
\end{center}
\end{figure}


To evaluate the prediction performance of the three models, we applied a
combination of the bootstrap and the cross-validation method to the data set.
For each bootstrap sample, we randomly divided the data into ten partitions.
Since the number of individuals is not large, we used nine folds of the data
to estimate the model and the remaining fold for the testing data set. We
calculated the mean squared prediction error (MSPE) for the testing data set.
The MSPEs for the two models over the 200 replications are reported as
boxplots in Figure 4. The means for MSPEs of the 200 replications for models
(5.1) and (5.2) are $3.6996$ and $3.4249$, respectively. The medians for MSPEs
of the 200 replications for models (5.1) and (5.2) are 3.5464 and 3.3421,
respectively. This figure shows that model (5.2) fits the data better than
model (5.1). We also calculated $95\%$ point-wise confidence intervals of the
estimated curves of $a(t)$ in model (5.1), $a(t)$ in model (5.2), and $g(u)$
in model (5.2), which are shown as (a), (b) and (c), respectively, in Figure
4. From Figure 5, it is evidenced that the functional slope for both models
are also identical, while $g(u)$ has a clear nonlinear feature. This also
confirms that model (5.2) is more flexible than model (5.1).

\begin{figure}[ptb]
\begin{center}
\resizebox{60mm}{50mm}{\includegraphics{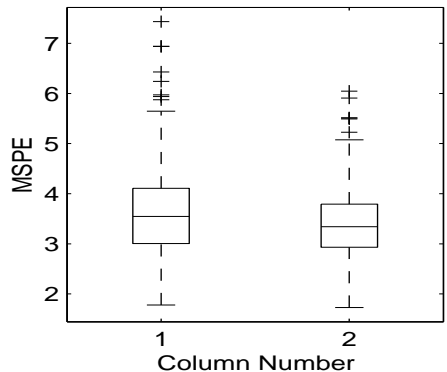}}
\caption{Boxplots for the mean squared prediction error
(MSPE) for three models. Label 1 is boxplot for model (5.1) and 2 is boxplot
for model (5.2).}
\end{center}
\end{figure}

\section{Summarizing remarks}

Functional data analysis is now very popular as it provides modern analytical
tools for data that are recorded as images or as a continuous phenomenon over
a period of time. Classical multivariate statistical tools may fail or may be
irrelevant in that context to take benefit from the underlying functional
structure of functional data. As a great variety of real data applications
involve functional phenomena, which may be represented as curves or more
complex objects, the demand of models and statistical tools for analyzing
functional data is ever more increasing.

The need for more comprehensive models that more adaptable motivated us to
propose and study a partial functional partially linear single index (PFPLSI)
model in this paper. The proposed PFPLSI model generalizes the standard
functional linear model, partial functional linear models and the partially
linear single-index model, among others. We have implemented functional
principal component analysis to estimate the slope function component of the
PFPLSI model, and the unknown link function of the single-index component has
been approximated by a B-spline function. To estimate the unknown parameters
in the proposed PFPLSI model, we have proposed a profile B-spline method. We
have derived the asymptotic properties, including the consistency and
asymptotic normality, of the proposed estimators of the unknown parameters.
The global convergence of the proposed estimator of the functional slope
function has also been established, and this convergence result has been shown
to be optimal in the minimax sense. A two-stage procedure was used to estimate
the unknown link function attaining the optimal global convergence rate of
convergence. We have also derived convergence rates of the mean squared
prediction error for a predictor. The lower prediction error demonstrates the
rationality of our modelling and the effectiveness of the proposed estimation
procedure. Monte Carlo studies conducted to examine the performance of the
proposed methodology demonstrate that the proposed estimators perform quite
satisfactorily and the theoretical results established seem to be valid.

An alternative approach to the PFPLSI model (1.1) that may be of interest is
\textit{functional linear quantile regression}. The functional linear quantile
regression where the conditional quantiles of the responses are modeled by a
set of scalar covariates and functional covariates. There may be several
advantages of using conditional quantiles instead of working with conditional
means. First, the quantile regression, in particular the median regression,
provides an alternative and complement to the mean regression, while being
resistant to outliers in the responses. In other words, it is more efficient
than the mean regression when the errors follow a distribution with heavy
tails. Second, the quantile regression is capable of dealing with
heteroscedasticity, that is the situations where variances depend on some
covariates. Finally, the quantile regression can give a more complete picture
on how the responses are affected by covariates; e.g., some tail behaviors of
the responses conditional on the covariates. For more details on quantile
regression, one may refer to the monograph of Koenker \cite{r16}. 
In view of the
model (1.1), we consider the following functional linear quantile regression:
for given $\tau\in(0,1)$,%
\[
Q_{\tau}(y|X,Z,W)=\int_{\mathcal{T}}a_{\tau}(t)X(t)dt+W^{T}\pmb{\alpha}_{0\tau
}+g(Z^{T}\pmb{\beta}_{0\tau})
\]
where $Q_{\tau}(y|X,Z,W)$ is the $\tau$-th conditional quantile of $Y$ given
the covariates $(X,Z,W).$ Although there is some reported work on functional
linear quantile regression in the literature, the above model has not been
studied yet. Further research is needed for these advancements.

\section*{Appendix: Proofs}

In this section we let $C>0$ denote a generic constant of which the value may
change from line to line. For a matrix $A=(a_{ij})$, set $\Vert A\Vert
_{\infty}=\max_{i}\sum_{j}|a_{ij}|$ and $|A|_{\infty}=\max_{i,j}|a_{ij}|$. For
a vector $v=(v_{1},\ldots,v_{k})^{T}$, set $\Vert v\Vert_{\infty}=\sum
_{j=1}^{k}|v_{j}|$ and $|v|_{\infty}=\max_{1\leq j\leq k}|v_{j}|$. We write
$Y_{i}=Y_{i}^{\ast}+\varepsilon_{i}$ with $Y_{i}^{\ast}=\int_{\mathcal{T}%
}a(t)X_{i}(t)dt+W_{i}^{T}\pmb{\alpha}_{0}+g(Z_{i}^{T}\pmb{\beta}_{0})$. Denote
$\check{Y}_{i}=Y_{i}^{\ast}-\frac{1}{n}\sum_{l=1}^{n}Y_{l}^{\ast}\tilde{\xi
}_{il}$, $\tilde{\varepsilon}_{i}=\varepsilon_{i}-\frac{1}{n}\sum_{l=1}%
^{n}\varepsilon_{l}\tilde{\xi}_{il}$ and $\check{\pmb{Y}}=(\check{Y}%
_{1},\ldots,\check{Y}_{n})^{T}$, $\tilde{\pmb{\varepsilon}}=(\tilde
{\varepsilon}_{1},\ldots,\tilde{\varepsilon}_{n})^{T}$. Then $\tilde{Y}%
_{i}=\check{Y}_{i}+\tilde{\varepsilon}_{i}$ and $\tilde{\pmb{Y}}%
=\check{\pmb{Y}}+\tilde{\pmb{\varepsilon}}$. Define
$\pmb{P}(\pmb{\beta})=I_{n}-\tilde{\pmb{B}}(\pmb{\beta})(\tilde{\pmb{B}}%
^{T}(\pmb{\beta})\tilde{\pmb{B}}(\pmb{\beta}))^{-1}\tilde{\pmb{B}}%
^{T}(\pmb{\beta})$, where $I_{n}$ is the $n\times n$ identity matrix. By
(3.5), (2.9) and (2.10), we have
\[
G_{n}(\pmb{\alpha},\pmb{\beta})=\frac{1}{n}[(\check{\pmb{Y}}-\tilde
{\pmb{W}}\pmb{\alpha})^{T}\pmb{P}(\pmb{\beta} )(\check{\pmb{Y}}-\tilde
{\pmb{W}}\pmb{\alpha})+2(\check{\pmb{Y}}-\tilde{\pmb{W}}\pmb{\alpha})^{T}%
\pmb{P}(\pmb{\beta} )\tilde{\pmb{\varepsilon}}+\tilde{\pmb{\varepsilon}}%
^{T}\pmb{P}(\pmb{\beta})\tilde{\pmb{\varepsilon} }].\tag{A.1}
\]

\textbf{Lemma\ A.1.} \ Suppose that Assumptions 1 to 4, 5' and 7 hold. Then
\[%
\begin{array}
[c]{l}%
\frac{1}{n}(\check{\pmb{Y}}-\tilde{\pmb{W}}\pmb{\alpha})^{T}(\check
{\pmb{Y}}-\tilde{\pmb{W}}\pmb{\alpha} )=\rho(\pmb{\alpha})+o_{p}(1),
\end{array}
\]
where $\rho(\pmb{\alpha})=(\pmb{\alpha}-\pmb{\alpha}_{0})^{T}E(VV^{T}%
)(\pmb{\alpha}-\pmb{\alpha}_{0})-2\pmb{b}_{0}^{T}E[\pmb{B}_{\pmb{\beta}_{0}%
}(Z^{T}\pmb{\beta}_{0})V^{T}](\pmb{\alpha}-\pmb{\alpha}_{0})+\pmb{b}_{0}%
^{T}\Gamma(\pmb{\beta}_{0},\pmb{\beta}_{0})\pmb{b}_{0},$ and $o_{p}(1)$ holds
uniformly for $\pmb{\alpha}$ in any bounded neighborhood of $\pmb{\alpha}_{0}$.

\textbf{Proof.} \ Define $\check{\xi}_{il}=\sum_{j=1}^{m}\frac{\xi_{lj}%
\xi_{ij}}{\lambda_{j}}$, $\check{Y}_{i1}=Y_{i}^{\ast}-\frac{1}{n}\sum
_{l=1}^{n}Y_{l}^{\ast}\check{\xi}_{il}$ and $\check{Y}_{i2}=\frac{1}{n}%
\sum_{l=1}^{n}Y_{l}^{\ast}(\tilde{\xi}_{il}-\check{\xi}_{il}).$ Then
$\check{Y}_{i}=\check{Y}_{i1}-\check{Y}_{i2}$ and
\[
\frac{1}{n}\check{\pmb{Y}}^{T}\check{\pmb{Y}}=\frac{1}{n}\sum_{i=1}^{n}%
(\check{Y}_{i1}^{2}-2\check{Y}_{i1}\check{Y}_{i2}+\check{Y}_{i2}%
^{2}).\tag{A.2}
\]
Denote $\check{Y}_{i21}=\sum_{j=1}^{m}\frac{1}{\lambda_{j}}[\frac{1}{n}%
\sum_{l=1}^{n}Y_{l}^{\ast}(\hat{\xi}_{lj}-\xi_{lj})]\xi_{ij}$, $\check
{Y}_{i22}=\sum_{j=1}^{m}(\frac{1}{\hat{\lambda}_{j}}-\frac{1}{\lambda_{j}%
})$ $(\frac{1}{n}\sum_{l=1}^{n}Y_{l}^{\ast}\hat{\xi}_{lj})\xi_{ij}$ and
$\check{Y}_{i23}=\sum_{j=1}^{m}\frac{1}{\hat{\lambda}_{j}}(\frac{1}{n}%
\sum_{l=1}^{n}Y_{l}^{\ast}\hat{\xi}_{lj})(\hat{\xi}_{ij}-\xi_{ij}).$ Then we
have
\[
\check{Y}_{i2}^{2}\leq3(\check{Y}_{i21}^{2}+\check{Y}_{i22}^{2}+\check
{Y}_{i23}^{2}).\tag{A.3}
\]
From Lemma 5.1 of Hall and Horowitz (2007) it follows that
\[
\hat{\xi}_{lj}-\xi_{lj}=\sum_{k\neq j}\frac{\xi_{lk}}{\hat{\lambda}%
_{j}-\lambda_{k}}\int\Delta\hat{\phi}_{j}\phi_{k}+\xi_{lj}\int(\hat{\phi}%
_{j}-\phi_{j})\phi_{j},\tag{A.4}
\]
where $\Delta=\hat{K}-K$. Then we obtain
\[%
\begin{array}
[c]{ll}%
\lbrack\frac{1}{n}\sum_{l=1}^{n}Y_{l}^{\ast}(\hat{\xi}_{lj}-\xi_{lj})]^{2} 
&\leq2(\sum_{k\neq j}\frac{\vec{\xi}_{k}}{\hat{\lambda}_{j}-\lambda_{k}}%
\int\Delta\hat{\phi}_{j}\phi_{k})^{2}+2(\vec{\xi}_{j}\int(\hat{\phi}_{j}%
-\phi_{j})\phi_{j})^{2}\\
 &\leq2[\sum_{k\neq j}\frac{\vec{\xi}_{k}^{2}}{(\hat{\lambda}_{j}-\lambda
_{k})^{2}}][\sum_{k=1}^{\infty}(\int\Delta\hat{\phi}_{j}\phi_{k})^{2}%
]+\\&2\vec{\xi}_{j}^{2}(\int(\hat{\phi}_{j}-\phi_{j})\phi_{j})^{2},
\end{array}
\]
where $\vec{\xi}_{j}=\frac{1}{n}\sum_{l=1}^{n}Y_{l}^{\ast}\xi_{lj}$. Lemma 6.1
of Cardot et al. (2007) yields that
\[
|\lambda_{j}-\lambda_{k}|\geq\lambda_{j}-\lambda_{j+1}\geq\lambda_{m}%
-\lambda_{m+1}\geq\lambda_{m}/(m+1)\geq\lambda_{m}/(2m)
\]
uniformly for $1\leq j\leq m$. From (5.2) of Hall and Horowitz (2007) we have
$\sup_{j\geq1}|\hat{\lambda}_{j}-\lambda_{j}|\leq|\Vert\Delta\Vert
|=O_{p}(n^{-1/2})$ and
\[%
\begin{array}
[c]{ll}%
(\int(\hat{\phi}_{j}-\phi_{j})\phi_{j})^{2}\leq\Vert\hat{\phi}_{j}-\phi
_{j}\Vert^{2}\leq C\frac{|\Vert\Delta\Vert|^{2}}{(\lambda_{j}-\lambda
_{j+1})^{2}}\leq C|\Vert\Delta\Vert|^{2}\lambda_{j}^{-2}j^{2}, &
\end{array}
\tag{A.5}
\]
where $|\Vert\Delta\Vert|=(\int_{\mathcal{T}}\int_{\mathcal{T}}\Delta
^{2}(s,t)dsdt)^{1/2}$. Using Parseval's identity, we obtain
\[
\sum_{k=1}^{\infty}(\int\Delta\hat{\phi}_{j}\phi_{k})^{2}=\int(\int\Delta
\hat{\phi}_{j})^{2}\leq|\Vert\Delta\Vert|^{2}=O_{p}(n^{-1}).
\]
Assumption 5' implies that $|\hat{\lambda}_{j}-\lambda_{j}|=o_{p}(\lambda
_{m}/m)$. Consequently, $\sum_{k\neq j}\frac{\vec{\xi}_{k}^{2}}{(\hat{\lambda
}_{j}-\lambda_{k})^{2}}$ $=\sum_{k\neq j}\frac{\vec{\xi}_{k}^{2}}{(\lambda
_{j}-\lambda_{k})^{2}}[1+o_{p}(1)]$, where $o_{p}(1)$ holds uniformly for
$1\leq j\leq m$. Using Lemma 6.2 of Cardot et al. (2007) and the fact that
$(\lambda_{j}-\lambda_{k})^{2}\geq(\lambda_{k}-\lambda_{k+1})^{2}$, we deduce
that
\[%
\begin{array}
[c]{l}%
\sum_{k\neq j}\frac{1}{(\lambda_{j}-\lambda_{k})^{2}}E(\vec{\xi}_{k}^{2})\\
\leq C\sum_{k\neq j}\frac{1}{(\lambda_{j}-\lambda_{k})^{2}}[n^{-1}\lambda
_{k}+{a_{k}^{\ast}}^{2}\lambda_{k}^{2}]\\
\leq C[\frac{1}{n(\lambda_{j}-\lambda_{j+1})}\sum_{k\neq j}\frac{\lambda_{k}%
}{|\lambda_{j}-\lambda_{k}|}+\sum_{k=1}^{j-1}\frac{\lambda_{k}^{2}{a_{k}%
^{\ast}}^{2}}{(\lambda_{k}-\lambda_{k+1})^{2}}+\\\sum_{k=j+1}^{2j}\frac
{j^{2}{a_{k}^{\ast}}^{2}}{(k-j)^{2}}+\sum_{k=2j+1}^{\infty}\frac{\lambda
_{k}^{2}{a_{k}^{\ast}}^{2}}{(\lambda_{j}-\lambda_{2j})^{2}}]\\
\leq C(n^{-1}\lambda_{j}^{-1}j^{2}\log j+1).
\end{array}
\]
where $a_{k}^{\ast}=a_{k}+\sum_{r=1}^{q}w_{rk}\alpha_{0r}$. Assumption 2
yields that
\[
\sum_{j=1}^{m}\lambda_{j}^{-2}j^{2}\log j\leq m^{-2}\lambda_{m}^{-2}\sum
_{j=1}^{m}j^{4}\log j\leq\lambda_{m}^{-2}m^{3}\log m
\]
and $\sum_{j=1}^{m}\lambda_{j}^{-1}\leq\lambda_{m}^{-1}m$. Therefore,
\[%
\begin{array}
[c]{ll}%
\frac{1}{n}\sum_{i=1}^{n}\check{Y}_{i21}^{2} & \leq(\sum_{j=1}^{m}\frac
{1}{\lambda_{j}}[\frac{1}{n}\sum_{l=1}^{n}Y_{l}^{\ast}(\hat{\xi}_{lj}-\xi
_{lj})]^{2})(\sum_{j=1}^{m}\frac{1}{n\lambda_{j}}\sum_{i=1}^{n}\xi_{ij}^{2})\\
& =O_{p}(n^{-2}\lambda_{m}^{-2}m^{4}\log m+n^{-1}\lambda_{m}^{-1}m^{2}).
\end{array}
\tag{A.6}
\]
Decomposing $\frac{1}{n}\sum_{l=1}^{n}Y_{l}^{\ast}\hat{\xi}_{lj}=\vec{\xi}%
_{j}+\frac{1}{n}\sum_{l=1}^{n}Y_{l}^{\ast}(\hat{\xi}_{lj}-\xi_{lj})$ and using
(A.6), we obtain
\[%
\begin{array}
[c]{ll}%
\frac{1}{n}\sum_{i=1}^{n}\check{Y}_{i22}^{2} & \leq C\sum_{j=1}^{m}\frac
{(\hat{\lambda}_{j}-\lambda_{j})^{2}}{\lambda_{j}^{3}}(\frac{1}{n}\sum
_{l=1}^{n}Y_{l}^{\ast}\hat{\xi}_{lj})^{2}[1+o_{p}(1)]\\&(\sum_{j=1}^{m}\frac
{1}{n\lambda_{j}}\sum_{i=1}^{n}\xi_{ij}^{2})\\
& =O_{p}(n^{-1}\lambda_{m}^{-1}m+n^{-3}\lambda_{m}^{-4}m^{4}\log
m+n^{-2}\lambda_{m}^{-3}m^{2}).
\end{array}
\tag{A.7}
\]
By (A.10) of Tang (2015), it holds that
\[
\Vert\hat{\phi}_{j}-\phi_{j}\Vert^{2}=O_{p}(n^{-1}j^{2}\log j)\tag{A.8}
\]
uniformly for $1\leq j\leq m$. Using (A.7) and (A.8), we obtain
\[%
\begin{array}
[c]{ll}%
\frac{1}{n}\sum_{i=1}^{n}\check{Y}_{i23}^{2} & \leq(\sum_{j=1}^{m}\frac
{1}{\hat{\lambda}^{2}}(\frac{1}{n}\sum_{l=1}^{n}Y_{l}^{\ast}\hat{\xi}%
_{lj})^{2})(\frac{1}{n}\sum_{i=1}^{n}\Vert X_{i}\Vert^{2})(\sum_{j=1}^{m}%
\Vert\hat{\phi}_{j}-\phi_{j}\Vert^{2})\\
& =O_{p}((n^{-1}m^{3}+n^{-3}\lambda_{m}^{-3}m^{6}\log m+n^{-2}\lambda_{m}%
^{-2}m^{4})\log m).
\end{array}
\tag{A.9}
\]
Then by (A.3), (A.6), (A.7), (A.9) and Assumption 5', we conclude that
\[
\frac{1}{n}\sum_{i=1}^{n}\check{Y}_{i2}^{2}=O_{p}(n^{-2}\lambda_{m}^{-2}%
m^{4}\log m+n^{-1}\lambda_{m}^{-1}m^{2})=o_{p}(h_{0}^{2}).\tag{A.10}
\]
Define $\xi_{j}^{\ast}=\frac{1}{n}\sum_{l=1}^{n}\lambda_{j}^{-1/2}\xi
_{lj}Y_{l}^{\ast}$. Since $E[\max_{1\leq j\leq m}(\xi_{j}^{\ast}-E(\xi
_{j}^{\ast}))^{2}]\leq$ \\ $\frac{1}{n}\sum_{j=1}^{m}E(\xi_{j}Y^{\ast})^{2}\leq
Cn^{-1}$, we then have $\max_{1\leq j\leq m}|\xi_{j}^{\ast}-E(\xi_{j}^{\ast
})|=O_{p}(n^{-1/2})$. Hence, we have
\[%
\begin{array}
[c]{ll}%
\frac{1}{n}\sum_{i=1}^{n}\check{Y}_{i1}^{2} & =\frac{1}{n}\sum_{i=1}^{n}%
{Y_{i}^{\ast}}^{2}-2\sum_{j=1}^{m}{\xi_{j}^{\ast}}^{2}+\sum_{j=1}^{m}%
\frac{{\xi_{j}^{\ast}}^{2}}{n\lambda_{j}}(\sum_{i=1}^{n}\xi_{ij}^{2}%
)+\\&\sum_{j\neq j^{\prime}}\xi_{j}^{\ast}\xi_{j^{\prime}}^{\ast}\bar{\xi
}_{jj^{\prime}}\\
& =\sum_{j=1}^{\infty}(a_{j}+\sum_{r=1}^{q}w_{rj}\alpha_{0r})^{2}\lambda
_{j}+E(V^{T}\pmb{\alpha}_{0}+g(Z^{T}\pmb{\beta}_{0}))^{2}\\
& -2\sum_{j=1}^{m}(a_{j}+\sum_{r=1}^{q}w_{rj}\alpha_{0r})^{2}\lambda_{j}%
+\\&\sum_{j=1}^{m}(a_{j}+\sum_{r=1}^{q}w_{rj}\alpha_{0r})^{2}\lambda_{j}%
+o_{p}(1)\\
& =E(V^{T}\pmb{\alpha}_{0}+g(Z^{T}\pmb{\beta}_{0}))^{2}+o_{p}(1),
\end{array}
\tag{A.11}
\]
where $\bar{\xi}_{jj^{\prime}}=\frac{1}{n(\lambda_{j}\lambda_{j^{\prime}%
})^{1/2}}\sum_{i=1}^{n}\xi_{ij}\xi_{ij^{\prime}}$. Combining (A.2), (A.10),
(A.11) and (3.1), we conclude that
\[
\frac{1}{n}\check{\pmb{Y}}^{T}\check{\pmb{Y}}=\pmb{\alpha}_{0}^{T}%
E(VV^{T})\pmb{\alpha}_{0}+2\pmb{b}_{0}^{T}E[\pmb{B}_{\pmb{\beta}_{0}}%
(Z^{T}\pmb{\beta}_{0})V^{T}]\pmb{\alpha}_{0}+\pmb{b}_{0}^{T}\Gamma
(\pmb{\beta}_{0},\pmb{\beta}_{0})\pmb{b}_{0}+o_{p}(1).\tag{A.12}
\]
Similar to the proof of (A.12), we obtain that
\[%
\begin{array}
[c]{l}%
\frac{1}{n}\tilde{\pmb{W}}^{T}\tilde{\pmb{W}}=E(VV^{T})+o_{p}(1),\\ \frac
{1}{n}\tilde{\pmb{Y}}^{T}\tilde{\pmb{W}}=\pmb{\alpha}_{0}E(V^{T}%
V)+\pmb{b}_{0}^{T}E[\pmb{B}_{\pmb{\beta}_{0}}(Z^{T}\pmb{\beta}_{0}%
)V]+o_{p}(1).
\end{array}
\]
Now Lemma A.1 follows from (A.12) and the preceding expression.

\textbf{Lemma\ A.2.} \ Under Assumptions 1, 4 and 5', it holds that
\[
\begin{array}
[c]{l}%
\sup_{\pmb{\beta}\in\Theta_{\rho_{0}}}\max_{1\leq j\leq m}\max_{1\leq k\leq
K_{\pmb{\beta}}}\lambda_{j}^{-\frac{1}{2}}|\frac{1}{n}\sum_{i=1}^{n}\xi
_{ij}B_{k\pmb{\beta}}^{(r)}(Z_{i}^{T}\pmb{\beta})|\\=o_{p}(n^{-\frac{1}{2}}%
h_{0}^{\frac{1}{4}-r}\log n),
\end{array}
\]%
\[
\begin{array}
[c]{l}%
\sup_{\pmb{\beta}\in\Theta_{\rho_{0}}}\max_{k,k^{\prime}}|\frac{1}{n}%
\sum_{i=1}^{n}B_{k\pmb{\beta}}(Z_{i}^{T}\pmb{\beta})B_{k^{\prime}%
\pmb{\beta}}(Z_{i}^{T}\pmb{\beta})-E[B_{k\pmb{\beta} }(Z_{i}^{T}%
\pmb{\beta})B_{k^{\prime}\pmb{\beta}}(Z_{i}^{T}\pmb{\beta})]|\\=o_{p}%
(n^{-\frac{1}{2}}h_{0}^{\frac{1}{2}}\log n),
\end{array}
\]%
\[
\begin{array}
[c]{l}%
\sup_{\pmb{\beta}\in\Theta_{\rho_{0}}}\max_{k,k^{\prime}}|\frac{1}{n}%
\sum_{i=1}^{n}B_{k\pmb{\beta}}^{\prime}(Z_{i}^{T}\pmb{\beta})B_{k^{\prime
}\pmb{\beta}}^{\prime}(Z_{i}^{T}\pmb{\beta})-E[B_{k\pmb{\beta}}^{\prime}%
(Z_{i}^{T}\pmb{\beta})B_{k^{\prime}\pmb{\beta}}^{\prime}(Z_{i}^{T}%
\pmb{\beta})]|\\=o_{p}(n^{-\frac{1}{2}}h_{0}^{-\frac{3}{2}}\log n),
\end{array}
\]
and
\[
\begin{array}
[c]{l}%
\sup_{\pmb{\beta}\in\Theta_{\rho_{0}}}\max_{k,k^{\prime}}|\frac{1}{n}%
\sum_{i=1}^{n}B_{k\pmb{\beta}}(Z_{i}^{T}\pmb{\beta})B_{k^{\prime}%
\pmb{\beta}}^{\prime\prime}(Z_{i}^{T}\pmb{\beta})-E[B_{k\pmb{\beta}}(Z_{i}%
^{T}\pmb{\beta})B_{k^{\prime}\pmb{\beta}}^{\prime\prime}(Z_{i}^{T}%
\pmb{\beta})]|\\=o_{p}(n^{-\frac{1}{2}}h_{0}^{-\frac{3}{2}}\log n)
\end{array}
\]
for $r=0,1,2$.

\textbf{Proof.} \ We give only the proof for the first step with $r=2$, as the
first step with $r=0,1$ and the other steps follow from similar arguments.
Define $\eta_{jki}(Z_{i}^{T}\pmb{\beta})=\lambda_{j}^{-1/2}\xi_{ij}%
B_{k\pmb{\beta}}^{\prime\prime}(Z_{i}^{T}\pmb{\beta})$. Applying Assumptions 1
and Lemma 5 of Kato \cite{r15},
we have $\max_{1\leq j\leq m,1\leq i\leq
n}|\lambda_{j}^{-1/2}\xi_{ij}|=O_{p}((mn)^{1/4})$. Hence, by Assumption 5',
for any $\varepsilon>0$ and $\epsilon>0$, there exists a positive constant
$\tilde{C}_{1}$ such that
\[
P\{\max_{1\leq j\leq m,1\leq i\leq n}|\lambda_{j}^{-1/2}\xi_{ij}|\geq\tilde
{C}_{1}n^{1/2}h_{0}^{1/4}(\log n)^{-1}\}<\epsilon/4.\tag{A.13}
\]
Using Assumptions 1 and the fact that $|B_{k\pmb{\beta}}^{\prime\prime}%
(Z_{i}^{T}\pmb{\beta})|\leq Ch_{0}^{-2}$, we obtain
\[%
\begin{array}
[c]{l}%
|E[\lambda_{j}^{-\frac{1}{2}}\xi_{ij}B_{k\pmb{\beta}}^{\prime\prime}(Z_{i}%
^{T}\pmb{\beta})I_{\{|\lambda_{j}^{-\frac{1}{2}}\xi_{ij}|\geq\tilde{C}%
_{1}n^{\frac{1}{2}}h_{0}^{\frac{1}{4}}(\log n)^{-1}\}}]|\\
\leq Cn^{-\frac{3}{2}}h_{0}^{-\frac{11}{4}}(\log n)^{3}E[\lambda_{j}%
^{-\frac{1}{2}}\xi_{ij}]^{4}<\varepsilon n^{-\frac{1}{2}}h_{0}^{-\frac{7}{4}%
}\log n/2.
\end{array}
\]
Denote
\[%
\begin{array}
[c]{ll}%
\tilde{\eta}_{jki}(Z_{i}^{T}\pmb{\beta}) & =\lambda_{j}^{-\frac{1}{2}}\xi
_{ij}B_{k\pmb{\beta}}^{\prime\prime}(Z_{i}^{T}\pmb{\beta})I_{\{|\lambda
_{j}^{-\frac{1}{2}}\xi_{ij}|<\tilde{C}_{1}n^{\frac{1}{2}}h_{0}^{\frac{1}{4}%
}(\log n)^{-1}\}}\\
& -E[\lambda_{j}^{-\frac{1}{2}}\xi_{ij}B_{k\pmb{\beta}}^{\prime\prime}%
(Z_{i}^{T}\pmb{\beta})I_{\{|\lambda_{j}^{-\frac{1}{2}}\xi_{ij}|<\tilde{C}%
_{1}n^{\frac{1}{2}}h_{0}^{\frac{1}{4}}(\log n)^{-1}\}}].
\end{array}
\]
Then we have
\[%
\begin{array}
[c]{l}%
P\{\sup_{\pmb{\beta}\in\Theta_{\rho_{0}}}\max_{j,k}|\frac{1}{n}\sum_{i=1}%
^{n}\eta_{jki}(Z_{i}^{T}\pmb{\beta})|\geq\varepsilon n^{-\frac{1}{2}}%
h_{0}^{-\frac{7}{4}}\log n\}\\
\leq P\{\max_{j,i}|\lambda_{j}^{-\frac{1}{2}}\xi_{ij}|\geq\tilde{C}%
_{1}n^{\frac{1}{2}}h_{0}^{\frac{1}{4}}(\log n)^{-1}\}\\
+P\{\sup_{\pmb{\beta}\in\Theta_{\rho_{0}}}\max_{j,k}|\frac{1}{n}\sum_{i=1}%
^{n}\tilde{\eta}_{jki}(Z_{i}^{T}\pmb{\beta})|\geq\varepsilon n^{-\frac{1}{2}%
}h_{0}^{-\frac{7}{4}}\log n/2\}.
\end{array}
\tag{A.14}
\]
Using the fact that $|B_{k\pmb{\beta}}^{\prime\prime}(Z_{i}^{T}%
\pmb{\beta})|\leq Ch_{0}^{-2},$ again we obtain
\[
|\tilde{\eta}_{jki}(Z_{i}^{T}\pmb{\beta})|\leq Cn^{\frac{1}{2}}h_{0}%
^{-\frac{7}{4}}(\log n)^{-1}.\tag{A.15}
\]
From Assumption 1, it follows that
\[
\sum_{i=1}^{n}E(\tilde{\eta}_{jki}^{2}(Z_{i}^{T}\pmb{\beta}))\leq
Cn\lambda_{j}^{-1}(E[B_{k\pmb{\beta}}^{\prime\prime4}(Z_{i}^{T}%
\pmb{\beta})]E(\xi_{j}^{4}))^{1/2}\leq Cnh_{0}^{-7/2}.\tag{A.16}
\]
For $\pmb{\beta}_{1}=(\beta_{11},\ldots,\beta_{1d})^{T}\in\Theta_{\rho_{0}}$
and $\pmb{\beta}_{2}=(\beta_{21},\ldots,\beta_{2d})^{T}\in\Theta_{\rho_{0}}$,
define $|\pmb{\beta}_{2}-\pmb{\beta}_{1}|=\max_{1\leq r\leq d-1}|\beta
_{2r}-\beta_{1r}|$. Since $\sum_{j=1}^{m}\frac{1}{n}\sum_{i=1}^{n}\lambda
_{j}^{-\frac{1}{2}}|\xi_{ij}|=O_{p}(m)$, then there exists a positive
$\tilde{C}_{2}$ such that
\[
P\{\sum_{j=1}^{m}\frac{1}{n}\sum_{i=1}^{n}\lambda_{j}^{-\frac{1}{2}}|\xi
_{ij}|\geq\tilde{C}_{2}m\}<\epsilon/4.\tag{A.17}
\]
From (2.6), for all $\pmb{\beta}\in\Theta_{\rho_{0}}$, the total of different
$B_{k\pmb{\beta}}(u)$ is not more than $(s+1)k_{n}$. Let $\Theta_{\rho_{0}}$
be divided into $N$ disjoint parts $\Theta_{\rho_{0}1},\cdots,\Theta_{\rho
_{0}N}$ such that for any $\pmb{\beta}\in\Theta_{\rho_{0}l},1\leq l\leq N$ and
any $1\leq j\leq m,1\leq k\leq(s+1)k_{n}$, when $\sum_{j=1}^{m}\frac{1}{n}%
\sum_{i=1}^{n}\lambda_{j}^{-\frac{1}{2}}|\xi_{ij}|<\tilde{C}_{2}m$,
\[%
\begin{array}
[c]{l}%
\sup_{\pmb{\beta}\in\Theta_{\rho_{0}l}}|\frac{1}{n}\sum_{i=1}^{n}\tilde{\eta
}_{jki}(Z_{i}^{T}\pmb{\beta})-\frac{1}{n}\sum_{i=1}^{n}\tilde{\eta}%
_{jki}(Z_{i}^{T}\pmb{\beta}_{l})|\\
\leq\sup_{\pmb{\beta}\in\Theta_{\rho_{0}l}}\lambda_{j}^{-\frac{1}{2}%
}\Big(\frac{1}{n}\sum_{i=1}^{n}|\xi_{ij}||B_{k\pmb{\beta}}^{\prime\prime
}(Z_{i}^{T}\pmb{\beta})-B_{k\pmb{\beta}}^{\prime\prime}(Z_{i}^{T}%
\pmb{\beta}_{l})|\\+E(|\xi_{ij}||B_{k\pmb{\beta}}^{\prime\prime}(Z_{i}%
^{T}\pmb{\beta})-B_{k\pmb{\beta}}^{\prime\prime}(Z_{i}^{T}\pmb{\beta}_{l}%
)|)\Big)\\
\leq\sup_{\pmb{\beta}\in\Theta_{\rho_{0}l}}Ch_{0}^{-3}\sum_{j=1}^{m}%
\Big(\frac{1}{n}\sum_{i=1}^{n}\lambda_{j}^{-\frac{1}{2}}|\xi_{ij}%
|+E(\lambda_{j}^{-\frac{1}{2}}|\xi_{ij}|)\Big)|\pmb{\beta}-\pmb{\beta}_{l}|\\
\leq Cmh_{0}^{-3}|\pmb{\beta}-\pmb{\beta}_{l}|<\varepsilon n^{-\frac{1}{2}%
}h_{0}^{-\frac{7}{4}}\log n/4.
\end{array}
\]
This can be done with $N=C(mn^{1/2}/(\varepsilon h_{0}^{5/4}\log n))^{d-1}$.
Using Bernstein inequality and (A.15), (A.16) and Assumption 5', for
sufficiently large $n$, it follows that
\[%
\begin{array}
[c]{l}%
P\Big(\sup_{\pmb{\beta}\in\Theta_{\rho_{0}}}\max_{j,k}|\frac{1}{n}\sum
_{i=1}^{n}\tilde{\eta}_{jki}(Z_{i}^{T}\pmb{\beta})|\geq\varepsilon
n^{-\frac{1}{2}}h_{0}^{-\frac{7}{4}}\log n/2,\\\sum_{j=1}^{m}\frac{1}{n}%
\sum_{i=1}^{n}\lambda_{j}^{-\frac{1}{2}}|\xi_{ij}|<\tilde{C}_{2}m\Big)\\
\leq P\Big(\cup_{l=1}^{N}\{\max_{j,k}|\frac{1}{n}\sum_{i=1}^{n}\tilde{\eta
}_{jki}(Z_{i}^{T}\pmb{\beta}_{l})|\geq\varepsilon n^{-\frac{1}{2}}%
h_{0}^{-\frac{7}{4}}\log n/4\}\Big)\\
\leq Cmk_{n}N\exp\Big\{-\frac{\varepsilon^{2}nh_{0}^{-\frac{7}{2}}(\log
n)^{2}}{32Cnh_{0}^{-7/2}+4Cn^{\frac{1}{2}}h_{0}^{-\frac{7}{4}}(\log
n)^{-1}\varepsilon n^{\frac{1}{2}}h_{0}^{-\frac{7}{4}}\log n}\Big\}<\epsilon
/2.
\end{array}
\]
Now Lemma A.2 follows from (A.13), (A.14), (A.17) and the preceding inequality.

\textbf{Lemma\ A.3.} \ Assume that Assumptions 1, 2, 4 and 5' hold. Then it
holds that
\[
\frac{1}{n}\tilde{\pmb{B}}^{T}(\pmb{\beta})\tilde{\pmb{B}}(\pmb{\beta})=\Gamma
(\pmb{\beta},\pmb{\beta} )+o_{p}(h_{0}^{2}),
\]
where $o_{p}(h_{0}^{2})$ holds uniformly for $1\leq k,k^{\prime}\leq
K_{\pmb{\beta} }$ and $\pmb{\beta}\in\Theta_{\rho_{0}}$.

\textbf{Proof}.\ Define
\[
\begin{array}
[c]{l}%
\tilde{B}_{k\pmb{\beta}1}(Z_{i}^{T}\pmb{\beta})=B_{k\pmb{\beta}}(Z_{i}%
^{T}\pmb{\beta})-\frac{1}{n}\sum_{l=1}^{n}B_{k\pmb{\beta}}(Z_{l}%
^{T}\pmb{\beta})\check{\xi}_{il},\\ \tilde{B}_{k\pmb{\beta}2}(Z_{i}%
^{T}\pmb{\beta})=\frac{1}{n}\sum_{l=1}^{n}B_{k\pmb{\beta}}(Z_{l}%
^{T}\pmb{\beta})(\tilde{\xi}_{il}-\check{\xi}_{il}).
\end{array}
\]
We decompose the $(k,k^{\prime})$th element of $\frac{1}{n}\tilde{\pmb{B}}%
^{T}(\pmb{\beta})\tilde{\pmb{B}}(\pmb{\beta})$ as
\[%
\begin{array}
[c]{ll}%
&\frac{1}{n}\sum_{i=1}^{n}\tilde{B}_{k\pmb{\beta}}(Z_{i}^{T}\pmb{\beta})\tilde
{B}_{k^{\prime}\pmb{\beta}}(Z_{i}^{T}\pmb{\beta}) \\& =\frac{1}{n}\sum_{i=1}%
^{n}\Big(\tilde{B}_{k\pmb{\beta}1}(Z_{i}^{T}\pmb{\beta})\tilde{B}_{k^{\prime
}\pmb{\beta}1}(Z_{i}^{T}\pmb{\beta})-\tilde{B}_{k\pmb{\beta}1}(Z_{i}%
^{T}\pmb{\beta})\tilde{B}_{k^{\prime}\pmb{\beta}2}(Z_{i}^{T}\pmb{\beta})\\
& -\tilde{B}_{k\pmb{\beta}2}(Z_{i}^{T}\pmb{\beta})\tilde{B}_{k^{\prime
}\pmb{\beta}1}(Z_{i}^{T}\pmb{\beta})+\tilde{B}_{k\pmb{\beta}2}(Z_{i}%
^{T}\pmb{\beta})\tilde{B}_{k^{\prime}\pmb{\beta}2}(Z_{i}^{T}\pmb{\beta})\Big).
\end{array}
\]
Applying the Cauchy-Schwarz inequality, Lemma A.2, (A.8) and Assumptions 2 and
5', we obtain
\[%
\begin{array}
[c]{l}%
\sup_{\pmb{\beta}\in\Theta_{\rho_{0}}}\max_{k}\frac{1}{n}\sum_{i=1}%
^{n}\Big(\sum_{j=1}^{m}\frac{1}{\lambda_{j}}[\frac{1}{n}\sum_{l=1}%
^{n}B_{k\pmb{\beta}}(Z_{l}^{T}\pmb{\beta})(\hat{\xi}_{lj}-\xi_{lj})]\xi
_{ij}\Big)^{2}\\
\leq\sup_{\pmb{\beta}\in\Theta_{\rho_{0}}}\max_{k}\Big(\sum_{j=1}^{m}\frac
{1}{\lambda_{j}}[\frac{1}{n}\sum_{l=1}^{n}B_{k\pmb{\beta}}(Z_{l}%
^{T}\pmb{\beta})(\hat{\xi}_{lj}-\xi_{lj})]^{2}\Big)\\\Big(\sum_{j=1}^{m}\frac
{1}{n\lambda_{j}}\sum_{i=1}^{n}\xi_{ij}^{2}\Big)\\
\leq(\sup_{\pmb{\beta}\in\Theta_{\rho_{0}}}\max_{k}\frac{1}{n}\sum_{l=1}%
^{n}B_{k\pmb{\beta}}^{2}(Z_{l}^{T}\pmb{\beta}))(\frac{1}{n}\sum_{l=1}^{n}\Vert
X_{l}\Vert^{2})(\sum_{j=1}^{m}\frac{\Vert\hat{\phi}_{j}-\phi_{j}\Vert^{2}%
}{\lambda_{j}})\\(\sum_{j=1}^{m}\frac{1}{n\lambda_{j}}\sum_{i=1}^{n}\xi_{ij}%
^{2})\\
=O_{p}(n^{-1}\lambda_{m}^{-1}m^{4}h_{0}\log m)=o_{p}(h_{0}^{3}).
\end{array}
\tag{A.18}
\]
Similar to the proof of (A.7), (A.9) and using Lemma A.2, we then deduce that
\[%
\begin{array}
[c]{l}%
\sup_{\pmb{\beta}\in\Theta_{\rho_{0}}}\max_{k}\frac{1}{n}\sum_{i=1}%
^{n}\Big(\sum_{j=1}^{m}(\frac{\hat{\xi}_{ij}}{\hat{\lambda}_{j}}-\frac
{\xi_{ij}}{\lambda_{j}})(\frac{1}{n}\sum_{l=1}^{n}B_{k\pmb{\beta}}(Z_{l}%
^{T}\pmb{\beta})\hat{\xi}_{lj})\Big)^{2}\\
=o_{p}(n^{-2}\lambda_{m}^{-2}mh_{0}^{1/2}(\log n)^{2})+o_{p}(n^{-2}\lambda
_{m}^{-1}m^{3}h_{0}^{1/2}(\log n)^{2})\\+O_{p}(n^{-2}\lambda_{m}^{-3}m^{4}%
h_{0}\log m)
+O_{p}(n^{-2}\lambda_{m}^{-2}m^{6}h_{0}(\log m)^{2})\\=o_{p}(h_{0}^{3}).
\end{array}
\tag{A.19}
\]
Using Lemma A.2 and Assumption 5', we conclude that
\[%
\begin{array}
[c]{ll}%
&\frac{1}{n}\sum_{i=1}^{n}\tilde{B}_{k\pmb{\beta}1}(Z_{i}^{T}\pmb{\beta})\tilde
{B}_{k^{\prime}\pmb{\beta}1}(Z_{i}^{T}\pmb{\beta}) \\& =\frac{1}{n}\sum
_{i=1}^{n}B_{k\pmb{\beta}}(Z_{i}^{T}\pmb{\beta})B_{k^{\prime}\pmb{\beta}}%
(Z_{i}^{T}\pmb{\beta})-2\sum_{j=1}^{m}\rho_{kj}\rho_{k^{\prime}j}\\
& +\sum_{j=1}^{m}\rho_{kj}\rho_{k^{\prime}j}(\frac{1}{n\lambda_{j}}\sum
_{i=1}^{n}\xi_{ij}^{2})+\sum_{j\neq j^{\prime}}\rho_{kj}\rho_{k^{\prime
}j^{\prime}}\bar{\xi}_{jj^{\prime}}\\
& =E[B_{k\pmb{\beta}}(Z^{T}\pmb{\beta})B_{k^{\prime}\pmb{\beta}}%
(Z^{T}\pmb{\beta})]+o_{p}(h_{0}^{2}),
\end{array}
\]
where $\rho_{kj}=\frac{1}{n\lambda_{j}^{1/2}}\sum_{l=1}^{n}\xi_{lj}%
B_{k\pmb{\beta}}(Z_{l}^{T}\pmb{\beta})$. Now Lemma A.3 follows from (A.18),
(A.19) and the preceding equation.

\textbf{Proof of Theorem\ 3.1.} By arguments similar to those used in the proof
of Lemmas A.1 and A.3, it follows that
\[
\frac{1}{n}\tilde{\pmb{B}}^{T}(\pmb{\beta})(\check{\pmb{Y}}-\tilde
{\pmb{W}}\pmb{\alpha})=\Pi(\pmb{\alpha},\pmb{\beta})+o_{p}(h_{0}).\tag{A.20}
\]
Using Lemma A.3, (A.20) and arguments similar to those used in the proof of
Lemma 1 of Tang (2013), we then deduce that
\[
\begin{array}
[c]{l}%
\frac{1}{n}(\check{\pmb{Y}}-\tilde{\pmb{W}}\pmb{\alpha})^{T}\tilde
{\pmb{B}}(\pmb{\beta})(\tilde{\pmb{B}}^{T}(\pmb{\beta})\tilde{\pmb{B}}%
(\pmb{\beta}))^{-1}\tilde{\pmb{B}}^{T}(\pmb{\beta})(\check{\pmb{Y}}%
-\tilde{\pmb{W}}\pmb{\alpha})\\=\Pi^{T}(\pmb{\alpha},\pmb{\beta})\Gamma
^{-1}(\pmb{\beta},\pmb{\beta})\Pi(\pmb{\alpha},\pmb{\beta})+o_{p}%
(1).\tag{A.21}
\end{array}
\]
Therefore, Lemma A.1 and (A.21) imply that
\[%
\begin{array}
[c]{ll}%
\frac{1}{n}(\check{\pmb{Y}}-\tilde{\pmb{W}}\pmb{\alpha})^{T}%
\pmb{P}(\pmb{\beta})(\check{\pmb{Y}}-\tilde{\pmb{W}}\pmb{\alpha}) &
=\rho(\pmb{\alpha})-\Pi^{T}(\pmb{\alpha},\pmb{\beta})\Gamma^{-1}%
(\pmb{\beta},\pmb{\beta})\Pi(\pmb{\alpha},\pmb{\beta})+o_{p}(1)\\
& =:\tilde{G}(\pmb{\alpha},\pmb{\beta})+o_{p}(1),
\end{array}
\tag{A.22}
\]
where $o_{p}(1)$ holds uniformly for $\pmb{\beta}\in\Theta_{\rho_{0}}$ and
$\pmb{\alpha}$ is in any bounded neighborhood of $\pmb{\alpha}_{0}$. Similar
to the proof of Lemmas A.1 and A.3, it holds that $\frac{1}{n}\tilde
{\pmb{\varepsilon}}^{T}\tilde{\pmb{\varepsilon}}=\sigma^{2}+o_{p}(1)$,
$\frac{1}{n}(\check{\pmb{Y}}-\tilde{\pmb{W}}\pmb{\alpha})^{T}\tilde
{\pmb{\varepsilon}}=o_{p}(h_{0})$ and $\frac{1}{n}\tilde{\pmb{B}}%
^{T}(\pmb{\beta})\tilde{\pmb{\varepsilon}}=o_{p}(h_{0})$. Similar to the proof
of (A.21) and (A.22), we further have $\frac{1}{n}(\check{\pmb{Y}}%
-\tilde{\pmb{W}}\pmb{\alpha})^{T}\pmb{P}(\pmb{\beta})\tilde{\pmb{\varepsilon}}%
=o_{p}(1)$ and $\frac{1}{n}\tilde{\pmb{\varepsilon}}^{T}%
\pmb{P}(\pmb{\beta})\tilde{\pmb{\varepsilon}}=\sigma^{2}+o_{p}(1)$. Therefore,
from (A.1), (A.22) and (3.2), it follows that
\[
G_{n}(\pmb{\alpha},\pmb{\beta})=G(\pmb{\alpha},\pmb{\beta})+o_{p}%
(1),\tag{A.23}
\]
where $o_{p}(1)$ holds uniformly for $\pmb{\beta}\in\Theta_{\rho_{0}}$ and
$\pmb{\alpha}$ is in any bounded neighborhood of $\pmb{\alpha}_{0}$. By the
fact that $(\hat{\pmb{\alpha}},\hat{\pmb{\beta}})$ is the minimizer of
$G_{n}(\pmb{\alpha},\pmb{\beta})$ and using (A.23), we have
\[
G_{n}(\hat{\pmb{\alpha}},\hat{\pmb{\beta}})\leq G_{n}(\pmb{\alpha}_{0}%
,\pmb{\beta}_{0})=G(\pmb{\alpha}_{0},\pmb{\beta}_{0})+o_{p}(1).\tag{A.24}
\]
By (A.1) and (A.22), we have that $\tilde{G}(\pmb{\alpha},\pmb{\beta})\geq0$
and $G(\pmb{\alpha},\pmb{\beta})\geq\sigma^{2}$. From (3.2), one obtains
$G(\pmb{\alpha}_{0},\pmb{\beta}_{0})=\sigma^{2}+o_{p}(1)$. Applying (A.23) and
(A.24), we obtain that $\sigma^{2}\leq G(\hat{\pmb{\alpha}},\hat
{\pmb{\beta}})=G_{n}(\hat{\pmb{\alpha}},\hat{\pmb{\beta}})+o_{p}(1)\leq
G(\pmb{\alpha}_{0},\pmb{\beta}_{0})+o_{p}(1).$ Therefore, $|G(\hat
{\pmb{\alpha}},\hat{\pmb{\beta}})-G(\pmb{\alpha}_{0},\pmb{\beta}_{0}%
)|=o_{p}(1)$; that is, $|G^{\ast}(\hat{\pmb{\theta}}_{-d})-G^{\ast
}(\pmb{\theta}_{0,-d})|=o_{p}(1)$. Since $G^{\ast}(\pmb{\theta}_{-d})$ is
locally convex at $\pmb{\theta}_{0,-d}$, it follows that $\hat{\pmb{\alpha}}%
-\pmb{\alpha}_{0}=o_{p}(1)$ and $\hat{\pmb{\beta}}_{-d}-\pmb{\beta}_{0,-d}%
=o_{p}(1)$. This completes the proof of (3.3).

From (A.11), Assumption 5 and the fact that $\lambda_{j}\leq C/(j\log j)$, we
have
\[
\sum_{j=m+1}^{\infty}(a_{j}+\sum_{r=1}^{q}w_{rj}\alpha_{0r})^{2}\lambda
_{j}\leq Cm^{-2\gamma}=o(h_{0}^{2}).\tag{A.25}
\]
Applying Assumption 5 and (A.25), we can easily prove that $\frac{1}{n}%
(\check{\pmb{Y}}-\tilde{\pmb{W}}\pmb{\alpha})^{T}(\check{\pmb{Y}}%
-\tilde{\pmb{W}}\pmb{\alpha})=\rho(\pmb{\alpha} )+o_{p}(h_{0}^{2})$ in Lemma
A.1, $\frac{1}{n}\tilde{\pmb{B}}^{T}(\pmb{\beta})\tilde{\pmb{B}}%
(\pmb{\beta})=\Gamma(\pmb{\beta},\pmb{\beta})+o_{p}(h_{0}^{4})$ in Lemma A.3
and $\frac{1}{n}\tilde{\pmb{B}}^{T}(\pmb{\beta})(\check{\pmb{Y}}%
-\tilde{\pmb{W}}\pmb{\alpha})=\Pi(\pmb{\alpha},\pmb{\beta} )+o_{p}(h_{0}^{3}%
)$. Consequently, it follows that $G_{n}(\pmb{\alpha}
,\pmb{\beta})=G(\pmb{\alpha},\pmb{\beta})+o_{p}(h_{0}^{2})$ and $|G(\hat
{\pmb{\alpha}},\hat{\pmb{\beta} })-G(\pmb{\alpha}_{0},\pmb{\beta}_{0}%
)|=o_{p}(h_{0}^{2})$. Now (3.4) follows from Assumption 8. This completes the
proof of Theorem 3.1.

\textbf{Lemma\ A.4.} \ Under Assumptions 1-7, it holds that
\[
\ddot{G}_{n}(\pmb{\theta}_{-d},\tilde{\pmb{b}}(\pmb{\theta}_{-d}%
))=2\Omega(\pmb{\theta}_{-d})+o_{p}(1),
\]
where $o_{p}(1)$ holds uniformly for $\pmb{\beta}\in\Theta_{\rho_{0}}$,
$\pmb{\alpha}$ is in any bounded neighborhood of $\pmb{\alpha}_{0}$ and
$\Omega(\pmb{\beta}_{-d})=(\pi_{kr})_{(q+d-1)\times(q+d-1)}$ with
\[
\pi_{kr}=E(V_{k}V_{r})-E[\pmb{B}(Z^{T}\pmb{\beta})V_{k}]^{T}\Gamma
^{-1}(\pmb{\beta} ,\pmb{\beta})E[\pmb{B}(Z^{T}\pmb{\beta})V_{r}%
],\ \ k,r=1,\ldots,q,\tag{A.26}
\]%
\[
\pi_{k(q+r)}=E[\dot{\pmb{B}}_{r}(Z^{T}\pmb{\beta})V_{k}]^{T}\bar
{\pmb{b}}(\pmb{\alpha},\pmb{\beta} )+E[\pmb{B}(Z^{T}\pmb{\beta})V_{k}%
]^{T}\breve{\pmb{b}}_{r}(\pmb{\alpha},\pmb{\beta}),\tag{A.27}
\]
for $k,=1,\ldots,q;r=1,\ldots,d-1$, and
\[%
\begin{array}
[c]{ll}%
\pi_{(q+k)(q+r)} & =[\bar{\pmb{b}}^{T}(\pmb{\alpha},\pmb{\beta})R_{rk}%
(\pmb{\beta},\pmb{\beta})+\breve{\pmb{b}}_{r}^{T}%
(\pmb{\alpha},\pmb{\beta})H_{k}(\pmb{\beta},\pmb{\beta})]\bar{\pmb{b}}%
(\pmb{\alpha},\pmb{\beta})-[\ddot{\Pi}_{kr}^{T}(\pmb{\alpha},\pmb{\beta})\\
& -\bar{\pmb{b}}^{T}(\pmb{\alpha},\pmb{\beta})M_{kr}%
(\pmb{\beta},\pmb{\beta})]\bar{\pmb{b}}(\pmb{\alpha},\pmb{\beta} )\\&+[\dot{\Pi
}_{k}^{T}(\pmb{\alpha},\pmb{\beta})-\bar{\pmb{b}}^{T}%
(\pmb{\alpha},\pmb{\beta})H_{k}(\pmb{\beta} ,\pmb{\beta})]\check{\pmb{b}}%
_{r}(\pmb{\alpha},\pmb{\beta}),
\end{array}
\tag{A.28}
\]
for $k,r=1,\ldots,d-1$, $\bar{\pmb{b}}(\pmb{\alpha},\pmb{\beta})=\Gamma
^{-1}(\pmb{\beta},\pmb{\beta} )\Pi(\pmb{\alpha},\pmb{\beta})$, $\check
{\pmb{b}}_{r}(\pmb{\alpha},\pmb{\beta})=-\Gamma^{-1}(\pmb{\beta}
,\pmb{\beta})$ $(H_{r}^{T}(\pmb{\beta},\pmb{\beta})+H_{r}%
(\pmb{\beta},\pmb{\beta}))\bar{\pmb{b}}(\pmb{\alpha},\pmb{\beta} )+\Gamma
^{-1}(\pmb{\beta},\pmb{\beta})\dot{\Pi}_{r}(\pmb{\alpha},\pmb{\beta})$,
$\dot{\Pi}_{r}(\pmb{\alpha},\pmb{\beta})=\frac{\partial\Pi
(\pmb{\alpha},\pmb{\beta})}{\partial\beta_{r}}$ and $\ddot{\Pi}_{kr}%
(\pmb{\alpha},\pmb{\beta})=\frac{\partial^{2}\Pi(\pmb{\alpha},\pmb{\beta})}%
{\partial\beta_{r}\beta_{k}}$, $M_{kr}(\pmb{\beta},\pmb{\beta}^{\prime})$ is a
$K_{n}\times K_{n}$ matrix whose $(l,l^{\prime})$th element is $E[B_{l}%
(Z^{T}\pmb{\beta})\ddot{B}_{l^{\prime}kr}(Z^{T}\pmb{\beta}^{\prime})]$ and
$\ddot{B}_{lkr}(Z^{T}\pmb{\beta})=\frac{\partial^{2}B_{l}(Z^{T}\pmb{\beta})}%
{\partial\beta_{k}\partial\beta_{r}}$.

\textbf{Proof.} \ \ Let $\tilde{\pi}_{kr}$ be the $(k,r)$th element of
$\ddot{G}_{n}(\pmb{\theta}_{-d},\tilde{\pmb{b}}(\pmb{\theta}_{-d}))$. From
(3.6) and (3.7), we have that
\[
\tilde{\pi}_{kr}=\frac{2}{n}[\tilde{\mathbf{W}}_{k}^{T}\tilde{\mathbf{W}}%
_{r}-\tilde{\mathbf{W}}_{k}^{T}\tilde{\pmb{B}}(\tilde{\pmb{B}}^{T}%
\tilde{\pmb{B}})^{-1}\tilde{\pmb{B}}^{T}\tilde{\mathbf{W}}_{r}%
],\ \ k,r=1,\ldots,q,\tag{A.29}
\]%
\[
\tilde{\pi}_{k(q+r)}=\frac{2}{n}[\tilde{\mathbf{W}}_{k}^{T}\dot{\tilde
{\pmb{B}}}_{r}\tilde{\pmb{b}}+\tilde{\mathbf{W}}_{k}^{T}\tilde{\pmb{B}}%
\dot{\pmb{b}}_{r}],\ \ k,=1,\ldots,q;r=1,\ldots,d-1,\tag{A.30}
\]%
\[
\tilde{\pi}_{(q+k)(q+r)}=\frac{2}{n}(\dot{\tilde{\pmb{B}}}_{r}\tilde
{\pmb{b}}+\tilde{\pmb{B}}\dot{\pmb{b}}_{r})^{T}\dot{\tilde{\pmb{B}}}_{k}%
\tilde{\pmb{b}}-\frac{2}{n}(\tilde{\pmb{Y}}-\tilde{\pmb{W}}\pmb{\alpha}-\tilde
{\pmb{B}}\tilde{\pmb{b}})^{T}(\ddot{\tilde{\pmb{B}}}_{kr}\tilde{\pmb{b}}%
+\dot{\tilde{\pmb{B}}}_{k}\dot{\pmb{b}}_{r}),\tag{A.31}
\]
for $k,r=1,\ldots,d-1$, where $\tilde{\mathbf{W}}_{k}=(\tilde{W}_{1k}%
,\ldots,\tilde{W}_{nk})^{T}$ for $k=1,\ldots,q$, $\tilde{\pmb{B}}%
=\tilde{\pmb{B}}(\pmb{\beta}_{-d})$, $\dot{\tilde{\pmb{B}}}_{r}=\dot
{\tilde{\pmb{B}}}_{r}(\pmb{\beta}_{-d})$ and $\tilde{\pmb{b}}=\tilde
{\pmb{b}}(\pmb{\alpha},\pmb{\beta}_{-d}),$ with for simplicity of notation,
$\dot{\pmb{b}}_{r}=\dot{\pmb{b}}_{r}(\pmb{\alpha},\pmb{\beta}_{-d}%
)=\frac{\partial\tilde{\pmb{b}}(\pmb{\alpha},\pmb{\beta}_{-d})}{\partial
\beta_{r}}$ and $\ddot{\tilde{\pmb{B}}}_{kr}=\ddot{\tilde{\pmb{B}}}%
_{kr}(\pmb{\beta}_{-d})=\frac{\partial^{2}\tilde{\pmb{B}}(\pmb{\beta}_{-d}%
)}{\partial\beta_{k}\partial\beta_{r}}$. Since $(\tilde{\pmb{B}}^{T}%
\tilde{\pmb{B}})^{-1}\tilde{\pmb{B}}^{T}\tilde{\pmb{B}}=I$, we then have
\[
\frac{\partial(\tilde{\pmb{B}}^{T}\tilde{\pmb{B}})^{-1}}{\partial\beta_{r}%
}\tilde{\pmb{B}}^{T}\tilde{\pmb{B}}+(\tilde{\pmb{B}}^{T}\tilde{\pmb{B}}%
)^{-1}\Big(\frac{\partial\tilde{\pmb{B}}^{T}}{\partial\beta_{r}}%
\tilde{\pmb{B}}+\tilde{\pmb{B}}^{T}\frac{\partial\tilde{\pmb{B}}}%
{\partial\beta_{r}}\Big)=0.
\]
Hence,
\[
\frac{\partial(\tilde{\pmb{B}}^{T}\tilde{\pmb{B}})^{-1}}{\partial\beta_{r}%
}=-(\tilde{\pmb{B}}^{T}\tilde{\pmb{B}})^{-1}(\dot{\tilde{\pmb{B}}}_{r}%
^{T}\tilde{\pmb{B}}+\tilde{\pmb{B}}^{T}\dot{\tilde{\pmb{B}}}_{r}%
)(\tilde{\pmb{B}}^{T}\tilde{\pmb{B}})^{-1}.
\]
Note that $\tilde{\pmb{b}}=(\tilde{\pmb{B}}^{T}\tilde{\pmb{B}})^{-1}%
\tilde{\pmb{B}}^{T}(\tilde{\pmb{Y}}-\tilde{\pmb{W}}\pmb{\alpha})$. We further
have
\[
\dot{\pmb{b}}_{r}=-(\tilde{\pmb{B}}^{T}\tilde{\pmb{B}})^{-1}[(\dot
{\tilde{\pmb{B}}}_{r}^{T}\tilde{\pmb{B}}+\tilde{\pmb{B}}^{T}\dot
{\tilde{\pmb{B}}}_{r})(\tilde{\pmb{B}}^{T}\tilde{\pmb{B}})^{-1}\tilde
{\pmb{B}}^{T}(\tilde{\pmb{Y}}-\tilde{\pmb{W}}\pmb{\alpha})-\dot{\tilde
{\pmb{B}}}_{r}^{T}(\tilde{\pmb{Y}}-\tilde{\pmb{W}}\pmb{\alpha})].\tag{A.32}
\]
Similar to the proof of Lemmas A.2 and A.3, we obtain that
\[
\frac{1}{n}(\dot{\tilde{\pmb{B}}}_{r}^{T}\tilde{\pmb{B}}+\tilde{\pmb{B}}%
^{T}\dot{\tilde{\pmb{B}}}_{r})=H_{r}^{T}(\pmb{\beta},\pmb{\beta})+H_{r}%
(\pmb{\beta},\pmb{\beta})+o_{p}(h_{0}^{3}).\tag{A.33}
\]
Furthermore, under Assumption 5, Lemma A.3 yields that $\frac{1}{n}%
\tilde{\pmb{B}}^{T}\tilde{\pmb{B}}=\Gamma(\pmb{\beta},\pmb{\beta})+o_{p}%
(h_{0}^{4})$. Similar to the proof of Lemma 1 of Tang \cite{r32},
we have
$|(\frac{K_{n}}{n}\tilde{\pmb{B}}^{T}\tilde{\pmb{B}})^{-1}-(K_{n}%
\Gamma(\pmb{\beta},\pmb{\beta}))^{-1}|_{\infty}=o_{p}(h_{0}^{3})$. By Lemma A.9
of Huang et al. \cite{r12.5},
we also have that $\Vert(\frac{K_{n}}{n}%
\tilde{\pmb{B}}^{T}\tilde{\pmb{B}})^{-1}\Vert_{\infty}\leq C$ and $\Vert
(K_{n}\Gamma(\pmb{\beta},\pmb{\beta}))^{-1}\Vert_{\infty}\leq C$. Using
(A.33), we have $\Vert H_{r}^{T}(\pmb{\beta},\pmb{\beta})+H_{r}%
(\pmb{\beta},\pmb{\beta})\Vert_{\infty}=O(1)$ and
\[
\Vert\frac{1}{n}(\dot{\tilde{\pmb{B}}}_{r}^{T}\tilde{\pmb{B}}+\tilde
{\pmb{B}}^{T}\dot{\tilde{\pmb{B}}}_{r})\Vert_{\infty}=\Vert H_{r}%
^{T}(\pmb{\beta},\pmb{\beta})+H_{r}(\pmb{\beta},\pmb{\beta})\Vert_{\infty
}+o_{p}(h_{0}^{2})=O_{p}(1).\tag{A.34}
\]
Similar to the proof of (A.20), we obtain $\frac{1}{n}\tilde{\pmb{B}}%
^{T}(\tilde{\pmb{Y}}-\tilde{\pmb{W}}\pmb{\alpha})=\Pi
(\pmb{\alpha},\pmb{\beta})+o_{p}(h_{0}^{3})$. Observe that $\Vert
\Pi(\pmb{\alpha},\pmb{\beta})\Vert_{\infty}=O(1)$ and hence $\Vert\frac{1}%
{n}\tilde{\pmb{B}}^{T}(\tilde{\pmb{Y}}-\tilde{\pmb{W}}\pmb{\alpha})\Vert
_{\infty}=O_{p}(1)$. Let $\vec{\pmb{B}}_{r}=\frac{1}{n}(\dot{\tilde{\pmb{B}}%
}_{r}^{T}\tilde{\pmb{B}}+\tilde{\pmb{B}}^{T}\dot{\tilde{\pmb{B}}}_{r})$,
$\vec{H}_{r}(\pmb{\beta},\pmb{\beta})=H_{r}^{T}(\pmb{\beta},\pmb{\beta})+H_{r}%
(\pmb{\beta},\pmb{\beta})$ and $\vec{\pmb{Y}}=\frac{1}{n}\tilde{\pmb{B}}%
^{T}(\tilde{\pmb{Y}}-\tilde{\pmb{W}}\pmb{\alpha})$. Then
\[%
\begin{array}
[c]{l}%
|(\frac{K_{n}}{n}\tilde{\pmb{B}}^{T}\tilde{\pmb{B}})^{-1}\vec{\pmb{B}}%
_{r}(\frac{K_{n}}{n}\tilde{\pmb{B}}^{T}\tilde{\pmb{B}})^{-1}\vec
{\pmb{Y}}-(K_{n}\Gamma(\pmb{\beta},\pmb{\beta}))^{-1}\vec{\pmb{B}}_{r}%
(\frac{K_{n}}{n}\tilde{\pmb{B}}^{T}\tilde{\pmb{B}})^{-1}\vec{\pmb{Y}}\\
\leq|(\frac{K_{n}}{n}\tilde{\pmb{B}}^{T}\tilde{\pmb{B}})^{-1}-(K_{n}%
\Gamma(\pmb{\beta},\pmb{\beta}))^{-1}|_{\infty}\Vert\vec{\pmb{B}}_{r}%
\Vert_{\infty}\Vert(\frac{K_{n}}{n}\tilde{\pmb{B}}^{T}\tilde{\pmb{B}}%
)^{-1}\Vert_{\infty}\Vert\vec{\pmb{Y}}\Vert_{\infty}\\
=o_{p}(h_{0}^{3})O_{p}(1)O_{p}(1)O_{p}(1)=o_{p}(h_{0}^{3})
\end{array}
\]
and
\[%
\begin{array}
[c]{l}%
|(K_{n}\Gamma(\pmb{\beta},\pmb{\beta}))^{-1}\vec{\pmb{B}}_{r}(\frac{K_{n}}%
{n}\tilde{\pmb{B}}^{T}\tilde{\pmb{B}})^{-1}\vec{\pmb{Y}}-(K_{n}\Gamma
(\pmb{\beta},\pmb{\beta}))^{-1}\vec{H}_{r}(\pmb{\beta},\pmb{\beta})(\frac
{K_{n}}{n}\tilde{\pmb{B}}^{T}\tilde{\pmb{B}})^{-1}\vec{\pmb{Y}}|_{\infty}\\
\leq\Vert(K_{n}\Gamma(\pmb{\beta},\pmb{\beta}))^{-1}\Vert_{\infty}%
|\vec{\pmb{B}}_{r}-\vec{H}_{r}(\pmb{\beta},\pmb{\beta})|_{\infty}\Vert
(\frac{K_{n}}{n}\tilde{\pmb{B}}^{T}\tilde{\pmb{B}})^{-1}\Vert_{\infty}%
\Vert\vec{\pmb{Y}}\Vert_{\infty}\\
=O(1)o_{p}(h_{0}^{3})O_{p}(1)O_{p}(1)=o_{p}(h_{0}^{3}).
\end{array}
\]
Furthermore, it holds that
\[%
\begin{array}
[c]{l}%
|(\frac{K_{n}}{n}\tilde{\pmb{B}}^{T}\tilde{\pmb{B}})^{-1}\vec{\pmb{B}}%
_{r}(\frac{K_{n}}{n}\tilde{\pmb{B}}^{T}\tilde{\pmb{B}})^{-1}\vec{\pmb{Y}}\\
-(K_{n}\Gamma(\pmb{\beta},\pmb{\beta}))^{-1}\vec{H}_{r}%
(\pmb{\beta},\pmb{\beta})(K_{n}\Gamma(\pmb{\beta},\pmb{\beta}))^{-1}%
\Pi(\pmb{\alpha},\pmb{\beta})|_{\infty}=o_{p}(h_{0}^{3})
\end{array}
\tag{A.35}
\]
Under Assumption 5, similar to the proof of (A.20), we deduce that
\[
\frac{1}{n}\dot{\tilde{\pmb{B}}}_{r}^{T}(\tilde{\pmb{Y}}-\tilde{\pmb{W}}%
\pmb{\alpha})=\dot{\Pi}_{r}(\pmb{\alpha},\pmb{\beta})+o_{p}(h_{0}%
^{2}).\tag{A.36}
\]
Similar to the proof of (A.35), we further deduce that
\[
|(\frac{K_{n}}{n}\tilde{\pmb{B}}^{T}\tilde{\pmb{B}})^{-1}(\frac{1}{n}%
\dot{\tilde{\pmb{B}}}_{r}^{T}(\tilde{\pmb{Y}}-\tilde{\pmb{W}}%
\pmb{\alpha}))-(K_{n}\Gamma(\pmb{\beta},\pmb{\beta}))^{-1}\dot{\Pi}%
_{r}(\pmb{\alpha},\pmb{\beta})|_{\infty}=o_{p}(h_{0}^{2}).\tag{A.37}
\]
Combining (A.32), (A.35) and (A.37), we then have
\[
|\dot{\pmb{b}}_{r}-\check{\pmb{b}}_{r}(\pmb{\alpha},\pmb{\beta})|_{\infty
}=o_{p}(h_{0}).\tag{A.38}
\]
By arguments similar to those used in the proof of (A.35), we further have
that
\[
\frac{1}{n}\dot{\pmb{b}}_{r}^{T}\pmb{B}^{T}\dot{\tilde{\pmb{B}}}_{k}%
\tilde{\pmb{b}}=\check{\pmb{b}}_{r}^{T}(\pmb{\alpha},\pmb{\beta})H_{k}%
(\pmb{\beta},\pmb{\beta})\bar{\pmb{b}}(\pmb{\alpha},\pmb{\beta})+o_{p}%
(1).\tag{A.39}
\]
Similar to the proof of (A.39), we obtain that
\[
\frac{1}{n}\tilde{\pmb{b}}^{T}\dot{\tilde{\pmb{B}}}_{r}^{T}\dot{\tilde
{\pmb{B}}}_{k}\tilde{\pmb{b}}=\bar{\pmb{b}}^{T}%
(\pmb{\alpha},\pmb{\beta})R_{rk}(\pmb{\beta},\pmb{\beta})\bar{\pmb{b}}%
(\pmb{\alpha},\pmb{\beta})+o_{p}(1),\tag{A.40}
\]%
and 
\[%
\begin{array}
[c]{ll}%
&\frac{1}{n}(\tilde{\pmb{Y}}-\tilde{\pmb{W}}\pmb{\alpha}-\tilde{\pmb{B}}%
\tilde{\pmb{b}})^{T}(\ddot{\tilde{\pmb{B}}}_{kr}\tilde{\pmb{b}}+\dot
{\tilde{\pmb{B}}}_{k}\dot{\pmb{b}}_{r}) \\& =[\ddot{\Pi}_{kr}^{T}%
(\pmb{\alpha},\pmb{\beta})-\bar{\pmb{b}}^{T}(\pmb{\alpha},\pmb{\beta})M_{kr}%
(\pmb{\beta},\pmb{\beta})]\bar{\pmb{b}}(\pmb{\alpha},\pmb{\beta})\\
& +[\dot{\Pi}_{k}^{T}(\pmb{\alpha},\pmb{\beta})-\bar{\pmb{b}}^{T}%
(\pmb{\alpha},\pmb{\beta})H_{k}(\pmb{\beta},\pmb{\beta})]\check{\pmb{b}}%
_{r}(\pmb{\alpha},\pmb{\beta})+o_{p}(1).
\end{array}
\]
Now (A.28) follows from (A.31), (A.39), (A.40) and the preceding expression.
Using the fact that $\frac{1}{n}\sum_{i=1}^{n}(W_{ik}-\frac{1}{n}\sum
_{l=1}^{n}W_{lk}\check{\xi}_{il})(W_{ir}-\frac{1}{n}\sum_{l=1}^{n}W_{lr}%
\check{\xi}_{il})=E(V_{k}V_{r})+o_{p}(1)$, (A.26) and (A.27) can be proved in
a similar fashion. This completes the proof of Lemma A.4.

\textbf{Lemma\ A.5.} \ Under Assumptions 1 to 3 and 5, it holds that
\[
\sum_{j=1}^{m}\lambda_{j}[a_{j}-\frac{1}{\hat{\lambda}_{j}}(\frac{1}{n}%
\sum_{l=1}^{n}\zeta_{l}\hat{\xi}_{lj})]^{2}=O_{p}(n^{-1}\lambda_{m}^{-1}m),
\]
where $\zeta_{l}=\sum_{q=1}^{\infty}a_{q}\xi_{lq}$.

\textbf{Proof} \ \ Set $S_{1}=\sum_{j=1}^{m}\lambda_{j}[a_{j}-\frac{1}%
{\lambda_{j}}(\frac{1}{n}\sum_{l=1}^{n}\zeta_{l}\xi_{lj})]^{2}$, $S_{2}%
=\sum_{j=1}^{m}\frac{1}{\lambda_{j}}[\frac{1}{n}\sum_{l=1}^{n}\zeta_{l}%
$ $(\hat{\xi}_{lj}-\xi_{lj})]^{2}$ and $S_{3}=\sum_{j=1}^{m}\lambda_{j}(\frac
{1}{\hat{\lambda}_{j}}-\frac{1}{\lambda_{j}})^{2}(\frac{1}{n}\sum_{l=1}%
^{n}\zeta_{l}\hat{\xi}_{lj})^{2}$. Note that
\[
\sum_{j=1}^{m}\lambda_{j}[a_{j}-\frac{1}{\hat{\lambda}_{j}}(\frac{1}{n}%
\sum_{l=1}^{n}\zeta_{l}\hat{\xi}_{lj})]^{2}\leq3(S_{1}+S_{2}+S_{3}%
).\tag{A.41}
\]
Since $E[a_{j}-\frac{1}{\lambda_{j}}(\frac{1}{n}\sum_{l=1}^{n}\zeta_{l}%
\xi_{lj})]=0$, then from Assumptions 1to 3, we obtain
\[
E(S_{1})=\sum_{j=1}^{m}\frac{1}{\lambda_{j}}Var(\frac{1}{n}\sum_{l=1}^{n}%
\zeta_{l}\xi_{lj})\leq\sum_{j=1}^{m}\frac{1}{n^{2}\lambda_{j}}\sum_{l=1}%
^{n}E(\zeta_{l}^{2}\xi_{lj}^{2})\leq Cm/n.\tag{A.42}
\]
Similar to the proof of (A.6), (A.7) and using Assumption 5, we deduce that
\[
S_{2}=O_{p}(n^{-2}\lambda_{m}^{-2}m^{3}\log m+n^{-1}\lambda_{m}^{-1}%
m)=O_{p}(n^{-1}\lambda_{m}^{-1}m)\tag{A.43}
\]
and
\[%
\begin{array}
[c]{ll}%
S_{3} & \leq C\sum_{j=1}^{m}\frac{(\hat{\lambda}_{j}-\lambda_{j})^{2}}%
{\lambda_{j}^{3}}\Big(\bar{\zeta}_{j}^{2}+[\frac{1}{n}\sum_{l=1}^{n}\zeta
_{l}(\hat{\xi}_{lj}-\xi_{lj})]^{2}\Big)[1+o_{p}(1)]\\
& =O_{p}(n^{-1}\lambda_{m}^{-1}+n^{-3}\lambda_{m}^{-4}m^{3}\log m+n^{-2}%
\lambda_{m}^{-3}m)=O_{p}(n^{-1}\lambda_{m}^{-1}).
\end{array}
\tag{A.44}
\]
Now Lemma A.5 follows from combining (A.41) to (A.44).

\textbf{Lemma\ A.6.} \ Denote
\[
\begin{array}
[c]{ll}%
\dot{g}_{0r}(Z_{i})=\frac{\partial g_{0}(Z_{i}^{T}\pmb{\beta})}{\partial
\beta_{r}}|_{\pmb{\beta}=\pmb{\beta}_{0}}\\=\sum_{k=1}^{K_{n}}b_{0k}%
B_{k}^{\prime}(Z_{i}^{T}\pmb{\beta} _{0})\Big(Z_{ir}-\frac{\beta_{0r}Z_{id}%
}{\sqrt{1-(\beta_{01}^{2}+\ldots+\beta_{0(d-1)}^{2})}}\Big)
\end{array}
\]
for $r=1,\ldots,d-1$ and $A_{ri}=\dot{g}_{0r}(Z_{i})-\frac{1}{n}\sum_{l=1}%
^{n}\dot{g}_{0r}(Z_{l})\tilde{\xi}_{il}$. Under Assumptions 1, 2, 4 and 5, it
holds that
\[
\sum_{j=1}^{m}\lambda_{j}^{-1}(\sum_{i=1}^{n}\xi_{ij}A_{ri})^{2}%
=O_{p}(nm+\lambda_{m}^{-2}m^{4}\log m).
\]

\textbf{Proof} \ \ Let $A_{ri}^{\ast}=\dot{g}_{0r}(Z_{i})-\sum_{j^{\prime}%
=1}^{m}\frac{1}{\lambda_{j^{\prime}}}(\frac{1}{n}\sum_{l=1}^{n}\dot{g}%
_{0r}(Z_{l})\xi_{lj^{\prime}})\xi_{ij^{\prime}}$. Observe that
\[%
\begin{array}
[c]{ll}%
(\sum_{i=1}^{n}\xi_{ij}A_{ri})^{2} & \leq4\Big(\sum_{i=1}^{n}\xi_{ij}%
A_{ri}^{\ast}\Big)^{2}\\
& +4\Big(\sum_{i=1}^{n}\xi_{ij}\sum_{j^{\prime}=1}^{m}\frac{1}{\lambda
_{j^{\prime}}}[\frac{1}{n}\sum_{l=1}^{n}\dot{g}_{0r}(Z_{l})(\hat{\xi
}_{lj^{\prime}}-\xi_{lj^{\prime}})]\xi_{ij^{\prime}}\Big)^{2}\\
& +4\Big(\sum_{i=1}^{n}\xi_{ij}\sum_{j^{\prime}=1}^{m}(\frac{1}{\hat{\lambda
}_{j^{\prime}}}-\frac{1}{\lambda_{j^{\prime}}})[\frac{1}{n}\sum_{l=1}^{n}%
\dot{g}_{0r}(Z_{l})\hat{\xi}_{lj^{\prime}}]\xi_{ij^{\prime}}\Big)^{2}\\
& +4\Big(\sum_{i=1}^{n}\xi_{ij}\sum_{j^{\prime}=1}^{m}\frac{1}{\hat{\lambda
}_{j^{\prime}}}[\frac{1}{n}\sum_{l=1}^{n}\dot{g}_{0r}(Z_{l})\hat{\xi
}_{lj^{\prime}}](\hat{\xi}_{ij^{\prime}}-\xi_{ij^{\prime}})\Big)^{2}\\
& =:4(T_{j1}+T_{j2}+T_{j3}+T_{j4}).
\end{array}
\tag{A.45}
\]
By direct computations and using Assumption 1, we obtain
\[%
\begin{array}
[c]{ll}%
E(\xi_{ij}^{2}{A_{ri}^{\ast}}^{2}) & \leq2E(\xi_{ij}^{2}\dot{g}_{0r}^{2}%
(Z_{i}))+2E[\xi_{ij}^{2}(\sum_{j^{\prime}=1}^{m}\frac{1}{\lambda_{j^{\prime}}%
}(\frac{1}{n}\sum_{l=1}^{n}\dot{g}_{0r}(Z_{l})\xi_{lj^{\prime}})\xi
_{ij^{\prime}})^{2}]\\
& \leq C(\lambda_{j}+m\lambda_{j}/n^{2}+(n-1)m\lambda_{j}/n^{2}+m^{2}%
\lambda_{j}/n^{2})\leq C\lambda_{j}%
\end{array}
\]
and
\[
|\sum_{i_{1}\neq i_{2}}E(\xi_{i_{1}j}\xi_{i_{2}j}{A_{ri_{1}}^{\ast}}%
{A_{ri_{2}}^{\ast}})|\leq C[(n-1)(n+2)\lambda_{j}/n+(n-1)m\lambda_{j}/n]\leq
Cn\lambda_{j}.
\]
Hence, it follows that
\[
E(T_{j1})=\sum_{i=1}^{n}E(\xi_{ij}^{2}{A_{ri}^{\ast}}^{2})+\sum_{i_{1}\neq
i_{2}}E(\xi_{i_{1}j}\xi_{i_{2}j}{A_{ri_{1}}^{\ast}}{A_{ri_{2}}^{\ast}})\leq
Cn\lambda_{j}.\tag{A.46}
\]
Similar to the proof of (A.6) and using Assumption 1, we have
\[
\sum_{j^{\prime}=1}^{m}\frac{1}{\lambda_{j^{\prime}}}[\frac{1}{n}\sum
_{l=1}^{n}\dot{g}_{0r}(Z_{l})(\hat{\xi}_{lj^{\prime}}-\xi_{lj^{\prime}}%
)]^{2}=O_{p}(n^{-2}\lambda_{m}^{-2}m^{3}\log m).
\]
Since $\sum_{j^{\prime}=1}^{m}\frac{1}{\lambda_{j^{\prime}}}E(\sum_{i=1}%
^{n}\xi_{ij}\xi_{ij^{\prime}})^{2}\leq Cn^{2}\lambda_{j}$, then
\[%
\begin{array}
[c]{ll}%
\sum_{j=1}^{m}\lambda_{j}^{-1}T_{j2} & \leq\Big(\sum_{j^{\prime}=1}^{m}%
\frac{1}{\lambda_{j^{\prime}}}[\frac{1}{n}\sum_{l=1}^{n}\dot{g}_{0r}%
(Z_{l})(\hat{\xi}_{lj^{\prime}}-\xi_{lj^{\prime}})]^{2}\Big)\\
& \times\Big(\sum_{j=1}^{m}\lambda_{j}^{-1}\sum_{j^{\prime}=1}^{m}\frac
{1}{\lambda_{j^{\prime}}}(\sum_{i=1}^{n}\xi_{ij}\xi_{ij^{\prime}})^{2}\Big)\\
& =O_{p}(n^{-2}\lambda_{m}^{-2}m^{3}\log m)O_{p}(n^{2}m)=O_{p}(\lambda
_{m}^{-2}m^{4}\log m).
\end{array}
\tag{A.47}
\]
Similar to the proof (A.7) and using Assumption 5, we deduce that
\[%
\begin{array}
[c]{ll}%
\sum_{j=1}^{m}\lambda_{j}^{-1}T_{j3} & \leq\Big(\sum_{j^{\prime}=1}^{m}%
\lambda_{j^{\prime}}(\frac{1}{\hat{\lambda}_{j^{\prime}}}-\frac{1}%
{\lambda_{j^{\prime}}})^{2}[\frac{1}{n}\sum_{l=1}^{n}\dot{g}_{0r}(Z_{l}%
)\hat{\xi}_{lj^{\prime}}]^{2}\Big)\\
& \times\Big(\sum_{j=1}^{m}\lambda_{j}^{-1}\sum_{j^{\prime}=1}^{m}\frac
{1}{\lambda_{j^{\prime}}}(\sum_{i=1}^{n}\xi_{ij}\xi_{ij^{\prime}})^{2}\Big)\\
& =O_{p}(\lambda_{m}^{-2}m^{2}+n^{-1}\lambda_{m}^{-4}m^{4}\log m)=O_{p}%
(\lambda_{m}^{-2}m^{2}\log m).
\end{array}
\tag{A.48}
\]
and
\[%
\begin{array}
[c]{ll}%
\sum_{j=1}^{m}\lambda_{j}^{-1}T_{j4} & \leq\Big(\sum_{j^{\prime}=1}^{m}%
\frac{1}{\lambda_{j^{\prime}}^{2}}[\frac{1}{n}\sum_{l=1}^{n}\dot{g}_{0r}%
(Z_{l})\hat{\xi}_{lj^{\prime}}]^{2}\Big)[1+o_{p}(1)]\\
& \times\Big(\sum_{j=1}^{m}\frac{1}{\lambda_{j}}\sum_{i=1}^{n}\xi_{ij}%
^{2}\Big)\Big(\sum_{j^{\prime}=1}^{m}\sum_{i=1}^{n}(\hat{\xi}_{ij^{\prime}%
}-\xi_{ij^{\prime}})^{2}\Big)\\
& =O_{p}(n^{-1}\lambda_{m}^{-1}m^{5}\log m+n^{-2}\lambda_{m}^{-3}m^{7}(\log
m)^{2})\\&=o_{p}(\lambda_{m}^{-2}m^{4}\log m).
\end{array}
\tag{A.49}
\]
Now Lemma A.6 follows from combining (A.45)-(A.49) and Assumption 5.

\textbf{Lemma\ A.7.} \ Under the Assumptions 1-3 and 5, it holds that
\[
n^{-1/2}|\sum_{j=1}^{m}\frac{1}{\hat{\lambda}_{j}}(\frac{1}{n}\sum_{l=1}%
^{n}\zeta_{l}\hat{\xi}_{lj})\sum_{i=1}^{n}(\hat{\xi}_{ij}-\xi_{ij}%
)A_{ri}|=o_{p}(1).
\]

\textbf{Proof} \ \ Let $\check{\zeta}_{j}=\frac{1}{n}\sum_{l=1}^{n}\zeta
_{l}\hat{\xi}_{lj}$. Applying the Cauchy-Schwarz inequality, we obtain
\[
\Big(\sum_{j=1}^{m}\frac{1}{\hat{\lambda}_{j}}\check{\zeta}_{j}\sum_{i=1}%
^{n}(\hat{\xi}_{ij}-\xi_{ij})A_{ri}\Big)^{2}\leq\Big(\sum_{j=1}^{m}\frac
{1}{\hat{\lambda}_{j}^{2}}\check{\zeta}_{j}^{2}\Big)\Big(\sum_{j=1}^{m}%
(\sum_{i=1}^{n}(\hat{\xi}_{ij}-\xi_{ij})A_{ri})^{2}\Big).
\]
Using (A.4), (A.5), Assumption 5, Parseval's identity and some arguments
similar to those used to prove Lemma A.6, we deduce that
\[%
\begin{array}
[c]{ll}%
\sum_{j=1}^{m}(\sum_{i=1}^{n}(\hat{\xi}_{ij}-\xi_{ij})A_{ri})^{2} & \\
\leq2\sum_{j=1}^{m}[(\sum_{k\neq j}(\hat{\lambda}_{j}-\lambda_{k})^{-1}%
\int\Delta\hat{\phi}_{j}\phi_{k}\sum_{i=1}^{n}\xi_{ik}A_{ri})^{2}&\\+(\sum
_{i=1}^{n}\xi_{ij}A_{ri})^{2}(\int(\hat{\phi}_{j}-\phi_{j})\phi_{j})^{2}] & \\
\leq C|\Vert\Delta\Vert|^{2}\sum_{j=1}^{m}[\sum_{k\neq j}(\hat{\lambda}%
_{j}-\lambda_{k})^{-2}(\sum_{i=1}^{n}\xi_{ik}A_{ri})^{2}+\lambda_{j}^{-2}%
j^{2}(\sum_{i=1}^{n}\xi_{ij}A_{ri})^{2}] & \\
=O_{p}(\lambda_{m}^{-1}m^{3}\log m+n^{-1}\lambda_{m}^{-3}m^{6}\log
m)=o_{p}(n). &
\end{array}
\]
Similar to the proof of (A.7) and using Assumption 5 , we obtain that
\[
\sum_{j=1}^{m}\frac{1}{\hat{\lambda}_{j}^{2}}\check{\zeta}_{j}^{2}%
=O_{p}(n^{-1}\lambda_{m}^{-1}m+1+n^{-2}\lambda_{m}^{-3}m^{3}\log
m+n^{-1}\lambda_{m}^{-2}m)=o_{p}(1).
\]
This completes the proof of Lemma A.7.

\textbf{Lemma\ A.8.} \ Set $\tilde{\zeta}_{i}=\zeta_{i}-\frac{1}{n}\sum
_{l=1}^{n}\zeta_{l}\tilde{\xi}_{il}$. Under Assumptions 1-4 and 5, it holds
that
\[
n^{-1/2}|\sum_{i=1}^{n}\tilde{\zeta}_{i}A_{ri}|=o_{p}(1).
\]
\textbf{Proof} \ Observe that
\[%
\begin{array}
[c]{ll}%
\sum_{i=1}^{n}\tilde{\zeta}_{i}A_{ri} & =\sum_{j=1}^{m}[a_{j}-\frac{1}%
{\hat{\lambda}_{j}}(\frac{1}{n}\sum_{l=1}^{n}\zeta_{l}\hat{\xi}_{lj}%
)]\sum_{i=1}^{n}\xi_{ij}A_{ri}\\
& -\sum_{j=1}^{m}\frac{1}{\hat{\lambda}_{j}}(\frac{1}{n}\sum_{l=1}^{n}%
\zeta_{l}\hat{\xi}_{lj})\sum_{i=1}^{n}(\hat{\xi}_{ij}-\xi_{ij})A_{ri}%
\\&+\sum_{j=m+1}^{\infty}a_{j}\sum_{i=1}^{n}\xi_{ij}A_{ri}.
\end{array}
\tag{A.50}
\]
Lemmas A.5, A.6 and Assumption 5 imply that
\[%
\begin{array}
[c]{l}%
n^{-\frac{1}{2}}|\sum_{j=1}^{m}[a_{j}-\frac{1}{\hat{\lambda}_{j}}(\frac{1}%
{n}\sum_{l=1}^{n}\zeta_{l}\hat{\xi}_{lj})]\sum_{i=1}^{n}\xi_{ij}A_{ri}|\\
\leq n^{-\frac{1}{2}}\Big(\sum_{j=1}^{m}\lambda_{j}[a_{j}-\frac{1}%
{\hat{\lambda}_{j}}(\frac{1}{n}\sum_{l=1}^{n}\zeta_{l}\hat{\xi}_{lj}%
)]^{2}\Big)^{\frac{1}{2}}\Big(\sum_{j=1}^{m}\lambda_{j}^{-1}(\sum_{i=1}^{n}%
\xi_{ij}A_{ri})^{2}\Big)^{\frac{1}{2}}\\
=O_{p}(n^{-1/2}\lambda_{m}^{-1/2}m+n^{-1}\lambda_{m}^{-3/2}m^{5/2}(\log
m)^{1/2})=o_{p}(1).
\end{array}
\tag{A.51}
\]
By arguments similar to those used in the proof of Lemma A.6 and using Lemma
6.1 of Cardot et al. \cite{r5},
we obtain that
\[%
\begin{array}
[c]{ll}%
&(\sum_{j=m+1}^{\infty}a_{j}\sum_{i=1}^{n}\xi_{ij}A_{ri})^{2} \\& \leq
(\sum_{j=m+1}^{\infty}a_{j}^{2})(\sum_{j=m+1}^{\infty}(\sum_{i=1}^{n}\xi
_{ij}A_{ri})^{2})\\
& =O_{p}(nm^{-2\gamma+1}+\lambda_{m}^{-2}m^{-2\gamma+4}\log m)\sum
_{j=m+1}^{\infty}\lambda_{j}\\&=o_{p}(n).
\end{array}
\tag{A.52}
\]
Now Lemma A.8 follows from combining (A.50)-(A.52) and Lemma A.7.

\textbf{Lemma\ A.9.} \ Suppose that Assumptions 1-5 hold. Then
\[
n^{-1/2}(\tilde{\pmb{Y}}-\tilde{\pmb{W}}\pmb{\alpha}_{0}-\tilde{\pmb{B}}%
(\pmb{\beta}_{0})\pmb{b}_{0})^{T}\dot{\tilde{\pmb{B}}}_{r}(\pmb{\beta}_{0}%
)\pmb{b}_{0}=n^{-1/2}\sum_{i=1}^{n}\dot{g}_{0r}(Z_{i})\varepsilon_{i}%
+o_{p}(1).
\]
\textbf{Proof} \ Using arguments similar to those used to prove Lemmas A.6 and
A.7, we deduce that $\sum_{i=1}^{n}A_{ri}^{2}=O_{p}(n)$, $n^{-1/2}\sum
_{i=1}^{n}\varepsilon_{i}(\frac{1}{n}\sum_{l=1}^{n}\dot{g}_{0r}(Z_{l}%
)\tilde{\xi}_{il})=o_{p}(1)$ and
\[
n^{-1/2}\sum_{i=1}^{n}(\frac{1}{n}\sum_{l=1}^{n}\varepsilon_{l}\tilde{\xi
}_{il})A_{ri}=o_{p}(1),\ \ \ n^{-1/2}\sum_{i=1}^{n}(\frac{1}{n}\sum_{l=1}%
^{n}R(Z_{l}^{T}\pmb{\beta}_{0})\tilde{\xi}_{il})A_{ri}=o_{p}(1).
\]
Hence
\[
n^{-1/2}\tilde{\pmb{\varepsilon}}^{T}\dot{\tilde{\pmb{B}}}_{r}(\pmb{\beta}_{0}%
)\pmb{b}_{0}=n^{-1/2}\sum_{i=1}^{n}\dot{g}_{0r}(Z_{i})\varepsilon_{i}%
+o_{p}(1).\tag{A.53}
\]
Using (3.1) and the assumption
$nh^{2p}\rightarrow0$, it follows that
$(\sum_{i=1}^{n}R(Z_{i}^{T}\pmb{\beta}_{0})A_{ri})^{2}$ $\leq(\sum_{i=1}^{n}%
R^{2}(Z_{i}^{T}\pmb{\beta}_{0}))$ $(\sum_{i=1}^{n}A_{ri}^{2})=o_{p}(n)$.
Consequently, we have
\[
n^{-1/2}\tilde{\pmb{R}}^{T}\dot{\tilde{\pmb{B}}}_{r}(\pmb{\beta}_{0}%
)\pmb{b}_{0}=o_{p}(1),\tag{A.54}
\]
where $\tilde{\pmb{R}}=(\tilde{R}(Z_{1}^{T}\pmb{\beta}_{0}),\ldots,\tilde
{R}(Z_{n}^{T}\pmb{\beta}_{0}))^{T}$ and $\tilde{R}(Z_{i}^{T}\pmb{\beta}_{0}%
)=R(Z_{i}^{T}\pmb{\beta} _{0})-$\\ $\frac{1}{n}\sum_{l=1}^{n}R(Z_{l}^{T}%
\pmb{\beta}_{0})\tilde{\xi}_{il}$. Now Lemma A.9 follows from Lemma A.8,
(A.53) and (A.54).

\textbf{Lemma\ A.10.} \ Under the assumptions of Theorem 2, it holds that
\[
\begin{array}
[c]{l}%
n^{-\frac{1}{2}}(\tilde{\pmb{b}}(\pmb{\alpha}_{0},\pmb{\beta}_{0}%
)-\pmb{b}_{0})^{T}\tilde{\pmb{B}}^{T}(\pmb{\beta}_{0})\dot{\tilde{\pmb{B}}%
}_{r}(\pmb{\beta}_{0})\pmb{b}_{0} \\=n^{-\frac{1}{2}}\pmb{\varepsilon}^{T}%
\pmb{B}(\pmb{\beta}_{0})\Gamma^{-1}(\pmb{\beta}_{0},\pmb{\beta}_{0}%
)H_{r}(\pmb{\beta}_{0},\pmb{\beta}_{0})\pmb{b}_{0} +o_{p}(1),
\end{array}
\]
where $\pmb{B}(\pmb{\beta}_{0})=(\pmb{B}(Z_{1}^{T}\pmb{\beta}_{0}%
),\ldots,\pmb{B}(Z_{n}^{T}\pmb{\beta} _{0}))^{T}$.

\textbf{Proof} \ Note that $\tilde{\pmb{b}}(\pmb{\alpha}_{0},\pmb{\beta}_{0}%
)-\pmb{b}_{0}=(\tilde{\pmb{B}}^{T}(\pmb{\beta}_{0})\tilde{\pmb{B}}%
(\pmb{\beta}_{0}))^{-1}\tilde{\pmb{B}}^{T}(\pmb{\beta}_{0})(\tilde
{\pmb{Y}}-\tilde{\pmb{W}}\pmb{\alpha}_{0}-\tilde{\pmb{B}}(\pmb{\beta}_{0}%
)\pmb{b}_{0})$. By (A.33), we obtain
\[%
\begin{array}
[c]{ll}%
|\frac{1}{n}\tilde{\pmb{B}}^{T}(\pmb{\beta}_{0})\dot{\tilde{\pmb{B}}}%
_{r}(\pmb{\beta}_{0})\pmb{b}_{0}|_{\infty} & =|H_{r}(\pmb{\beta}_{0}%
,\pmb{\beta}_{0})\pmb{b}_{0}|_{\infty}+o_{p}(h_{0}^{2})\\
& \leq\max_{1\leq k\leq K_{n}}E[B_{k}(Z^{T}\pmb{\beta}_{0})|\dot{g}%
_{0r}(Z)|]+o_{p}(h_{0}^{2})\\&=O_{p}(h_{0})
\end{array}
\]
Similar to Lemma A.9, we have $\Vert n^{-\frac{1}{2}}%
(\tilde{\pmb{Y}}-\tilde{\pmb{W}}\pmb{\alpha}_{0}-\pmb{\varepsilon}-\tilde
{\pmb{B}}(\pmb{\beta}_{0})\pmb{b}_{0})^{T}\tilde{\pmb{B}}(\pmb{\beta}_{0}%
)\Vert_{\infty}=o_{p}(1)$ and $\Vert n^{-\frac{1}{2}}\pmb{\varepsilon}^{T}%
\tilde{\pmb{B}}(\pmb{\beta}_{0})\Vert_{\infty}=O_{p}(K_{n}^{1/2})$. Hence
\[%
\begin{array}
[c]{l}%
n^{-\frac{1}{2}}|(\tilde{\pmb{Y}}-\tilde{\pmb{W}}\pmb{\alpha}_{0}%
-\pmb{\varepsilon}-\tilde{\pmb{B}}(\pmb{\beta}_{0})\pmb{b}_{0})^{T}%
\tilde{\pmb{B}}(\pmb{\beta}_{0})(\tilde{\pmb{B}}^{T}(\pmb{\beta}_{0}%
)\tilde{\pmb{B}}(\pmb{\beta}_{0}))^{-1}\tilde{\pmb{B}}^{T}(\pmb{\beta}_{0}%
)\dot{\tilde{\pmb{B}}}_{r}(\pmb{\beta} _{0})\pmb{b}_{0}|\\
\leq K_{n}\Vert n^{-\frac{1}{2}}|(\tilde{\pmb{Y}}-\tilde{\pmb{W}}%
\pmb{\alpha}_{0}-\pmb{\varepsilon} -\tilde{\pmb{B}}(\pmb{\beta}_{0}%
)\pmb{b}_{0})^{T}\tilde{\pmb{B}}(\pmb{\beta}_{0})\Vert_{\infty}\Vert
(\frac{K_{n}}{n}\tilde{\pmb{B}}^{T}(\pmb{\beta}_{0})\tilde{\pmb{B}}%
(\pmb{\beta}_{0}))^{-1}\Vert_{\infty}\\
\times|\frac{1}{n}\tilde{\pmb{B}}^{T}(\pmb{\beta}_{0})\dot{\tilde{\pmb{B}}%
}_{r}(\pmb{\beta}_{0})\pmb{b}_{0}|_{\infty}=K_{n}o_{p}(1)O_{p}(1)O_{p}%
(h_{0})=o_{p}(1).
\end{array}
\]
Using arguments similar to those used in the proof of (A.35), we can deduce
that
\[
|(\tilde{\pmb{B}}^{T}(\pmb{\beta}_{0})\tilde{\pmb{B}}(\pmb{\beta}_{0}%
))^{-1}\tilde{\pmb{B}}^{T}(\pmb{\beta} _{0})\dot{\tilde{\pmb{B}}}%
_{r}(\pmb{\beta}_{0})\pmb{b}_{0}-\Gamma^{-1}(\pmb{\beta}_{0},\pmb{\beta}_{0}%
)H_{r}(\pmb{\beta}_{0},\pmb{\beta}_{0})\pmb{b}_{0}|_{\infty}=o_{p}(h_{0}).
\]
Hence
\[%
\begin{array}
[c]{l}%
|n^{-\frac{1}{2}}\pmb{\varepsilon}^{T}\tilde{\pmb{B}}(\pmb{\beta}_{0}%
)[(\tilde{\pmb{B}}^{T}(\pmb{\beta} _{0})\tilde{\pmb{B}}(\pmb{\beta}_{0}%
))^{-1}\tilde{\pmb{B}}^{T}(\pmb{\beta}_{0})\dot{\tilde{\pmb{B}}}%
_{r}(\pmb{\beta}_{0})\pmb{b}_{0}-\Gamma^{-1}(\pmb{\beta}_{0},\pmb{\beta}_{0}%
)H_{r}(\pmb{\beta}_{0},\pmb{\beta} _{0})\pmb{b}_{0}]|\\
\leq\Vert n^{-\frac{1}{2}}\pmb{\varepsilon}^{T}\tilde{\pmb{B}}(\pmb{\beta}_{0}%
)\Vert_{\infty}|(\tilde{\pmb{B}}^{T}(\pmb{\beta}_{0})\tilde{\pmb{B}}%
(\pmb{\beta}_{0}))^{-1}\tilde{\pmb{B}}^{T}(\pmb{\beta} _{0})\dot
{\tilde{\pmb{B}}}_{r}(\pmb{\beta}_{0})\pmb{b}_{0}-\\\Gamma^{-1}(\pmb{\beta}_{0}%
,\pmb{\beta}_{0})H_{r}(\pmb{\beta}_{0},\pmb{\beta}_{0})\pmb{b}_{0}|_{\infty}\\
=O_{p}(K_{n}^{1/2})o_{p}(h_{0})=o_{p}(1).
\end{array}
\]
Using arguments similar to those used to prove Lemmas A.6 and A.7, we deduce
that
\[
\Vert n^{-\frac{1}{2}}\pmb{\varepsilon}^{T}(\tilde{\pmb{B}}(\pmb{\beta}_{0}%
)-\pmb{B}(\pmb{\beta} _{0}))\Vert_{\infty}=n^{-\frac{1}{2}}\sum_{k=1}^{K_{n}%
}|\sum_{i=1}^{n}\varepsilon_{i}[\frac{1}{n}\sum_{l=1}^{n}B_{k}(Z_{l}%
^{T}\pmb{\beta}_{0})\tilde{\xi}_{il}]|=o_{p}(1).
\]
Therefore,
\[%
\begin{array}
[c]{l}%
|n^{-\frac{1}{2}}\pmb{\varepsilon}^{T}(\tilde{\pmb{B}}(\pmb{\beta}_{0}%
)-\pmb{B}(\pmb{\beta} _{0}))\Gamma^{-1}(\pmb{\beta}_{0},\pmb{\beta}_{0}%
)H_{r}(\pmb{\beta}_{0},\pmb{\beta}_{0})\pmb{b}_{0}]|\\
\leq K_{n}\Vert n^{-\frac{1}{2}}\pmb{\varepsilon}^{T}(\tilde{\pmb{B}}%
(\pmb{\beta}_{0})-\pmb{B}(\pmb{\beta}_{0}))\Vert_{\infty}\Vert(K_{n}%
\Gamma(\pmb{\beta}_{0},\pmb{\beta}_{0}))^{-1}\Vert_{\infty}|H_{r}%
(\pmb{\beta}_{0},\pmb{\beta}_{0})\pmb{b}_{0}|_{\infty}\\
=K_{n}o_{p}(1)O_{p}(1)O_{p}(h_{0})=o_{p}(1).
\end{array}
\]
This completes the proof of Lemma A.10.

\textbf{Proof of Theorem\ 3.2.} \ From Lemma A.4 and Assumption 8, we have
\[
\ddot{G}_{n}(\pmb{\beta}_{-d}^{\ast},\tilde{\pmb{b}}(\pmb{\beta}_{-d}^{\ast
}))=2\Omega(\pmb{\beta}_{-d}^{\ast})+o_{p}(1)=2\Omega(\pmb{\beta}_{0,-d}%
)+o_{p}(1)=2\Omega_{0}+o_{p}(1).\tag{A.55}
\]
Note that $(\tilde{\pmb{Y}}-\tilde{\pmb{W}}\pmb{\alpha}_{0}-\tilde
{\pmb{B}}(\pmb{\beta}_{0})\tilde{\pmb{b}}(\pmb{\alpha}_{0},\pmb{\beta}_{0}%
))^{T}\dot{\tilde{\pmb{B}}}_{r}(\pmb{\beta}_{0})\tilde{\pmb{b}}(\pmb{\alpha}
_{0},\pmb{\beta}_{0})$ can be written as
\[%
\begin{array}
[c]{l}%
(\tilde{\pmb{Y}}-\tilde{\pmb{W}}\pmb{\alpha}_{0}-\tilde{\pmb{B}}%
(\pmb{\beta}_{0})\tilde{\pmb{b}}(\pmb{\alpha}_{0},\pmb{\beta}_{0}))^{T}%
\dot{\tilde{\pmb{B}}}_{r}(\pmb{\beta}_{0})\tilde{\pmb{b}}(\pmb{\alpha}_{0}%
,\pmb{\beta} _{0})\\
=(\tilde{\pmb{Y}}-\tilde{\pmb{W}}\pmb{\alpha}_{0}-\tilde{\pmb{B}}%
(\pmb{\beta}_{0})\pmb{b}_{0})^{T}\dot{\tilde{\pmb{B}}}_{r}(\pmb{\beta}_{0}%
)\pmb{b}_{0}-(\tilde{\pmb{b}}(\pmb{\alpha}_{0},\pmb{\beta}_{0})-\pmb{b}_{0}%
)^{T}\tilde{\pmb{B}}^{T}(\pmb{\beta}_{0})\dot{\tilde{\pmb{B}}}_{r}%
(\pmb{\beta}_{0})\pmb{b}_{0}\\
+(\tilde{\pmb{Y}}-\tilde{\pmb{W}}\pmb{\alpha}_{0}-\tilde{\pmb{B}}%
(\pmb{\beta}_{0})\pmb{b}_{0})^{T}\dot{\tilde{\pmb{B}}}_{r}(\pmb{\beta}_{0}%
)(\tilde{\pmb{b}}(\pmb{\alpha}_{0},\pmb{\beta}_{0})-\pmb{b}_{0})\\
-(\tilde{\pmb{b}}(\pmb{\alpha}_{0},\pmb{\beta}_{0})-\pmb{b}_{0})^{T}%
\tilde{\pmb{B}}^{T}(\pmb{\beta}_{0})\dot{\tilde{\pmb{B}}}_{r}(\pmb{\beta}_{0}%
)(\tilde{\pmb{b}}(\pmb{\alpha}_{0},\pmb{\beta}_{0})-\pmb{b}_{0}).
\end{array}
\tag{A.56}
\]
Similar to Lemma A.9, we have $\Vert n^{-\frac{1}{2}}%
(\tilde{\pmb{Y}}-\tilde{\pmb{W}}\pmb{\alpha}_{0}-\tilde{\pmb{B}}%
(\pmb{\beta}_{0})\pmb{b}_{0})^{T}\dot{\tilde{\pmb{B}}}_{r}(\pmb{\beta}_{0}%
)\Vert_{\infty}=O_{p}(K_{n}^{3/2})$ and $|n^{-\frac{1}{2}}(\tilde
{\pmb{Y}}-\tilde{\pmb{W}}\pmb{\alpha}_{0}-\tilde{\pmb{B}}(\pmb{\beta}_{0}%
)\pmb{b}_{0})^{T}\tilde{\pmb{B}}(\pmb{\beta}_{0})|_{\infty}=O_{p}(h_{0}%
^{1/2})$. Hence
\[%
\begin{array}
[c]{ll}%
|\tilde{\pmb{b}}(\pmb{\alpha}_{0},\pmb{\beta}_{0})-\pmb{b}_{0}|_{\infty} &
\leq\frac{K_{n}}{n}\Vert(\frac{K_{n}}{n}\tilde{\pmb{B}}^{T}(\pmb{\beta}_{0}%
)\tilde{\pmb{B}}(\pmb{\beta}_{0}))^{-1}\Vert_{\infty}\\|\tilde{\pmb{B}}%
^{T}(\pmb{\beta}_{0})(\tilde{\pmb{Y}}-\tilde{\pmb{W}}\pmb{\alpha}_{0}%
-\tilde{\pmb{B}}(\pmb{\beta}_{0})\pmb{b}_{0})|_{\infty}\\
& =O_{p}(n^{-\frac{1}{2}}h_{0}^{-\frac{1}{2}}).
\end{array}
\]
Further, we deduce that
\[%
\begin{array}
[c]{l}%
n^{-\frac{1}{2}}|(\tilde{\pmb{Y}}-\tilde{\pmb{W}}\pmb{\alpha}_{0}%
-\tilde{\pmb{B}}(\pmb{\beta}_{0})\pmb{b}_{0})^{T}\dot{\tilde{\pmb{B}}}%
_{r}(\pmb{\beta}_{0})(\tilde{\pmb{b}}(\pmb{\alpha}_{0},\pmb{\beta}_{0}%
)-\pmb{b}_{0})|\\
\leq\Vert n^{-\frac{1}{2}}(\tilde{\pmb{Y}}-\tilde{\pmb{W}}\pmb{\alpha}_{0}%
-\tilde{\pmb{B}}(\pmb{\beta} _{0})\pmb{b}_{0})^{T}\dot{\tilde{\pmb{B}}}%
_{r}(\pmb{\beta}_{0})\Vert_{\infty}|\tilde{\pmb{b}}(\pmb{\alpha}_{0}%
,\pmb{\beta}_{0})-\pmb{b}_{0}|_{\infty}\\
=O_{p}(n^{-1/2}h_{0}^{-2})=o_{p}(1).
\end{array}
\tag{A.57}
\]
Applying (A.34), we have
\[%
\begin{array}
[c]{l}%
n^{-\frac{1}{2}}|(\tilde{\pmb{b}}(\pmb{\alpha}_{0},\pmb{\beta}_{0}%
)-\pmb{b}_{0})^{T}\tilde{\pmb{B}}^{T}(\pmb{\beta}_{0})\dot{\tilde{\pmb{B}}%
}_{r}(\pmb{\beta}_{0})(\tilde{\pmb{b}}(\pmb{\alpha}_{0},\pmb{\beta}
_{0})-\pmb{b}_{0})|\\
\leq n^{\frac{1}{2}}K_{n}\Vert\frac{1}{n}\tilde{\pmb{B}}^{T}(\pmb{\beta}_{0}%
)\dot{\tilde{\pmb{B}}}_{r}(\pmb{\beta}_{0})\Vert_{\infty}|\tilde
{\pmb{b}}(\pmb{\alpha}_{0},\pmb{\beta}_{0})-\pmb{b}_{0}|_{\infty}^{2}%
\\=O_{p}(n^{-1/2}h_{0}^{-2})=o_{p}(1).
\end{array}
\tag{A.58}
\]
Now $(\tilde{\pmb{Y}}-\tilde{\pmb{W}}\hat{\pmb{\alpha}}-\tilde{\pmb{B}}%
(\hat{\pmb{\beta}}_{-d})\hat{\pmb{b}})^{T}\tilde{\pmb{W}}_{k}$ can be written
as
\[
\begin{array}
[c]{l}%
(\tilde{\pmb{Y}}-\tilde{\pmb{W}}\hat{\pmb{\alpha}}-\tilde{\pmb{B}}%
(\hat{\pmb{\beta}}_{-d})\hat{\pmb{b}})^{T}\tilde{\pmb{W}}_{k}\\=(\tilde
{\pmb{Y}}-\tilde{\pmb{W}}\hat{\pmb{\alpha}}-\tilde{\pmb{B}}(\pmb{\beta}_{0}%
)\pmb{b}_{0})^{T}\tilde{\pmb{W}}_{k}-(\tilde{\pmb{b}}(\pmb{\alpha}_{0}%
,\pmb{\beta}_{0})-\pmb{b}_{0})^{T}\tilde{\pmb{B}}^{T}(\pmb{\beta}_{0}%
)\tilde{\pmb{W}}_{k}%
\end{array}
\]
for $k=1,\ldots,q$. Similar to the proof of Lemma A.9, we deduce that
\[
n^{-\frac{1}{2}}(\tilde{\pmb{Y}}-\tilde{\pmb{W}}\hat{\pmb{\alpha}}%
-\tilde{\pmb{B}}(\pmb{\beta}_{0})\pmb{b}_{0})^{T}\tilde{\pmb{W}}_{k}%
=n^{-\frac{1}{2}}\pmb{\varepsilon}^{T}\tilde{\pmb{W}}_{k}+o_{p}(1).
\]
We decompose $\pmb{\varepsilon}^{T}\tilde{\pmb{W}}_{k}$ into three terms as
\[%
\begin{array}
[c]{ll}%
\pmb{\varepsilon}^{T}\tilde{\pmb{W}}_{k} & =\sum_{i=1}^{n}\varepsilon
_{i}\Big(W_{ik}-\sum_{j=1}^{m}\frac{E(W_{lk}\xi_{j})}{\lambda_{j}}\xi
_{ij}\Big)-\\&\sum_{i=1}^{n}\varepsilon_{i}\sum_{j=1}^{m}\frac{\xi_{ij}}%
{\lambda_{j}}\Big(\frac{1}{n}\sum_{l=1}^{n}W_{lk}\xi_{lj}-E(W_{lk}\xi
_{j})\Big)\\
& -\sum_{i=1}^{n}\varepsilon_{i}\frac{1}{n}\sum_{l=1}^{n}W_{lk}(\tilde{\xi
}_{il}-\check{\xi}_{ij}).
\end{array}
\]
Similar to the proof of Lemma A.8, we have $\sum_{i=1}^{n}\varepsilon_{i}%
\frac{1}{n}\sum_{l=1}^{n}W_{lk}(\tilde{\xi}_{il}-\check{\xi}_{ij})=o_{p}(n)$.
Since
\[
\sum_{i=1}^{n}\varepsilon_{i}\Big(W_{ik}-\sum_{j=1}^{m}\frac{E(W_{lk}\xi_{j}%
)}{\lambda_{j}}\xi_{ij}\Big)=\sum_{i=1}^{n}\varepsilon_{i}V_{i}+\sum_{i=1}%
^{n}\varepsilon_{i}\sum_{j=m+1}^{\infty}w_{kj}\xi_{ij},
\]
$\sum_{i=1}^{n}\varepsilon_{i}\sum_{j=1}^{m}\frac{\xi_{ij}}{\lambda_{j}%
}\Big(\frac{1}{n}\sum_{l=1}^{n}W_{lk}\xi_{lj}-E(W_{lk}\xi_{j})\Big)=o_{p}(n)$
and $\sum_{i=1}^{n}\varepsilon_{i}\sum_{j=m+1}^{\infty}w_{kj}\xi_{ij}%
$ $=o_{p}(n)$, it follows that $n^{-\frac{1}{2}}\pmb{\varepsilon} ^{T}%
\tilde{\pmb{W}}_{k}=n^{-\frac{1}{2}}\pmb{\varepsilon}^{T}\pmb{V}_{k}+o_{p}%
(1)$, where $\pmb{V}_{k}=(V_{1k},\ldots,V_{nk})^{T}$. Similar to the proof of
Lemma A.10, we have
\[
\begin{array}
[c]{l}%
n^{-\frac{1}{2}}(\tilde{\pmb{b}}(\pmb{\alpha}_{0},\pmb{\beta}_{0}%
)-\pmb{b}_{0})^{T}\tilde{\pmb{B}}^{T}(\pmb{\beta}_{0})\tilde{\pmb{W}}%
_{k}\\=n^{-\frac{1}{2}}\pmb{\varepsilon}^{T}\pmb{B}(\pmb{\beta} _{0})\Gamma
^{-1}(\pmb{\beta}_{0},\pmb{\beta}_{0})E(\pmb{B}(Z^{T}\pmb{\beta}_{0}%
)W_{k})+o_{p}(1).
\end{array}
\]
Hence
\[
\begin{array}
[c]{l}%
n^{-\frac{1}{2}}(\tilde{\pmb{Y}}-\tilde{\pmb{W}}\hat{\pmb{\alpha}}%
-\tilde{\pmb{B}}(\hat{\pmb{\beta}}_{-d})\hat{\pmb{b}})^{T}\tilde{\pmb{W}}%
_{k}\\=n^{-\frac{1}{2}}\pmb{\varepsilon}^{T}(\pmb{V}_{k}-\pmb{B}(\pmb{\beta}_{0}%
)\Gamma^{-1}(\pmb{\beta}_{0},\pmb{\beta}_{0})E(\pmb{B}(Z^{T}\pmb{\beta}
_{0})W_{k}))+o_{p}(1).\tag{A.59}
\end{array}
\]
Now (3.9) follows from (A.55)-(A.59), Lemmas A.9 and A.10, and the Central
Limit Theorem. This completes the proof of Theorem 3.2.

\textbf{Lemma\ A.11.} \ Under the assumptions of Theorem 3.3, it holds that
\[
\|\tilde{\pmb{b}}(\hat{\pmb{\alpha}},\hat{\pmb{\beta}})-\pmb{b}_{0}%
\|^{2}=O_{p}(n^{-1}K_{n}^{2}).
\]

\textbf{Proof.} \ From Assumption 6 and Lemma A.3, all the eigenvalues of
$(\frac{K_{n}}{n}\tilde{\pmb{B}}^{T}(\hat{\pmb{\beta}})\tilde{\pmb{B}}%
(\hat{\pmb{\beta}}))^{-1}$ are bounded away from zero and infinity, except
possibly on an event whose probability tends to zero. We then have
\[%
\begin{array}
[c]{ll}%
\Vert\tilde{\pmb{b}}(\hat{\pmb{\alpha}},\hat{\pmb{\beta}})-\pmb{b}_{0}%
\Vert^{2}\leq CK_{n}^{2}\Vert\tilde{\pmb{B}}^{T}(\hat{\pmb{\beta}}%
)(\tilde{\pmb{Y}}-\tilde{\pmb{W}}\hat{\pmb{\alpha}}-\tilde{\pmb{B}}%
(\hat{\pmb{\beta}})\pmb{b}_{0})\Vert^{2}/n^{2}, &
\end{array}
\tag{A.60}
\]
where $\Vert a\Vert=(a_{1}^{2}+\ldots+a_{k}^{2})^{1/2}$ for a vector
$a=(a_{1},\ldots,a_{k})^{T}$. Let $\pmb{F}(\pmb{\alpha},\pmb{\beta}_{-d}%
)=\tilde{\pmb{B}}^{T}(\pmb{\beta}_{-d})(\tilde{\pmb{Y}}-\tilde{\pmb{W}}%
\pmb{\alpha}-\tilde{\pmb{B}}(\pmb{\beta}_{-d})\pmb{b}_{0})$. By a Taylor
expansion, we have that
\[
\pmb{F}(\hat{\pmb{\alpha}},\hat{\pmb{\beta}}_{-d})=\pmb{F}(\pmb{\alpha}_{0}%
,\pmb{\beta}_{0,-d})-\tilde{\pmb{B}}^{T}(\pmb{\beta}_{-d}^{\star}%
)\tilde{\pmb{W}}(\hat{\pmb{\alpha}}-\pmb{\alpha}_{0})+\frac{\partial
\pmb{F}}{\partial\pmb{\beta}_{-d}}\Big
|_{\pmb{\beta}_{-d}=\pmb{\beta}_{-d}^{\star}}(\hat{\pmb{\beta}}_{-d}%
-\pmb{\beta}_{0,-d}),\tag{A.61}
\]
where $({\pmb{\alpha}^{\star}}^{T},{\pmb{\beta}_{-d}^{\star}}^{T})^{T}$ is
between $({\hat{\pmb{\alpha}}}^{T},{\hat{\pmb{\beta}}_{-d}}^{T})^{T}$ and
$(\pmb{\alpha}_{0}^{T},\pmb{\beta}_{0,-d}^{T})^{T},$ and
\[
\frac{\partial\pmb{F}}{\partial\beta_{r}}\Big|_{\pmb{\beta}_{-d}%
=\pmb{\beta}_{-d}^{\star}}=\dot{\tilde{\pmb{B}}}_{r}^{T}(\pmb{\beta}_{-d}%
^{\star})(\tilde{\pmb{Y}}-\tilde{\pmb{W}}\pmb{\alpha}^{\star}-\tilde
{\pmb{B}}(\pmb{\beta}_{-d}^{\star})\pmb{b}_{0})-\tilde{\pmb{B}}^{T}%
(\pmb{\beta}_{-d}^{\star})\dot{\tilde{\pmb{B}}}_{r}(\pmb{\beta}_{-d}^{\star
})\pmb{b}_{0}.
\]
Similar to the proof of (A.33) and (A.36), we obtain that
\[
\begin{array}
[c]{l}%
\Vert\frac{1}{n}\dot{\tilde{\pmb{B}}}_{r}^{T}(\pmb{\beta}_{-d}^{\star}%
)(\tilde{\pmb{Y}}-\tilde{\pmb{W}}\pmb{\alpha}^{\star}-\tilde{\pmb{B}}%
(\pmb{\beta}_{-d}^{\star})\pmb{b}_{0})\Vert^{2}\\=\Vert E[\dot{\pmb{B}}%
_{r}(Z^{T}\pmb{\beta}_{-d}^{\star})V^{T}](\pmb{\alpha}^{\star}%
-\pmb{\alpha}_{0})\Vert^{2}+o_{p}(1)=o_{p}(1)
\end{array}
\]
and $\Vert\frac{1}{n}\tilde{\pmb{B}}^{T}(\pmb{\beta}^{\star})\dot
{\tilde{\pmb{B}}}_{r}(\pmb{\beta}^{\star})\pmb{b}_{0}\Vert^{2}=O_{p}(1)$. From
Theorem 3.2, it holds that $\Vert\hat{\pmb{\beta}}_{-d}$ $-\pmb{\beta}_{0,-d}%
\Vert^{2}=O_{p}(n^{-1})$. Hence
\[
\Big\|\frac{\partial\pmb{F}}{\partial\pmb{\beta}_{-d}}\Big|_{\pmb{\beta}_{-d}%
=\pmb{\beta}_{-d}^{\star}}(\hat{\pmb{\beta}}_{-d}-\pmb{\beta}_{0,d}%
)\Big\|^{2}\leq\sum_{r=1}^{d-1}\Big\|\frac{\partial\pmb{F}}{\partial\beta_{r}%
}\Big|_{\pmb{\beta}_{-d}=\pmb{\beta}_{-d}^{\star}}\Big\|^{2}\Vert
\hat{\pmb{\beta}}-\pmb{\beta}_{0}\Vert^{2}=o_{p}(n).\tag{A.62}
\]
It is easy to prove that $\Vert\tilde{\pmb{B}}^{T}(\pmb{\beta}_{-d}^{\star
})\tilde{\pmb{W}}(\hat{\pmb{\alpha}}-\pmb{\alpha}_{0})\Vert^{2}=o_{p}(n)$. By
arguments similar to those used to prove Lemma A.9, we can prove that
$\Vert\pmb{F}(\pmb{\alpha}_{0},\pmb{\beta}_{0,-d})\Vert^{2}=O_{p}(n)$. Now
Lemma A.11 follows from (A.60)-(A.62). This completes the proof of Lemma A.11.

\textbf{Lemma\ A.12.} \ Define $\check{a}_{j}=\frac{1}{\hat{\lambda}_{j}%
}E[(Y-W^{T}\pmb{\alpha}_{0}-g(Z^{T}\pmb{\beta}_{0}))\xi_{j}]$. Under the
assumptions of Theorem 3.3, it holds that
\[
\sum_{j=1}^{\tilde{m}}(\hat{a}_{j}-\check{a}_{j})^{2}=O_{p}(n^{-1}\tilde
{m}\lambda_{\tilde{m}}^{-1} +n^{-2}\tilde{m}\lambda_{\tilde{m}}^{-2}\sum
_{j=1}^{\tilde{m}}a_{j}^{2}\lambda_{j}^{-2}j^{3}).
\]

\textbf{Proof.} \ Note that $E[(Y-W^{T}\pmb{\alpha}_{0}-g(Z^{T}\pmb{\beta}_{0}%
))\xi_{j}]=a_{j}\lambda_{j}$. Define $I_{1}=\frac{1}{n}\sum_{i=1}^{n}%
[Y_{i}-W_{i}^{T}\pmb{\alpha}_{0}-g(Z_{i}^{T}\pmb{\beta}_{0})]\xi_{ij}%
-a_{j}\lambda_{j}$, $I_{2}=\frac{1}{n}\sum_{i=1}^{n}[Y_{i}-W_{i}%
^{T}\pmb{\alpha}_{0}-g(Z_{i}^{T}\pmb{\beta}_{0})] $ $(\hat{\xi}_{ij}-\xi_{ij})$
and $I_{3}=\frac{1}{n}\sum_{i=1}^{n}[W_{i}^{T}(\hat{\pmb{\alpha}}%
-\pmb{\alpha}_{0})+(\hat{g}(Z_{i}^{T}\hat{\pmb{\beta}})-g(Z_{i}^{T}%
\pmb{\beta}_{0}))]\hat{\xi}_{ij}$. Then we have
\[
\sum_{j=1}^{\tilde{m}}(\hat{a}_{j}-\check{a}_{j})^{2}\leq3\sum_{j=1}%
^{\tilde{m}}\lambda_{j}^{-2}(I_{1}^{2}+I_{2}^{2}+I_{3}^{2})[1+o_{p}%
(1)],\tag{A.63}
\]
where $o_{p}(1)$ holds uniformly for $j=1,\ldots,\tilde{m}$. Since
$E(I_{1})=0$ and $E(I_{1}^{2})\leq\frac{1}{n}[\sum_{k=1}^{\infty}a_{k}%
^{2}E(\xi_{k}^{2}\xi_{j}^{2})+\sigma^{2}\lambda_{j}]\leq C\lambda_{j}/n$, we
obtain that
\[
\sum_{j=1}^{\tilde{m}}\lambda_{j}^{-2}I_{1}^{2}=O_{p}(n^{-1}\sum_{j=1}%
^{\tilde{m}}\lambda_{j}^{-1})=O_{p}(n^{-1}\tilde{m}\lambda_{\tilde{m}}%
^{-1}).\tag{A.64}
\]
Let $M(t)=E[(Y_{i}-W_{i}^{T}\pmb{\alpha}_{0}-g(Z_{i}^{T}\pmb{\beta}_{0}%
))X_{i}(t)]=\sum_{k=1}^{\infty}a_{k}\lambda_{k}\phi_{k}(t)$. Then
\[%
\begin{array}
[c]{ll}%
I_{2}^{2} & \leq2\int_{\mathcal{T}}\Big(\frac{1}{n}\sum_{i=1}^{n}[Y_{i}%
-W_{i}^{T}\pmb{\alpha}_{0}-g(Z_{i}^{T}\pmb{\beta}_{0})]X_{i}(t)-M(t)\Big)^{2}%
dt\Vert\hat{\phi}_{j}-\phi_{j}\Vert^{2}\\
& +2\Big(\int_{\mathcal{T}}M(t)(\hat{\phi}_{j}(t)-\phi_{j}(t))dt\Big)^{2}.
\end{array}
\]
Applying Assumption 1, it holds that
\[%
\begin{array}
[c]{l}%
E\Big(\int_{\mathcal{T}}\Big(\frac{1}{n}\sum_{i=1}^{n}[Y_{i}-W_{i}%
^{T}\pmb{\alpha}_{0}-g(Z_{i}^{T}\pmb{\beta}_{0})]X_{i}(t)-M(t)\Big)^{2}%
dt\Big)\\
\leq\frac{1}{n}\int_{\mathcal{T}}E([Y_{i}-W_{i}^{T}\pmb{\alpha}_{0}%
-g(Z_{i}^{T}\pmb{\beta}_{0})]^{2}X_{i}^{2}(t))dt=O(n^{-1}).
\end{array}
\]
From (A.8), we obtain $\sum_{j=1}^{\tilde{m}}\lambda_{j}^{-2}\Vert\hat{\phi
}_{j}-\phi_{j}\Vert^{2}=O_{p}(n^{-1}\tilde{m}^{3}\lambda_{\tilde{m}}^{-2}%
\log\tilde{m})$. By arguments similar to those used in the proof of (5.15) of
Hall and Horowitz \cite{r11},
it follows that
\[
\begin{array}
[c]{l}%
\sum_{j=1}^{\tilde{m}}\lambda_{j}^{-2}\Big(\int_{\mathcal{T}}M(t)(\hat{\phi
}_{j}(t)-\phi_{j}(t))dt\Big)^{2}\\=O_{p}(\frac{\tilde{m}}{n\lambda_{\tilde{m}}%
}+\frac{\tilde{m}}{n^{2}\lambda_{\tilde{m}}^{2}}\sum_{j=1}^{\tilde{m}}%
a_{j}^{2}\lambda_{j}^{-2}j^{3}+\frac{\tilde{m}^{3}\log\tilde{m}}{n^{2}%
\lambda_{\tilde{m}}^{2}}).
\end{array}
\]
Hence, using the assumption that $n^{-1/2}\tilde{m}\lambda_{\tilde{m}%
}\rightarrow0$, we obtain
\[
\sum_{j=1}^{\tilde{m}}\lambda_{j}^{-2}I_{2}^{2}=O_{p}(n^{-1}\tilde{m}%
\lambda_{\tilde{m}}^{-1}+n^{-2}\tilde{m}\lambda_{\tilde{m}}^{-2}\sum
_{j=1}^{\tilde{m}}a_{j}^{2}\lambda_{j}^{-2}j^{3}).\tag{A.65}
\]
Define $I_{31}=\frac{1}{n}\sum_{i=1}^{n}[\hat{g}(Z_{i}^{T}\hat{\pmb{\beta}}%
)-g_{0}(Z_{i}^{T}\hat{\pmb{\beta}})]\hat{\xi}_{ij}$, $I_{32}=\frac{1}{n}%
\sum_{i=1}^{n}[W_{i}^{T}(\hat{\pmb{\alpha}}-\pmb{\alpha}_{0})+(g_{0}(Z_{i}%
^{T}\hat{\pmb{\beta} })-g(Z_{i}^{T}\pmb{\beta}_{0}))]\hat{\xi}_{ij}$, 
$L_{j}=(l_{jkk^{\prime}})_{K_{n}\times K_{n}}$ with $l_{jkk^{\prime}}%
=(\frac{1}{n}\sum_{i=1}^{n}B_{k}(Z_{i}^{T}\hat{\pmb{\beta}})\hat{\xi}%
_{ij}) $ $(\frac{1}{n}\sum_{i=1}^{n}B_{k^{\prime}}(Z_{i}^{T}\hat{\pmb{\beta}}%
)\hat{\xi}_{ij})$. We write $\frac{1}{n}\sum_{i=1}^{n}B_{k}(Z_{i}^{T}%
\hat{\pmb{\beta}})\hat{\xi}_{ij}=\frac{1}{n}\sum_{i=1}^{n}[B_{k}(Z_{i}%
^{T}\pmb{\beta}_{0})\xi_{ij}+(B_{k}(Z_{i}^{T}\hat{\pmb{\beta} })-B_{k}%
(Z_{i}^{T}\pmb{\beta}_{0}))\xi_{ij}+B_{k}(Z_{i}^{T}\hat{\pmb{\beta}})(\hat
{\xi}_{ij}-\xi_{ij})]$. Then we have
\[
\begin{array}
[c]{l}%
|L_{j}|_{\infty}=\max_{k,k^{\prime}}|l_{jkk^{\prime}}|\\\leq\sum_{k=1}^{K_{n}%
}\Big(\frac{1}{n}\sum_{i=1}^{n}B_{k}(Z_{i}^{T}\pmb{\beta}_{0})\xi
_{ij}\Big)^{2}+\frac{C}{n}\sum_{i=1}^{n}[h_{0}^{-2}\Vert\hat{\pmb{\beta}}%
-\pmb{\beta}_{0}\Vert^{2}\xi_{ij}^{2}+(\hat{\xi}_{ij}-\xi_{ij})^{2}].
\end{array}
\]
Simple calculations yield $\sum_{k=1}^{K_{n}}E\Big(\frac{1}{n}\sum_{i=1}%
^{n}B_{k}(Z_{i}^{T}\pmb{\beta}_{0})\xi_{ij}\Big)^{2}\leq Cn^{-1}\lambda_{j}$.
Applying Lemma A.11, we obtain that $\Vert\tilde{\pmb{b}}(\hat{\pmb{\alpha}}%
,\hat{\pmb{\beta} })-\pmb{b}_{0}\Vert_{\infty}^{2}\leq K_{n}\Vert
\tilde{\pmb{b}}(\hat{\pmb{\alpha}},\hat{\pmb{\beta} })-\pmb{b}_{0}\Vert
^{2}=O_{p}(n^{-1}K_{n}^{3})$. Hence, under the assumptions of Theorem 3.3, it
holds that
\[%
\begin{array}
[c]{ll}%
\sum_{j=1}^{\tilde{m}}\lambda_{j}^{-2}I_{31}^{2} & \leq\sum_{j=1}^{\tilde{m}%
}\lambda_{j}^{-2}\Vert\tilde{\pmb{b}}(\hat{\pmb{\alpha}},\hat{\pmb{\beta}}%
)-\pmb{b}_{0}\Vert_{\infty}^{2}\cdot|L_{j}|_{\infty}\\
& =O_{p}(n^{-2}\tilde{m}\lambda_{\tilde{m}}^{-1}h_{0}^{-3}+n^{-2}\tilde
{m}\lambda_{\tilde{m}}^{-1}h_{0}^{-5}+n^{-2}\tilde{m}^{3}\lambda_{\tilde{m}%
}^{-2}h_{0}^{-3}\log\tilde{m})\\
& =O_{p}(n^{-1}\tilde{m}\lambda_{\tilde{m}}^{-1}).
\end{array}
\tag{A.66}
\]
Using a Taylor expansion, Theorem 3.2, and the assumption that $nh_{0}%
^{2p}\rightarrow0$, we deduce that
\[%
\begin{array}
[c]{ll}%
\sum_{j=1}^{\tilde{m}}\lambda_{j}^{-2}I_{32}^{2} & \leq\Big(\sum_{j=1}%
^{\tilde{m}}\frac{1}{n\lambda_{j}^{2}}\sum_{i=1}^{n}\hat{\xi}_{ij}%
^{2}\Big)\\&\Big(\frac{1}{n}\sum_{i=1}^{n}[W_{i}^{T}(\hat{\pmb{\alpha}}%
-\pmb{\alpha} _{0})+(g_{0}(Z_{i}^{T}\hat{\pmb{\beta}})-g(Z_{i}^{T}%
\pmb{\beta}_{0}))]^{2}\Big)\\
& =O_{p}(\tilde{m}\lambda_{\tilde{m}}^{-1}+n^{-1}\tilde{m}^{3}\lambda
_{\tilde{m}}^{-2}\log\tilde{m})O_{p}(n^{-1}+h_{0}^{2p})\\&=O_{p}(n^{-1}\tilde
{m}\lambda_{\tilde{m}}^{-1}).
\end{array}
\tag{A.67}
\]
Now Lemma A.12 follows from combining (A.63)-(A.67).

\textbf{Proof of Theorem\ 3.3.} \ Note that
\[
\begin{array}
[c]{lll}%
\int_{\mathcal{T}}[\hat{a}(t)-a(t)]^{2}dt&\leq C\Big(\sum_{j=1}^{\tilde{m}%
}(\hat{a}_{j}-\check{a}_{j})^{2}+\sum_{j=1}^{\tilde{m}}(\check{a}_{j}%
-a_{j})^{2}+\\&\tilde{m}\sum_{j=1}^{\tilde{m}}a_{j}^{2}\Vert\hat{\phi}_{j}%
-\phi_{j}\Vert^{2}+\sum_{j=\tilde{m}+1}^{\infty}a_{j}^{2}\Big)\tag{A.68}
\end{array}
\]
and
\[%
\begin{array}
[c]{ll}%
\sum_{j=1}^{\tilde{m}}(\check{a}_{j}-a_{j})^{2}=\sum_{j=1}^{\tilde{m}}%
\frac{(\hat{\lambda}_{j}-\lambda_{j})^{2}}{\lambda_{j}^{2}}a_{j}^{2}%
[1+o_{p}(1)]=O_{p}(n^{-1}\lambda_{\tilde{m}}^{-1}\sum_{j=1}^{\tilde{m}}%
a_{j}^{2}\lambda_{j}^{-1}). &
\end{array}
\tag{A.69}
\]
Assumption 3 implies that $\tilde{m}\sum_{j=1}^{\tilde{m}}a_{j}^{2}\Vert
\hat{\phi}_{j}-\phi_{j}\Vert^{2}=O_{p}(\tilde{m}n^{-1}\sum_{j=1}^{\tilde{m}%
}a_{j}^{2}j^{2}\log j)$ $=o_{p}(\tilde{m}/n)$ and $\sum_{j=\tilde{m}+1}^{\infty
}a_{j}^{2}=O(\tilde{m}^{-2\gamma+1})$. Now (3.10) follows from Lemma A.12,
(A.68) and (A.69). This completes the proof of Theorem 3.3.

\textbf{Proof of Theorem\ 3.4.} \ From Assumption 6 and Lemma A.3, all the
eigenvalues of $(\frac{K^{\ast}_{n}}{n}\pmb{B}^{\ast T}(\hat{\pmb{\beta}
})\pmb{B}^{\ast}(\hat{\pmb{\beta}}))^{-1}$ are bounded away from zero and
infinity, except possibly on an event whose probability tends to zero. Similar
to (3.1), there exists a spline function $g^{\ast}(u)=\sum_{k=1}^{K^{\ast}_{n}%
}b_{0k}^{\ast}B_{k}^{\ast}(u)$ such that
\[
\sup_{u\in\lbrack U_{\pmb{\beta}_{0}},U^{\pmb{\beta}_{0}}]}|g(u)-g^{\ast}(u)|\leq Ch^{p}.\tag{A.70}
\]
Let $\pmb{b}^{\ast}_{0}=(b^{\ast}_{01},\ldots,b^{\ast}_{0K^{\ast}_{n}})^{T}$.
Using the properties of B-splines (de Boor 1978), we obtain
\[%
\begin{array}
[c]{ll}%
\int_{U_{\pmb{\beta}_{0}}}^{U^{\pmb{\beta}_{0}}}(\hat{g}(u)-g(u))^{2}du & \leq
C(\Vert\pmb{b}^{\ast}(\hat{\pmb{\alpha}},\hat{\pmb{\beta}})-\pmb{b}^{\ast}%
_{0}\Vert^{2}/K^{\ast}_{n}+h^{2p}).\\
&
\end{array}
\tag{A.71}
\]
Using arguments similar to those used to prove Lemma A.11 and using the fact
that
\[
\sum_{k=1}^{K^{\ast}_{n}}(\sum_{i=1}^{n}B^{\ast}_{k}(Z_{i}\pmb{\beta} _{0}%
)R^{\ast}(Z_{i}\pmb{\beta}_{0}))^{2}=O_{p}(n^{2}h^{2p+1}),
\]
where $R^{\ast}(u)=g(u)-\sum_{k=1}^{K^{\ast}_{n}}b^{\ast}_{0k}B^{\ast}_{k}(u)$,
one can prove that
\[
\Vert\pmb{b}^{\ast}(\hat{\pmb{\alpha}},\hat{\pmb{\beta}})-\pmb{b}^{\ast}%
_{0}\Vert^{2}=O_{p}(n^{-1}{K^{\ast}_{n}}^{2})+O_{p}(h^{2p-1}).\tag{A.72}
\]
Now (3.12) follows from (A.71) and the fact that $h=O({K^{\ast}_{n}}^{-1})$. This
completes the proof of Theorem 3.4.

\textbf{Proof of Theorem\ 3.5.} \ Observe that
\[%
\begin{array}
[c]{l}%
\mbox{MSPE}\leq3\{\Vert\hat{a}-a\Vert_{K}^{2}+(\hat{\pmb{\alpha}}%
-\pmb{\alpha}_{0})^{T}E(WW^{T})(\hat{\pmb{\alpha}}-\pmb{\alpha}_{0}%
)+\\E([\hat{g}(Z_{n+1}^{T}\hat{\pmb{\beta} })-g(Z_{n+1}^{T}\pmb{\beta}_{0}%
)]^{2}|\mathcal{S})\},
\end{array}
\tag{A.73}
\]
where $\Vert\hat{a}-a\Vert_{K}^{2}=\int_{\mathcal{T}}\int_{\mathcal{T}%
}K(s,t)[\hat{a}(s)-a(s)][\hat{a}(t)-a(t)]dsdt$. Under the assumptions of
Theorem 3.5, using arguments similar to those used in the proof of Theorem 2 of
Tang (2015), we deduce that
\[
\Vert\hat{a}-a\Vert_{K}^{2}=O_{p}(n^{-(\delta+2\gamma-1)/(\delta+2\gamma
)}).\tag{A.74}
\]
Write
\[
\hat{g}(Z_{n+1}^{T}\hat{\pmb{\beta}})-g(Z_{n+1}^{T}\pmb{\beta}_{0})=\hat
{g}(Z_{n+1}^{T}\hat{\pmb{\beta}})-g^{\ast}(Z_{n+1}^{T}\hat{\pmb{\beta}}%
)+g^{\ast}(Z_{n+1}^{T}\hat{\pmb{\beta}})-g(Z_{n+1}^{T}\pmb{\beta}_{0}).
\]
Using a Taylor expansion, Theorems 3.2 and 3.4, (A.71), and the property of
B-spline function, we obtain
\[%
\begin{array}
[c]{ll}%
&E([\hat{g}(Z_{n+1}^{T}\hat{\pmb{\beta}})-g^{\ast}(Z_{n+1}^{T}\hat
{\pmb{\beta}})]^{2}|\mathcal{S}) \\& \leq2E([\hat{g}(Z_{n+1}^{T}\pmb{\beta}_{0}%
)-g^{\ast}(Z_{n+1}^{T}\pmb{\beta}_{0})]^{2}|\mathcal{S})\\
& +Ch^{-2}(\sum_{k=1}^{K^{\ast}_{n}}|\hat{b}_{k}-b^{\ast}_{0k}|)^{2}%
(\hat{\pmb{\beta}}-\pmb{\beta}_{0})^{T}E(ZZ^{T})(\hat{\pmb{\beta}}%
-\pmb{\beta}_{0})\\
& =O_{p}(n^{-\frac{2p}{2p+1}})+O_{p}(n^{-2}h^{-5}+n^{-1}h^{2p-4}%
)=O_{p}(n^{-\frac{2p}{2p+1}}).
\end{array}
\]
Using a Taylor expansion, Theorem 3.2 and (A.70), we also obtain
\[%
\begin{array}
[c]{ll}%
E([g^{\ast}(Z_{n+1}^{T}\hat{\pmb{\beta}})-g(Z_{n+1}^{T}\pmb{\beta}_{0}%
]^{2}|\mathcal{S}) & \leq2E([g^{\ast}(Z_{n+1}^{T}\pmb{\beta}_{0})-g(Z_{n+1}%
^{T}\pmb{\beta}_{0})]^{2}|\mathcal{S})\\
& +C(\hat{\pmb{\beta}}-\pmb{\beta}_{0})^{T}E(ZZ^{T})(\hat{\pmb{\beta}}%
-\pmb{\beta}_{0})=O_{p}(h^{2p}).
\end{array}
\]
Hence, $E([\hat{g}(Z_{n+1}^{T}\hat{\pmb{\beta}})-g(Z_{n+1}^{T}\pmb{\beta}_{0}%
)]^{2}|\mathcal{S})=O_{p}(n^{-2p/(2p+1)})$. Now (3.14) follows from (A.73),
(A.74) and Theorem 3.2. This completes the proof of Theorem 3.5.

\medskip


\section*{Acknowledgments} 


Part of data collection and sharing for this project was funded by the
Alzheimer's Disease Neuroimaging Initiative (ADNI) (National Institutes of
Health Grant U01 AG024904). ADNI is funded by the National Institute on Aging,
the National Institute of Biomedical Imaging and Bioengineering, and through
generous contributions from the following: Alzheimer's Association;
Alzheimer's Drug Discovery Foundation; BioClinica, Inc.; Biogen Idec Inc.;
Bristol-Myers Squibb Company; Eisai Inc.; Elan Pharmaceuticals, Inc.; Eli
Lilly and Company; F. Hoffmann-La Roche Ltd and its affilated company
Genentech, Inc.; GE Healthcare; Innogenetics, N.V.; IXICO Ltd.; Janssen
Alzheimer Immunotherapy Research \& Development, LLC.; Johnson \& Johnson
Pharmaceutical Research \& Development LLC.; Medpace, Inc.; Merck \& Co.,
Inc.; Meso Scale Diagnostics, LLC.; NeuroRx Research; Novartis Pharmaceuticals
Corporation; Pfizer Inc.; Piramal Imaging; Servier; Synarc Inc.; and Takeda
Pharmaceutical Company. The Canadian Institutes of Health Research is
providing funds to support ADNI clinical sites in Canada. Private sector
contributions are facilitated by the Foundation for the National Institutes of
Health (\url{www.fnih.org}). The grantee organization is the Northern
California Institute for Research and Education, and the study is coordinated
by the Alzheimer's Disease Cooperative Study at the University of California,
San Diego. ADNI data are disseminated by the Laboratory for Neuro Imaging at
the University of California, Los Angeles.



\end{document}